\pdfoutput=1
\RequirePackage{ifpdf}
\ifpdf 
\documentclass[pdftex]{sigma}
\else
\documentclass{sigma}
\fi

\usepackage[all]{xy}
\usepackage{mathrsfs}
\usepackage{enumitem}
\usepackage[all]{xy}
\usepackage{stmaryrd}
\usepackage{tensor}

\numberwithin{equation}{section}

\newtheorem{Theorem}{Theorem}[section]
\newtheorem{Corollary}[Theorem]{Corollary}
\newtheorem{Lemma}[Theorem]{Lemma}
\newtheorem{Proposition}[Theorem]{Proposition}
\newtheorem{Claim}[Theorem]{Claim}

 { \theoremstyle{definition}
\newtheorem{Definition}[Theorem]{Definition}
\newtheorem{Note}[Theorem]{Note}
\newtheorem{Example}[Theorem]{Example}
\newtheorem{Remark}[Theorem]{Remark} }

\newcommand{\rmap}{\to}
\newcommand{\X}{\ensuremath{\mathfrak{X}}}

\DeclareMathOperator{\Ann}{Ann}
\newcommand{\pr}{\operatorname{pr}}

\renewcommand{\hom}{\operatorname{Hom}}

\newcommand{\Ker}{\operatorname{Ker}} 
\renewcommand{\Im}{\operatorname{Im}} 
\newcommand{\rank}{\operatorname{rk}} 

\newcommand{\de}{\mathrm{d}} 

\newcommand{\JM}{(P,L,\{\cdot,\cdot\})} 
\newcommand{\CM}{(M,H)} 
\newcommand{\SM}{(S,\omega)} 
\newcommand{\PM}{(P,\Lambda)} 
\newcommand{\cotan}{\mathrm{T}^*} 
\newcommand{\HH}{\mathrm{H}}
\newcommand{\SYM}{L^* \setminus \{0\}} 

\newcommand{\pg}{\mathbb{P}(\mathfrak{g}^*)} 
\newcommand{\RR}{\mathbb{R}}
\newcommand{\ZZ}{\mathbb{Z}}

\begin{document}

\allowdisplaybreaks

\newcommand{\arXivNumber}{1406.2138}

\renewcommand{\thefootnote}{}

\renewcommand{\PaperNumber}{033}

\FirstPageHeading

\ShortArticleName{Contact Isotropic Realisations of Jacobi Manifolds via Spencer Operators}

\ArticleName{Contact Isotropic Realisations of Jacobi Manifolds\\ via Spencer Operators\footnote{This paper is a~contribution to the Special Issue ``Gone Fishing''. The full collection is available at \href{http://www.emis.de/journals/SIGMA/gone-fishing2016.html}{http://www.emis.de/journals/SIGMA/gone-fishing2016.html}}}

\Author{Mar\'ia Amelia SALAZAR~$^\dag$ and Daniele SEPE~$^\ddag$}

\AuthorNameForHeading{M.A.~Salazar and D.~Sepe}

\Address{$^\dag$~IMPA, Estrada Dona Castorina 110, Rio de Janeiro, Brazil 22460-320}
\EmailD{\href{mailto:mariasalazarp@gmail.com}{mariasalazarp@gmail.com}}
\URLaddressD{\url{http://w3.impa.br/~salazarp/}}

\Address{$^\ddag$~Universidade Federal Fluminense, Instituto de Matem\'atica, Departamento de Matem\'atica\\
 \hphantom{$^\ddag$}~Aplicada, Rua M\'ario Santos Braga S/N, Campus do Valonguinho, Niter\'oi, Brazil 24020-140}
\EmailD{\href{mailto:danielesepe@id.uff.br}{danielesepe@id.uff.br}}
\URLaddressD{\url{https://sites.google.com/site/danielesepe/}}

\ArticleDates{Received October 07, 2016, in f\/inal form May 17, 2017; Published online May 25, 2017}

\Abstract{Motivated by the importance of symplectic isotropic realisations in the study of Poisson manifolds, this paper investigates the local and global theory of contact isotropic realisations of Jacobi manifolds, which are those of minimal dimension. These arise naturally when considering multiplicity-free actions in contact geometry, as shown in this paper. The main results concern a classif\/ication of these realisations up to a~suitable notion of isomorphism, as well as establishing a relation between the existence of symplectic and contact isotropic realisations for Poisson manifolds. The main tool is the classical Spencer operator which is related to Jacobi structures via their associated Lie algebroid, which allows to generalise previous results as well as providing more conceptual proofs for existing ones.}


\Classification{53D10; 53D17; 53D20; 37J15}

{\small \tableofcontents}

\renewcommand{\thefootnote}{\arabic{footnote}}
\setcounter{footnote}{0}

\section{Introduction}\label{sec:introduction}
Poisson structures can be considered as inf\/initesimal objects by thinking of them as Lie algebroid structures on cotangent bundles together with the canonical symplectic form. One approach to understand the geometry of a Poisson manifold is to consider its `representations', which correspond to {\it symplectic realisations} (cf.~\cite{CrainicFernandes:jdd,Weinstein}). The representations of smallest dimension, which are known as {\it isotropic}, play an important role, for they enjoy special geometric properties, such as the existence of a {\it $\ZZ$-affine structure} transverse to the symplectic foliation (cf.~\cite{PMCT2, dd}). Moreover, they are related to multiplicity-free actions of compact Lie groups on symplectic manifolds (cf.~\cite{dd,gui_sja}), as well as to the recently introduced notion of `compactness' in Poisson geometry (cf.\ \cite{PMCT1,PMCT2}). However, regarding Poisson manifolds as being {\it Jacobi}, the associated Lie algebroid structure is def\/ined on the f\/irst jet bundle of the manifold, which is endowed with a canonical contact form. This point of view sheds a dif\/ferent light on Poisson geometry. For instance, in this case, the `representations' are {\it contact} realisations. Continuing with the above analogy, a natural question is to study geometric properties of those representations of smallest dimension and to investigate whether they are connected with Hamiltonian Lie group actions in contact geometry as well as a~slightly more general notion of `compactness' for Poisson and Jacobi manifolds. The aim of the present paper is to begin to address the above question, with a view to continue this line of research in future papers.

A {\it Jacobi structure} on a manifold $M$ is a real line bundle $L \to M$ whose space of sec\-tions~$\Gamma(L)$ is endowed with a local Lie bracket $\{\cdot, \cdot\}$ (cf.\ Def\/inition~\ref{defn:jacobi_structure}); these generalise Poisson structures and contact manifolds (cf.\ Section~\ref{sec:basic-notions}). A~Jacobi manifold $\JM$ is encoded in a Lie algebroid structure on the f\/irst jet bundle $J^1L$ which is compatible with the {\it classical Spencer operator} (cf.\ Note~\ref{obs:prop_cso}). {\it Contact realisations} of~$\JM$ are representations of the Lie algebroid $J^1L$: these are surjective submersions $\phi \colon \CM \to \JM$, where $\CM$ is a contact manifold and~$\phi$ is a Jacobi morphism satisfying a transversality condition (cf.\ Def\/inition~\ref{defn:CR}). The main object of study of this paper are contact {\it isotropic} realisations, which can be thought of as being of smallest dimension (cf.\ Def\/inition~\ref{defn:CIR} and Note~\ref{item:po}).

As motivation, it is shown that these objects arise naturally when considering contact ma\-ni\-folds and their `integrable' symmetries; some of the results may be of independent in\-te\-rest. For instance, contact isotropic realisations arise when considering the contact analogue of multiplicity-free Hamiltonian actions on symplectic manifolds (cf.\ Def\/inition~\ref{defn:mult_free}); these are the non-abelian analogues of contact toric manifolds studied in~\cite{lerman_contact_toric,wolbert}. The main results of this paper concern the classif\/ication of contact isotropic realisations to a suitable notion of isomorphism. These results are analogous to those of~\cite{dd} for symplectic isotropic realisations of Poisson manifolds and generalise~\cite{Banyaga,jovanovic}. In particular,
\begin{itemize}\itemsep=0pt
\item given a contact isotropic realisation of a Jacobi manifold, we show the existence of a {\it $\ZZ$-projective structure} transversal to the foliation on the Jacobi manifold (cf.\ Section~\ref{sec:period-lattices-as} and Appendix~\ref{apendixC});
\item we provide cohomological criteria to determine whether there exists a~contact isotropic realisation inducing a given transversal $\ZZ$-projective structure and to classify, up to isomorphism, all such realisations (cf.\ Theorem~\ref{thm:main}).
\end{itemize}

The methods employed in this paper are completely dif\/ferent from those of \cite{Banyaga,dd,jovanovic,lerman_contact_toric,wolbert}, for the main tool used here is the classical Spencer operator and its properties (developed in greater generality in~\cite{Maria}), which allows to deal with Jacobi structures def\/ined on non-trivial line bundles as well as providing a more conceptual approach to the case of Jacobi brackets on trivial line bundles. For instance, en route to proving Theorem~\ref{thm:main}, a~{\it local} contact classif\/ication of contact isotropic realisations is attained, providing a dif\/ferent, more intrinsic proof to \cite[Theorem~4]{jovanovic}\footnote{The relation between contact isotropic realisations and integrable systems on contact manifolds is going to be explored in depth in a separate paper.}.

Lastly, we go back to Poisson manifolds, for which there is a natural question: what relation is there between existence of a contact isotropic realisation and that of a symplectic isotropic realisation? The methods developed in this paper allow to tackle the above question in the case in which a Poisson manifold admits transversal $\ZZ$-af\/f\/ine and $\ZZ$-projective structures which are `related'. First, it is shown that a transversal $\ZZ$-af\/f\/ine structure satisfying an `integrality' condition (cf.\ Def\/inition~\ref{defn:strong_Z_aff}) induces naturally a~transversal $\ZZ$-projective structure (cf.\ Corollary~\ref{cor:strong_Z_proj}). Intuitively speaking, this relation should be viewed as analogous to the following construction: given an integral symplectic manifold, it is possible to construct a~Boothby--Wang-type contact manifold (cf.\ Example~\ref{sec:motiv-preq-circle}). In fact, the last main result of this paper generalises the above construction: given transversal $\ZZ$-af\/f\/ine and $\ZZ$-projective structures (denoted by~$\Xi$ and~$\Sigma$ respectively) on a Poisson manifold $\PM$ related as above, there exists a contact isotropic reali\-sa\-tion of $\PM$ inducing~$\Sigma$ if and only if there exists a symplectic isotropic reali\-sa\-tion of~$\PM$ inducing~$\Xi$ whose total space has an integral symplectic form (cf.\ Theorem~\ref{thm:cir_vs_sir}).

\looseness=-1 The structure of the paper is as follows. Section~\ref{sec:motiv-mult-free} introduces the basic notions used throughout the paper, sets the notation and provides motivating families of examples to study contact isotropic realisations. Seeing as the approach to contact and Jacobi manifolds taken in this paper follows~\cite{Maria-Crainic} and dif\/fers from the standard one taken in most other works in the literature, Sections~\ref{sec:basic-notions} and~\ref{sec:lie-algebr-assoc} provide a summary of some of the notions and of the results of~\cite{Maria-Crainic}. Section~\ref{sec:cont-isotr-real} def\/ines contact isotropic realisations and establishes their basic properties, while Section~\ref{sec:smooth-invar-cirs} tackles the problem of classifying contact isotropic realisations up to isomorphism. Throughout Sections~\ref{sec:motiv-mult-free}--\ref{sec:smooth-invar-cirs}, the case of Jacobi structures with trivial coef\/f\/icients is presented as an example to illustrate the various notions introduced. The case of Poisson manifolds is considered in Section~\ref{sec:sympl-vers-cont}, which moreover compares symplectic and contact isotropic realisations. In particular, Section~\ref{sec:cont-isotr-real-1} interprets the main results of Sections~\ref{sec:cont-isotr-real} and~\ref{sec:smooth-invar-cirs} in the case in which the contact structures are co-oriented and the Jacobi structures are def\/ined on trivial line bundles. For readers who are not familiar with the techniques for the general case, it may be helpful to read Section~\ref{sec:cont-isotr-real-1} alongside Sections~\ref{sec:cont-isotr-real} and~\ref{sec:smooth-invar-cirs}. To make the transition to the general case, it is useful to think that many results follow from choosing local trivialisations of the line bundle, thus reducing the problem to the case of Jacobi structures over a trivial line bundle; the classical Spencer operator and the associated language can be viewed as a way to bypass the above approach. Furthermore, to make the exposition smoother, some longer or technical proofs of results used throughout the paper have been placed in appendices at the end. Appendix~\ref{sec:proofs-prop-regul} deals with the geometric properties of the Jacobi manifolds that admit contact isotropic realisations. Appendices~\ref{sec:proofs-sect-refs} and~\ref{apendixC} provide proofs for some results of Section~\ref{sec:motiv-famil-exampl}, and of Sections~\ref{sec:period-lattices-as} and~\ref{sec:class-cir-over}, respectively.

Throughout this note, there are two types of comments, labelled {\bf Note} and {\bf Remark} respectively; those with the former label are central to the problems studied in this paper, while those with the latter may be skipped at a f\/irst reading.

{\bf Notation and conventions.} Throughout the paper, the intersection of two subsets $U_i \cap U_j$ is denoted by $U_{ij}$. All line bundles considered in this paper are real unless otherwise stated.

\section{Basic notions and motivation} \label{sec:motiv-mult-free}
The aim of this section is to recall fundamental notions regarding Jacobi structures, to establish notation, and to provide a family of motivating examples for the rest of the paper. Sections~\ref{sec:basic-notions} and~\ref{sec:lie-algebr-assoc} follow the approach to contact and Jacobi manifolds of~\cite{Maria-Crainic}, which should be considered as the main reference for any detail missing below. More details regarding Jacobi structures and their properties can be found in~\cite{jacobi,Dazord-Lich-Marle91,Lichnerowicz,Kirillov} amongst others.

\subsection{Jacobi structures: def\/initions and examples}\label{sec:basic-notions}
\begin{Definition}\label{defn:jacobi_structure} A {\it Jacobi structure} on a manifold $P$ is a pair $(L,\{\cdot,\cdot\})$ consisting of a line bundle $L\to P$, and a~{\em local} Lie bracket $\{\cdot,\cdot\} \colon \Gamma(L) \times \Gamma(L) \to \Gamma(L)$, i.e.,
\begin{gather*}
 \operatorname{supp}(\{u,v\})\subset \operatorname{supp}(u)\cap \operatorname{supp} (v)\qquad \forall \, u,v\in\Gamma(L),
\end{gather*}
where $\operatorname{supp}(u)$ denotes the support of $u$. A {\it Jacobi manifold} is a triple $(P,L,\{\cdot,\cdot\})$, where $(L,\{\cdot,\cdot\})$ is a Jacobi structure on~$P$.
\end{Definition}

Jacobi structures simultaneously generalise contact and Poisson structures.

\begin{Definition}\label{defn:contact_structure} A {\it contact structure} on a manifold $M$ is a smooth hyperplane distribution $H\subset TM$ whose curvature map
 \begin{gather}\label{curvature}
 c\colon \ H\times H\to T M/H,
 \end{gather}
def\/ined on sections by $c(X,Y):=[X,Y] \text{ mod } H$ is f\/ibre-wise non-degenerate. A {\it contact manifold} is a pair $(M,H)$, where $H$ is a~contact structure on~$M$.
\end{Definition}

\begin{Note}\label{obs:ct_mfs} The majority of works in the literature on contact geometry concentrates on the case in which $H = \ker \theta$, for some 1-form $\theta \in \Omega^1(M)$ whose dif\/ferential $\de \theta$ makes $H \to M$ into a symplectic vector bundle. In this case the bundle $TM/ H$ is trivial and the above map~$c$ coincides with $\de \theta|_H$. Such contact structures are henceforth referred to as being {\em co-oriented}. More generally, a contact structure~$H$ on~$M$ can be encoded equivalently as the kernel of the 1-form $\theta \in \Omega^1(M,L)$ given by projection onto $L: = T M/H$. This is henceforth referred to as the {\it $($canonical$)$ generalised contact form} of~$\CM$. Throughout this paper, both points of views are used interchangeably.
\end{Note}

The following family of examples of contact manifolds plays a~prominent role throughout this paper.

\begin{Example}\label{ex:jet_bundle}
Let $\pi\colon L\to P$ a line bundle and (by abuse of notation) denote by $\pi \colon J^1L\to P$ the {\it first jet bundle}
\begin{gather*} J^1 L|_x=\big\{j^1_xu\,|\, u\in\Gamma(L)\big\}.\end{gather*}
The Cartan contact form $ \theta_{\mathrm{can}}\in\Omega^1(J^1L,\operatorname{pr}^*L)$ def\/ines a contact structure as in Note~\ref{obs:ct_mfs}. It is def\/ined by
\begin{gather*}
 \theta_{\mathrm{can}, j^1_xu}=D_x(\operatorname{pr}-u\circ\pi)\colon \ T_{j^1_xu}\big(J^1L\big)\to L_{x},
\end{gather*}
with $\operatorname{pr}\colon J^1L\to L$, $j^1u_x\mapsto u(x)$, and where we have used the canonical identif\/ication $L_x\simeq T_{u(x)}(L_x).$ The Cartan contact form detects {\it holonomic} sections of~$J^1L$, i.e., those of the form $j^1u\colon x\mapsto j^1_xu$, for some $u\in\Gamma(L)$, in the sense that a~section~$\xi$ of~$J^1L$ is holonomic if and only if $\xi^*\theta_\text{can}=0$.
\end{Example}

\begin{Example}\label{exm:ctc_jacobi}
A contact manifold $\CM$ comes equipped with a well-known natural Jacobi structure $(L = T M/H, \{\cdot,\cdot\})$ which can be described as follows (cf.~\cite{Maria-Crainic,Dazord} and references therein). A {\it Reeb vector field} of~$\CM$ is any vector f\/ield~$R$ satisfying
 \begin{gather*}[R,\Gamma(H)]\subset\Gamma(H);\end{gather*}
the vector spaces of Reeb vector f\/ields is denoted by $\X_{\mathrm{Reeb}}(M,H)$. The map
 \begin{gather*} 
 \X_{\mathrm{Reeb}}(M,H) \to \Gamma(L), \qquad R \mapsto \theta(R)
 \end{gather*}
is a vector space isomorphism (cf.\ \cite[Lemma~2.2 and Corollary~2.3]{Maria-Crainic}). The inverse image of $u\in\Gamma(L)$ under the above isomorphism, denoted by~$R_u$, is called the {\it Reeb vector field associated to~$u$}. The space~$\X_{\mathrm{Reeb}}(M,H)$ is closed under the Lie bracket of vector f\/ields, hence it induces a~Lie bracket on~$\Gamma(L)$, explicitly given by
 \begin{gather} \label{eq:2}
 \{u,v\} := [R_u,R_v] \quad \text{mod} \ H,
 \end{gather}
for $u,v \in \Gamma(L)$. In the case in which $TM/H \to M$ is trivial, i.e., in the co-orientable case (cf.\ Note~\ref{obs:ct_mfs}), the above Jacobi structure can be described easily as follows. Given a contact $1$-form $\theta \in \Omega^1(M)$ with $H = \ker \theta$, let $R_{\theta} \equiv R_1$ be {\em the} Reeb vector f\/ield associated to (the function $1$ and to) $\theta$. Then the space of Reeb vector f\/ields is nothing but the space of Reeb vector f\/ields associated to all contact forms def\/ining $H$ and the induced bracket is well-known (cf.\ Example~\ref{obs:kir} below).
\end{Example}

\begin{Example}\label{obs:kir} A {\it Poisson structure} on a manifold $P$ is a bivector f\/ield $\Lambda \in \X^2(P)$, which satisf\/ies $\llbracket \Lambda,\Lambda \rrbracket=0$, where $\llbracket \cdot, \cdot \rrbracket$ is the Schouten bracket. We say that the pair $(P,\Lambda)$ is a~{\it Poisson manifold}. In this case $\Lambda$ induces a Jacobi structure $\{\cdot,\cdot\}$ on the trivial line bundle $\RR_P$ by observing that $\Lambda$ determines a~local Lie bracket on the sections of the trivial bundle $\RR_P \to P$. By abuse of notation, the induced Jacobi manifold is also
denoted by $(P,\Lambda)$. More generally, a~Jacobi structure $(\RR_P,\{\cdot, \cdot\})$ on $P$ is completely determined by a pair
$(\Lambda, R) \in \X^2(P) \times \X(P)$, satisfying
\begin{gather*} 
 \llbracket \Lambda,\Lambda \rrbracket =2R\wedge\Lambda,\qquad \llbracket \Lambda,R \rrbracket=0.
\end{gather*}
The Lie bracket on $\Gamma(\RR_P) = C^{\infty}(P)$ is given by
\begin{gather*}
 \{f,g\} := \Lambda(\de f, \de g) + f (Rg) - g(Rf),
\end{gather*}
for $f,g \in C^{\infty}(P)$ (cf.~\cite{Lich78}).
\end{Example}

Morphisms between Jacobi manifolds are def\/ined as follows.

\begin{Definition}\label{defn:jacobi_maps}
 Let $(N, L_N,\{\cdot,\cdot\}_N)$ and $(P, L_P, \{\cdot,\cdot\}_P)$ be Jacobi manifolds such that there exists an isomorphism $F\colon \phi^*L_P \to
 L_N$. A map $\phi \colon N \to P$ is said to be {\it Jacobi} with {\it bundle component} $F$ if for all $u,v \in \Gamma(L_P)$
 \begin{gather*}
 \{F \circ \phi^*u, F \circ \phi^*v\}_N = F \circ \phi^*\{u,v\}_P.
 \end{gather*}
\end{Definition}

\begin{Remark}
 If the line bundles $L_N$ and $L_P$ are assumed to be trivial in Def\/inition~\ref{defn:jacobi_maps}, then the above notion of Jacobi map with bundle component is often referred to in the literature as a {\em conformal} Jacobi morphism (cf.\ \cite[Section~1.6]{Dazord-Lich-Marle91}).
\end{Remark}

\begin{Example}\label{exm:symplectisation}
Let $\CM$ be a contact manifold with associated line bundle $L \to M$, where the notation is as in Example~\ref{exm:ctc_jacobi}. The submanifold $L^* \setminus \{0\} \hookrightarrow T^*M$ is symplectic, where $T^*M$ is endowed with the canonical symplectic form (cf.\ \cite[Def\/inition~2.3]{lerman_contact_toric}). Let $\Omega$ denote the symplectic form on $L^* \setminus \{0\}$; the symplectic manifold $(\SYM, \Omega)$ is said to be the {\it symplectisation}\footnote{This notion dif\/fers slightly with that in the literature on co-orientable contact structures, where a connected component of $\SYM$ is declared to be the symplectisation of~$\CM$.} of~$\CM$. The projection $\operatorname{pr} \colon (\SYM,\Omega) \to \CM$ is a Jacobi map with bundle component $F_{\operatorname{pr}} \colon \operatorname{pr}^* L \to \RR_{\SYM}$
 given by $F_{\operatorname{pr}}(\alpha, u) = \alpha(X)$, where $u = X \text{ mod } H$.
\end{Example}

\begin{Example}\label{exm:proj_dual_lie_algebra} Let $\mathfrak{g}^*$ denote the dual of a f\/inite-dimensional real Lie algebra; endow it with the standard linear Poisson structure $\Lambda$. Its projectivisation $\pg : =\frac{\mathfrak{g}^* \setminus \{0\}}{\RR^*}$ admits a natural Jacobi structure $(O(1),\{\cdot,\cdot\})$, where $O(1) \to \pg$ is the dual of the tautological line bundle. The bracket $\{\cdot,\cdot\}$ on $ \Gamma(O(1))$ can be def\/ined as follows. Recall that $\Gamma(O(1))$ can be identif\/ied with the vector space of smooth homogeneous functions of degree 1 on $\mathfrak{g}^* \setminus \{0\}$, denoted by $\mathrm{C}^{\infty}_{h}(\mathfrak{g}^* \setminus \{0\})$. Since $\mathfrak{g}^* \setminus \{0\} \subset \mathfrak{g}^*$ is an open subset which is the complement of a~symplectic leaf, $\Lambda$ restricts to a~Poisson structure on $\mathfrak{g}^* \setminus \{0\}$, also denoted by $\Lambda$. Linearity of $\Lambda$ implies that $\mathrm{C}^{\infty}_{h}(\mathfrak{g}^* \setminus \{0\})$ is a Lie subalgebra of $(\mathrm{C}^{\infty}(\mathfrak{g}^* \setminus \{0\}), \{\cdot,\cdot\}_{\Lambda})$. This def\/ines a~Lie bracket on $\Gamma(O(1))$ which is manifestly local and, hence, a Jacobi structure $(O(1),\{\cdot,\cdot\})$ as claimed. If $\mathbb{S}(\mathfrak{g}^*) = \frac{\mathfrak{g}^* \setminus \{0\}}{\RR^+}$, the covering map $q \colon \mathbb{S}(\mathfrak{g}^*) \to \pg$ can be used to def\/ine a Jacobi structure $(q^*(O(1)),\{\{\cdot,\cdot\}\})$. The vector bundle $q^*(O(1)) \to \mathbb{S}(\mathfrak{g}^*)$ is trivialisable and a~trivialisation can be obtained by picking a metric on $\mathfrak{g}$; once such a choice is made, this example recovers the Jacobi structures considered in~\cite{deleon,Lich78}, and \cite[Example~2.3]{Zambon}. Finally, the natural projection $\pi \colon (\mathfrak{g}^* \setminus \{0\},\Lambda) \to (\pg, O(1),\{\cdot,\cdot\})$ is a~Jacobi map with bundle component $F_{\pi}\colon \pi^*(O(1)) \to \RR_{\mathfrak{g}^* \setminus \{0\}}$, given by $F_{\pi}(x,\eta) = \eta(x)$.
\end{Example}

\subsection{The Lie algebroid and Spencer operator of a Jacobi manifold}\label{sec:lie-algebr-assoc}
In analogy with Poisson manifolds, the geometric structure of a Jacobi manifold can be completely encoded by a Lie algebroid together with a Spencer operator (cf.\ \cite{Maria-Crainic,Maria,Dazord,Kerbrat}). In general, a {\it Lie algebroid} over a manifold $M$ is a vector bundle $\mathsf{A} \to M$, together with a Lie bracket $[\cdot,\cdot]$ on $\Gamma(\mathsf{A})$, and a~vector bundle map $\rho\colon \mathsf{A}\to TM$, called the {\it anchor map}, satisfying the compatibility condition
\begin{gather*}[\alpha,f\beta]=f[\alpha,\beta]+L_{\rho(\alpha)}(f)\beta\end{gather*}
for all $f\in C^{\infty}(M),$ and all $\alpha,\beta \in \Gamma(\mathsf{A})$. Given a Jacobi manifold $(P,L,\{\cdot,\cdot\})$, the f\/irst jet bundle $J^1 L \to P$ can be endowed with the structure of a Lie algebroid as follows (cf.~\cite{Maria-Crainic} for a~proof).

\begin{Proposition}\label{prop:alg}
Given a Jacobi manifold $(P, L,\{\cdot,\cdot\})$, there exists a~Lie algebroid structure on $J^1 L \to P$ which is uniquely characterised by the following properties:
 \begin{enumerate}[label={\rm (\Roman*)}, ref={\rm (\Roman*)}]\itemsep=0pt
 \item\label{item:A} The anchor map $\rho\colon J^1 L\to T P$ satisfies
 \begin{gather*}\{u,fv\}=f\{u,v\}+L_{\rho(j^1u)}(f)v,\end{gather*}
 for all $u,v\in\Gamma(L),$ $f\in C^\infty(P)$;
 \item\label{item:B} The Lie bracket on $\Gamma(J^1 L)$ satisfies
 \begin{gather*}[j^1u,j^1v]=j^1\{u,v\},\qquad \forall \, u,v\in\Gamma(L).\end{gather*}
 \end{enumerate}
 This Lie algebroid is henceforth referred to as the associated Lie algebroid to $(P,L,\{\cdot,\cdot\})$.
\end{Proposition}

\begin{Note}\label{obs:trivial_case} If $L=\RR_P$, the above Lie algebroid agrees with that def\/ined in \cite[Def\/inition~1.4]{jacobi}, where the identif\/ication $j^1_xf\mapsto (\de_x f,f(x))$ of $J^1\RR_P$ with $\cotan P\times\RR $ is used. In this case the anchor $\rho\colon \cotan P\times\RR\to TP$ becomes
 \begin{gather*}\rho(\omega,\lambda)=\Lambda^\sharp(\omega)+\lambda R,\end{gather*}
where $(\Lambda,R)$ are as in Example~\ref{obs:kir}. In particular, if the underlying structure is Poisson (i.e., $R=0$), then the kernel of the anchor contains $0 \oplus \RR$.
\end{Note}

\begin{Note}\label{obs:Jac_maps}
The condition for a map to be Jacobi can be formulated in terms of anchors of Lie algebroids. Suppose that $\phi\colon (N, L_N,\{\cdot,\cdot\}_N) \to (P, L_P, \{\cdot,\cdot\}_P)$ is a Jacobi map with bundle component $F$, then the following diagram commutes
 \begin{gather*}
 \xymatrix{J^1L_N \ar[r]^-{\rho_N} & TN \ar[d]^-{D\phi} \\
 J^1L_P \ar[u]^-{F \circ \phi^*} \ar[r]_-{\rho_P} & TP.}
 \end{gather*}
\end{Note}

\begin{Remark}\label{rk:leaves}
 The image of the anchor of a Lie algebroid def\/ines a (singular) foliation on the base manifold. Unlike what happens in Poisson geometry, the leaves of (the foliation induced by the Lie algebroid associated to) a Jacobi manifold may be even or odd dimensional. If the induced foliation is regular, i.e., integrable in the sense of Frobenius, the underlying Jacobi manifold is said to be {\it regular}.
\end{Remark}

In analogy with what happens for Poisson manifolds, the Lie algebroid associated to a Jacobi manifold comes with compatible extra structure. Fix any line bundle $L \to P$ (not necessarily with a Jacobi bracket); its f\/irst jet bundle $J^1L$ f\/its into a short exact sequence of vector bundles
\begin{gather*}
 0\rmap \cotan P \otimes L \stackrel{i}{\to} J^1L \stackrel{\operatorname{pr}}{\to} L\rmap 0,
\end{gather*}
where $ i(df\otimes u)= f j^1(u)- j^1(fu)$ for $f \in C^{\infty}(P)$ and $u \in \Gamma(L)$. While the above sequence is {\em not} canonically split, the map $u\mapsto j^1u$ gives a canonical splitting at the level of sections:
\begin{gather*} \Gamma\big(J^1L\big)\cong \Gamma(L)\oplus\Omega^1(P;L), \end{gather*}
known as the {\it Spencer decomposition}.

\begin{Definition}\label{defn:cso}
The {\it classical Spencer operator} associated to $L \to P$ is the projection $\mathrm{D}\colon$ $\Gamma(J^1L) \to\Omega^1(P;L)$.
\end{Definition}

\begin{Example}\label{exm:spencer_trivial}
Suppose that $L = \RR_P$, i.e., it is trivial; then using the identif\/ication of $J^1 \RR_P \cong \cotan P \oplus \RR_P$ of Note~\ref{obs:trivial_case}, identify an element of $\Gamma(J^1 \RR_P)$ with a pair $(\eta,f)$, where $\eta \in \Omega^1(P)$ and $f \in C^{\infty}(P)$. Then the classical Spencer operator is given by $\mathrm{D}(\eta,f) = \de f - \eta$.
\end{Example}

\begin{Note}\label{obs:prop_cso}
The classical Spencer operator $\mathrm{D}$ associated to $L \to P$ is completely determined by the following two conditions:
 \begin{itemize}\itemsep=0pt
 \item for any $u \in \Gamma(L)$, $\mathrm{D}(j^1u)=0$;
 \item for any $X \in \X(P)$, $\alpha \in \Gamma(J^1L)$ and $f \in C^{\infty}(P)$,
 \begin{gather*}\mathrm{D}_X(f\alpha)=f\mathrm{D}_X(\alpha)+df(X)\operatorname{pr}(\alpha),\end{gather*}
 i.e., the Leibniz identity holds.
 \end{itemize}
Moreover, the $C^{\infty}(P)$-module structure on $\Gamma(J^1L) $ induced by the Spencer decomposition is the following
 \begin{gather*} f\cdot (u, \phi)= (f u, \phi+ df\otimes u),\end{gather*}
for any $f \in C^{\infty}(P)$, $u \in \Gamma(L)$ and $\phi \in \Omega^1(P;L)$. Finally, the classical Spencer operator associated to $L \to P$ is related to the Cartan contact form $\theta_{\mathrm{can}}$ on $J^1L$ (cf.\ Example~\ref{ex:jet_bundle}) as follows. If $\alpha \in \Gamma (J^1L )$, then
\begin{gather*}
\mathrm{D}(\alpha) = \alpha^* \theta_{\mathrm{can}}.
\end{gather*}
\end{Note}

Suppose that $(P, L ,\{\cdot, \cdot\})$ is a Jacobi manifold. In this case, the classical Spencer operator is compatible with the Lie algebroid structure on $J^1L \to P$ def\/ined in Proposition~\ref{prop:alg} in the following sense. The Jacobi bracket def\/ines a f\/lat $J^1L$-connection on $L \to P$, i.e., an operator $\nabla\colon \Gamma(J^1L)\times\Gamma(L)\to \Gamma(L)$ which is $C^\infty(P)$-linear in the f\/irst component, satisf\/ies the Leibniz identity $\nabla_\alpha(fu)=f\nabla_\alpha(u)+L_{\rho(\alpha)}(f)u$, and the f\/latness equation $\nabla_{[\alpha,\beta]}=\nabla_\alpha\nabla_\beta-\nabla_\beta\nabla_\alpha$, for any $f \in C^{\infty}(P)$, $\alpha,\beta \in \Gamma(J^1L)$
and $u \in \Gamma(L)$. It is uniquely def\/ined by the formula
\begin{gather}\label{eq:nabla}
 \nabla\colon \ \Gamma\big(J^1L\big)\times\Gamma(L)\to\Gamma(L),\qquad \nabla_{j^1u}(v)= \{u,v\}.
\end{gather}
The compatibility of $\mathrm{D}$ with the Lie algebroid structure on~$J^1L \to P$ can be encoded in the following two equalities
\begin{gather}
 \label{horizontal}\mathrm{D}_{\rho(\alpha)}(\alpha') =\nabla_{\alpha'}(\pr(\alpha))+ pr([\alpha, \alpha']), \\
 \label{vertical} \mathrm{D}_X[\alpha,\alpha']=\nabla_\alpha(\mathrm{D}_X\alpha')- \mathrm{D}_{[\rho(\alpha),X]}\alpha'-\nabla_{\alpha'}(\mathrm{D}_X\alpha)+\mathrm{D}_{[\rho(\alpha'),X]}\alpha,
\end{gather}
for any $\alpha,\alpha'\in\Gamma(J^1L)$ and $X \in \X(P)$ (cf.~\cite{Maria} for a general notion of Spencer operators). Henceforth, the classical Spencer operator $\mathrm{D}$ associated to a line bundle $L \to P$ endowed with a Jacobi bracket is referred to as {\it the Spencer operator}
associated to the Jacobi manifold $(P,L,\{\cdot,\cdot\})$.

\subsection{Motivating example: Multiplicity-free actions on contact manifolds}\label{sec:motiv-famil-exampl}

Throughout this section, let $G$ denote a compact connected Lie group, and denote symplectic and contact manifolds by $\SM$ and $\CM$ respectively. If $\mathfrak{g} = \operatorname{Lie}(G)$, $G$ acts on $\mathfrak{g}^*$ by the coadjoint action, which is henceforth understood to be {\em the} $G$-action on~$\mathfrak{g}^*$. Proofs of the main results of this subsection can be found in Appendix~\ref{sec:proofs-sect-refs}. First the most basic notions of symmetry in symplectic and contact geometry are recalled for completeness.

\begin{Definition}\label{defn:ham_ctct}\quad
 \begin{itemize}\itemsep=0pt
 \item An action $G \curvearrowright \SM$ is said to be {\it Hamiltonian} if there exists a smooth $G$-equivariant map $\chi\colon S \to \mathfrak{g}^*$, called {\it moment map}, such that for all $\xi \in \mathfrak{g} = \operatorname{Lie}(G)$,
 \begin{gather*}
 \omega(X_{\xi}, -) = \de \langle \chi, \xi \rangle,
 \end{gather*}
 where $X_{\xi}$ is the vector f\/ield on $M$ induced by $\xi$ and $\langle \cdot,\cdot \rangle$ denotes the standard pairing between $\mathfrak{g}$ and $\mathfrak{g}^*$.
 \item An action $G \curvearrowright \CM$ is said to be {\it contact} if it preserves the contact structure~$H$, i.e., for all $g \in G$, $p \in M$, $D_pg(H_p) = H_{g \cdot p}$.
 \end{itemize}
\end{Definition}

\begin{Note}\label{obs:contact_action_lift}
Suppose that $G \curvearrowright \CM$ is contact. Its lift to the cotangent bundle $G \curvearrowright (T^*M,\omega_{\mathrm{can}})$ def\/ined by $(g \cdot \alpha)(X) = \alpha(D g^{-1}(X))$, where $\alpha \in T^*M$ and $X \in TM$, is Hamiltonian with moment map $\mu\colon T^*M \to \mathfrak{g}^*$ given by $\mu(\alpha)(\xi) = \alpha (D \operatorname{pr} (X_{\xi} ) )$, where $\xi \in \mathfrak{g}$, $X_{\xi}$ is the vector f\/ield on $T^*M$ by $\xi$, and $\operatorname{pr} \colon T^*M \to M$ is the projection. This action preserves $L^* \setminus \{0\}$, thus inducing a Hamiltonian action on the symplectisation $(\SYM,\Omega)$ whose moment map $\mu \colon \SYM \to \mathfrak{g}^*$ is $\RR^*$-equivariant, i.e., for all $t \in
\RR^*$ and all $\alpha \in \SYM$, $\mu(t\alpha) = t\mu(\alpha)$.
\end{Note}

For the purposes at hand, it is useful to recall some simple notions associated to Lie group actions.

\begin{Definition}\label{defn:action}
 An action of $G$ on a manifold $N$ is said to be
 \begin{itemize}\itemsep=0pt
 \item {\it locally free} at $p \in N$ if the stabiliser $G_p := \{ g \in G \,|\, g \cdot p = p \}$ is discrete;
 \item {\it effective} if $\bigcap\limits_{p \in N} G_p = \{e\}$.
 \end{itemize}
\end{Definition}

The following def\/inition introduces a notion of `integrable' group actions of compact Lie groups on symplectic and contact manifolds (cf.\ \cite[Def\/inition 5.1.11]{gui_sja} for the symplectic case and \cite[Def\/inition~2.11]{lerman_contact_toric} for torus actions on contact manifolds).

\begin{Definition}\label{defn:mult_free}
Let $G$ be a compact Lie group. An action $G \curvearrowright \SM$ (respectively $G \curvearrowright \CM$) is said to be {\it multiplicity-free} if it is Hamiltonian (respectively contact), ef\/fective, locally free at some point, and if $\dim S = \dim G + \operatorname{rk} G$ (respectively $ \dim M = \dim G + \operatorname{rk} G - 1$), where $\operatorname{rk} G$ denotes the rank of $G$, i.e., the dimension of its maximal torus.
\end{Definition}

In the context of Hamiltonian actions of compact Lie groups, multiplicity-free actions can be thought of as being `integrable'. More precisely, if a compact Lie group $G$ acts on~$\SM$ so that the action is locally free at some point, then
\begin{gather*}
 \dim S \geq \dim G + \operatorname{rk} G
\end{gather*}
(cf.\ \cite[Theorem~5.1.6]{gui_sja}). The condition of being multiplicity-free corresponds to considering the case in which $\dim S$ is completely determined by the group acting (e.g., when $G = \mathbb{T}^n$, then $\dim S = 2n$). Moreover, the following lemma illustrates the relation between the notions of multiplicity-free actions in contact and symplectic geometry by considering symplectisations.

\begin{Lemma}\label{prop:lift}
If a contact action $G \curvearrowright \CM$ is multiplicity-free then the induced Hamiltonian action on $(\SYM,\Omega)$ is too.
\end{Lemma}

\begin{proof} To show that the induced action on $(\SYM,\Omega)$ is multiplicity-free, it suf\/f\/ices to check that the action is locally free at some point and ef\/fective, as the condition on the dimension of~\smash{$\SYM$} and on~$G$ is automatically satisf\/ied. Both claims follow by observing that the projection $\operatorname{pr} \colon \SYM \to M$ is $G$-equivariant, thus implying that, for all $\alpha \in \SYM$, \smash{$G_{\alpha} \subset G_{\operatorname{pr}(\alpha)}$}.
\end{proof}

\begin{Example}\label{exm:prequant_mult_free}
Let $G$ be a (compact) simply connected Lie group acting on a closed symplectic manifold $\SM$ in a multiplicity-free fashion. Denote the moment map by $\chi \colon \SM \to \mathfrak{g}^*$ and, for each $\xi \in \mathfrak{g}$, set $\chi^{\xi}:=\langle \chi,\xi \rangle \colon S \to \RR$. Suppose further that~$\omega$ is integral and consider the prequantum circle bundle $ \phi \colon \CM \to \SM$ as in Example~\ref{sec:motiv-preq-circle}. The inf\/initesimal $\mathfrak{g}$-action on~$S$ can be lifted to an inf\/initesimal $\mathfrak{g}$-action on~$M$, by setting $\xi \mapsto R_{\phi^* \chi^{\xi}}$, where $R_{\phi^* \chi^{\xi}}$ is the Reeb vector f\/ield associated to the function $\phi^*\chi^{\xi}$. This holds because $\phi \colon \CM \to \SM$ is a~{\em contact (isotropic) realisation} (cf.\ Def\/initions~\ref{defn:CR} and~\ref{defn:CIR}). Since~$M$ is compact, this action can be integrated to a $G$-action on $M$ which is, in fact, contact; moreover, $\phi$ is $G$-equivariant. This $G$-action commutes with the (contact!) $S^1$-action coming from the principal $S^1$-bundle structure on $\phi \colon M \to S$; this follows because the functions $\phi^*\chi^{\xi}$ are basic. This yields a contact $G \times S^1$-action on $\CM$ which is multiplicity-free since the $G$-action on $\SM$ is and because $\phi \colon \CM \to \SM$ is $G$-equivariant.
\end{Example}

\begin{Example}\label{exm:S3}
Not all examples of multiplicity-free actions in contact geometry arise as in Example~\ref{exm:prequant_mult_free}. For instance, consider $M = S^3 \cong \mathrm{SU}(2)$ endowed with a left-invariant contact structure $H$ (cf.~\cite{boothby_wang}). The $\mathrm{SU}(2)$-action on $\CM$ given by left-translations is multiplicity-free and cannot arise as in Example~\ref{exm:prequant_mult_free} because $\mathrm{SU}(2)$ is simply-connected.
\end{Example}

Multiplicity-free actions on contact manifolds satisfy a strong property which generalises \cite[Lemma~2.12]{lerman_contact_toric}.

\begin{Proposition}\label{lemma:mu_misses_zero}
Let $G \curvearrowright \CM$ be a multiplicity-free action and denote the moment map of the induced Hamiltonian action on $(\SYM,\Omega)$ by $\mu \colon \SYM \to \mathfrak{g}^*$. Then $\mu ( \SYM ) \subset \mathfrak{g}^* \setminus \{0\}$.
\end{Proposition}

The proof of Proposition \ref{lemma:mu_misses_zero} is postponed to Appendix~\ref{sec:proofs-sect-refs}. An immediate consequence of Proposition~\ref{lemma:mu_misses_zero} and Note~\ref{obs:contact_action_lift} is the following.

\begin{Corollary}\label{cor:Jacobi_mm}
Associated to a multiplicity-free action $G \curvearrowright \CM$ is a smooth map $\phi \colon M \to \pg$ which makes the following diagram commutative
 \begin{gather*} 
 \xymatrix{L^*\setminus \{0\} \ar[d]_-{\operatorname{pr}} \ar[r]^-{\mu}
 & \mathfrak{g}^* \setminus \{0\} \ar[d]^-{\pi} \\
 M \ar[r]_-{\phi} & \pg.}
 \end{gather*}
\end{Corollary}

This motivates introducing the following notion.

\begin{Definition}\label{defn:jmp}
Given a multiplicity-free action $G \curvearrowright \CM$, the map $\phi \colon M \to \pg$ as in Corollary \ref{cor:Jacobi_mm} is called the {\it Jacobi moment map} associated to the action.
 \end{Definition}

\begin{Example}\label{exm:mom_map_prequant}
Fix notation as in Example~\ref{exm:prequant_mult_free}. Suppose that $G \times S^1 \curvearrowright \CM$ is multiplicity-free, where $\CM$ is the total space of a principal $S^1$-bundle of a closed symplectic manifold $\SM$ endowed with a multiplicity-free action of the compact, simply connected Lie group~$G$. The symplectisation of $\CM$ can be written as $(M \times \RR^*, \de(t\alpha))$, where $\alpha \in \Omega^1(M)$ is a~$G \times S^1$-invariant contact form. The moment map for the lifted Hamiltonian $G \times S^1$-action on \smash{$(M \times \RR^*, \de(t\alpha))$} can be written as $\mu(m,t) = (t\mu(\Phi(m)),t) \in \mathfrak{g}^* \times \RR$, where $\Phi \colon M \to S$ is the projection and $\mu \colon S \to \mathfrak{g}^*$ is the moment map of the Hamiltonian $G$-action. Since $t \in \RR^*$ in the above formula, it follows that the image of $\mu$ does not contain $(0,0) \in \mathfrak{g} \times \RR$. In this case, the Jacobi moment map $\phi \colon M \to \mathbb{P}(\mathfrak{g}^* \times \RR)$ can be written as $\phi(m) = [\mu(\Phi(m)):1] \in \mathbb{P} (\mathfrak{g}^* \times \RR )$.
\end{Example}

The next theorem illustrates properties of the Jacobi moment map, which serve as motivation for the rest of this paper (cf.\ \cite[Proposition~3.4, Lemma~3.5 and Proposition~5.1]{Zambon} for statements analogous to~\ref{item:5} and~\ref{item:6} in the case of co-oriented contact
manifolds).

\begin{Theorem}\label{lemma:properties_Jac_mm}
The Jacobi moment map $\phi \colon \CM \to (\pg,O(1),\{\cdot,\cdot\} )$ associated to a~mul\-ti\-plicity-free contact action $ G \curvearrowright \CM$ is Jacobi with bundle component $F_{\phi} \colon \phi^*(O(1)) \to L$, where $F_{\phi}(p,\eta) \in L_p$ is the unique element such that, for all $\alpha \in L^*_p \setminus \{0\}$, $\eta(\mu(\alpha)) = \alpha ( F_{\phi}(p,\eta) )$, and $\mu \colon L^* \setminus \{0\} \to \mathfrak{g}^*$ is the moment map of the lifted action. Moreover, if $M_{\mathrm{prin}} \subset M$ denotes the open, dense subset consisting of points whose stabiliser is
 discrete, then for all $p \in M_{\mathrm{prin}}$,
 \begin{enumerate}[label={\rm (J\arabic*)}, ref={\rm (J\arabic*)}]\itemsep=0pt
 \item\label{item:5} $D_p\phi$ is onto;
 \item \label{item:6} $\ker D_p \phi$ is transversal to the contact distribution, i.e.,
 \begin{gather*}
 \ker D_p \phi + H_p = T_p M;
 \end{gather*}
 \item \label{item:14} $\ker D_p \phi \subset (\ker D_p \phi )^{\perp}$, where $ ( \ker D_p \phi )^{\perp} = \rho (F_{\phi} (\phi^* (J^1O(1) ) ) )$ and $\rho \colon J^1L \to TM$ is the anchor map of the Lie algebroid associated to~$\CM$.
 \end{enumerate}
\end{Theorem}

In light of Corollary \ref{cor:Jacobi_mm} and Theorem~\ref{lemma:properties_Jac_mm}, the theory of multiplicity-free actions on contact manifolds brings about Jacobi maps satisfying properties \ref{item:5}--\ref{item:14}. The rest of the paper concentrates on studying properties of such maps; this can be seen as both a f\/irst step towards understanding properties of multiplicity-free actions in the contact setting and a generalisation
of the aforementioned actions.

\section{Contact isotropic realisations of Jacobi manifolds}\label{sec:cont-isotr-real}

This paper studies the classif\/ication of a special type of contact realisations of Jacobi manifolds (cf.\ \cite[Def\/inition~3.7]{Zambon}), called {\it contact isotropic realisations} (or CIR for short). The aim of this section is to def\/ine such realisations and to give their basic properties.

\begin{Definition}\label{defn:CR}
A {\it contact realisation} of a Jacobi manifold $\JM$ is a contact manifold $\CM$, together with a surjective submersion $\phi\colon \CM \to \JM$ satisfying the following properties:
 \begin{enumerate}[label=(CR\arabic*),ref=(CR\arabic*),leftmargin=*]\itemsep=0pt
 \item \label{item:IR1}$\phi$ is a Jacobi map with bundle component $F\colon \phi^*(L)\to L_M$, where $L_M:= TM/H$ (cf.\ Def\/inition~\ref{defn:jacobi_maps})
 \item \label{item:IR2} $H = \ker \theta$ is transversal to $\phi$, i.e.,
 \begin{gather*}
 H+\Ker D\phi=\mathrm{T}M.
 \end{gather*}
 \end{enumerate}
\end{Definition}

Henceforth, whenever referring to a contact realisation $\phi\colon \CM \to \JM$, the dimensions of $M$ and of $P$ are f\/ixed to be $2n+1$ and $2n+1-k$ respectively.

\begin{Definition}\label{defn:CIR}
A {\it contact isotropic realisation} (CIR for short) of a Jacobi manifold $\JM$ is a contact realisation $\phi\colon \CM\to \JM$ with connected, compact f\/ibers satisfying the isotropic condition:
 \begin{enumerate}[label=(I),ref=(I)]\itemsep=0pt
 \item \label{item:IR3} $\Ker D\phi\subset (\Ker D\phi)^{\perp}$, where $ (\Ker D\phi)^{\perp} := \rho_M(F(\phi^*J^1L))$ is the {\em pseudo-orthogonal distribution} of $\ker D\phi$.
 \end{enumerate}
\end{Definition}

The isomorphism $F$ is henceforth assumed to be the identity to simplify the notation and the exposition.

\begin{Note}\label{item:po} For Jacobi manifolds with trivial bundle component $L=\RR_P$, condition~\ref{item:IR3} is equivalent to requiring that the subbundle $ \ker D\phi \cap H\subset H$ is isotropic. This can be seen directly by using property~\ref{item:IR3}, and the description of $H=\ker\theta$, by a contact form $\theta\in\Omega^1(M)$ (which trivialises $TM/H$). In this case, the Jacobi structure $(\Lambda_M, R_{\theta})$ is determined by $\theta$ and its Reeb vector f\/ield $R_{\theta}$ (cf.\ Example~\ref{exm:ctc_jacobi}), and $-c^\sharp:=\Lambda_M^\sharp\colon T^*M\to H\subset TM$ is uniquely def\/ined by the equations
 \begin{gather*}i_{\Lambda_M^\sharp(\omega)}\theta=0,\qquad i_{\Lambda_M^\sharp(\omega)}\de\theta=-\omega+\omega(R_{\theta})\theta\end{gather*}
 (see Note~\eqref{obs:trivial_case}). Thus,
 \begin{gather*}(\Ker D\phi)^{\perp}=\rho_M(\phi^*(T^*P)\times\RR)=\Lambda_M^\sharp((\ker
 D\phi)^\circ)\oplus \RR \langle R_{\theta} \rangle=(D\phi \cap H)^{\perp}\oplus \RR \langle R_{\theta} \rangle, \end{gather*}
where $(\ker D\phi)^\circ\subset T^*M$ denotes the annihilator of $\ker D\phi$, and, by abuse of notation, $(D\phi \cap H)^{\perp}$ denotes the symplectic orthogonal of $D\phi \cap H$ in $(H,\de\theta)$. This explains the terminology used in Def\/inition~\ref{defn:CIR}.
\end{Note}

\begin{Example}[prequantum circle bundles over symplectic manifolds]\label{sec:motiv-preq-circle} Let $\SM$ be a symplectic manifold whose symplectic form $\omega$ is {\it integral}, i.e., it def\/ines a cohomology class $[\omega] \in \HH^2(S;\ZZ)$. It is well-known that the total space of the principal $S^1$-bundle $\phi \colon M \to S$ with Chern class~$[\omega]$ admits a~contact structure $H = \ker \theta$, where $\theta \in \Omega^1(M)$ is a~connection 1-form whose curvature is~$\omega$ (cf.\ \cite[Section~7.2]{geiges_book} and \cite[Example~3.48]{mcduff_salamon}). Using Note~\eqref{item:po}, one can check that prequantum circle bundles are examples of CIR.
\end{Example}

\begin{Example}[contact non-commutative integrable systems, cf.\ \cite{jovanovic}]\label{sec:cont-isotr-real-2}
Given a contact manifold $\CM$, say that $X \in \X(M)$ is an {\em infinitesimal automorphism} of $\CM$ if it is a Reeb vector f\/ield (cf.\ Example~\ref{exm:ctc_jacobi}). Fix a contact manifold $\CM$ and an inf\/initesimal automorphism $X$ of $\CM$. Following \cite[Section~5.1]{jovanovic}, the dynamical system $\dot{x}(t) = X(x(t))$ is said to be {\em $($non-Hamiltonian$)$ completely integrable} if there exists an open, dense subset $M_{\mathrm{reg}} \subset M$, a proper surjective submersion $\phi\colon M_{\mathrm{reg}} \to P$, and an abelian Lie algebra~$\mathcal{X}$ of inf\/initesimal symmetries of $\CM$ such that $X$ is tangent to the f\/ibres of $\phi$ and the latter are the orbits of~$\mathcal{X}$. Moreover, such a system is said to be {\em contact non-commutative $($Hamiltonian$)$ integrable} if the foliation induced by~$\phi$ is transversal to $\mathcal{H}$ and satisf\/ies property~\ref{item:IR3} (cf.\ \cite[Def\/initions~2 and~3]{jovanovic}). Contact non-commutative ntegrable systems (with connected f\/ibres) can be viewed as examples of CIRs: \cite[Theorem~3]{jovanovic} implies that if $\phi \colon M_{\mathrm{reg}} \to P$ is the associated proper surjective submersion, then there exists a unique line bundle $L \to P$, a unique bundle map $F \colon \phi^*(L) \to L_M$, and a unique Jacobi structure $\{\cdot,\cdot\}$ making $\phi \colon \CM \to \JM$ satisfy property~\ref{item:IR1}.
\end{Example}

The existence of a CIR imposes geometric restrictions on the underlying Jacobi manifold, as illustrated below.

\begin{Lemma}\label{lemma:regular} If $\phi\colon \CM\to \JM$ is a CIR, then the Jacobi structure on $P$ is regular with corank equal to $k = \dim \ker D\phi$, and all its leaves are even dimensional.
\end{Lemma}

\begin{proof} Fix points $p \in P$ and $m \in M$ with $m \in \phi^{-1}(p)$. By property~\ref{item:IR1} and Note~\ref{obs:Jac_maps}, have that
\begin{gather}\label{contando_rangos}
\rank \rho_{p}\big(J_p^1(L)\big)=\rank \big(\rho_{M,m} \circ \phi^*\big(J^1L\big)\big) -\rank\ker \big(D\phi|_{\rho_{M,m}\circ\phi^*(J^1L)}\big).
\end{gather}
On the one hand, property \ref{item:IR2} implies that $\rho_M\circ\phi^* \colon J^1L \to T M$ is injective: indeed, if $c \colon H \times H \to L$ is the curvature map of equation~\eqref{curvature}, it can be shown that for $(\bar u,\bar\eta)\in\Gamma(L_M)\oplus \Omega^1(M;L_M)\simeq \Gamma(J^1L_M)$,
 \begin{gather*}
 \rho_M(\bar u,\bar\eta)=R_{\bar u}-c^\sharp(\bar\eta|_H).
 \end{gather*}
Hence, if $\rho_M(\phi^*u,\phi^*\eta)=0$, for $u\in\Gamma(L)$ and $\eta\in\Omega^1(P;L)$, then $R_{(\phi^*u)}(m)=c^\sharp(\circ(\phi^*\eta|_H))(m)$ $\in H$. By def\/inition of Reeb vector f\/ields, this means that
\begin{gather*}
 (\phi^*u)(m)=R_{(\phi^*u)}(m)\quad {\rm mod} \ H=0,
\end{gather*}
implying in turn that $\phi^*\eta_m|_H=0$. Condition~\ref{item:IR2} gives that $D\phi|_{H_m}\colon H_m\to T_{\phi(m)}P$ is surjective (as $\phi$ is a surjective submersion); thus $\phi^*\eta_m|_H=0$ holds if and only if $\eta_{\phi(m)}=0$, thus proving injectivity of $\rho_M \circ \phi^*$.

Condition \ref{item:IR3} gives that $\rank\ker (D\phi|_{\rho_M\circ\phi^*J^1L})=\rank \ker D\phi=k$. By def\/inition $\rank J^1L=2n-k+2$, while injectivity of $\rho_M \circ \phi^*$ implies that
\begin{gather*}
 \rank \rho_{M,m}(\phi^* J^1L) = 2n-k+2;
\end{gather*}
equation \eqref{contando_rangos} yields $\rank \rho_{P,p}(J_p^1(L)))=2n+2-k - k = 2n+2-2k$. Since $p \in P$ is arbitrary, the proof of the lemma is completed.
\end{proof}

\begin{Remark}\label{rk:mult_free}
Suppose that $G \curvearrowright \CM$ is multiplicity-free (cf.\ Def\/inition \ref{defn:mult_free}); by assumption~$G$ is compact, thus implying that its Lie algebra~$\mathfrak{g}$ is of compact type. Thus all symplectic leaves of the linear Poisson structure $\Lambda$ on~$\mathfrak{g}^*$ are compact, as these are precisely coadjoint orbits. Moreover, since they are contained in the level set of a~quadratic Casimir, they are transversal to the Euler vector f\/ield $E \in \X^1(\mathfrak{g}^* \setminus \{0\})$. By \cite[Remarque~2.6]{Dazord-Lich-Marle91}, this implies that all leaves of the Jacobi manifold $\left(\mathbb{P}(\mathfrak{g}^*), O(1), \{\cdot,\cdot\}\right)$ of Example~\ref{exm:proj_dual_lie_algebra} are even dimensional.
\end{Remark}

\begin{Note}\label{obs:prop_jac_even}
Suppose that $\JM$ is regular and that all of its leaves are even dimensional; denote the induced foliation by $\mathcal{F}$. Then
 \begin{itemize}\itemsep=0pt
 \item the kernel of the anchor of $J^1L \to P$ is a bundle of abelian Lie algebras f\/itting in a short exact sequence of vector bundles
 \begin{gather} \label{eq:38}
 0 \to \nu^* \otimes L \to \ker \rho \to L \to 0,
 \end{gather}
 where $\nu^* = \left( T \mathcal{F}\right)^{\circ}$ is the conormal bundle to $\mathcal{F}$;
 \item there is a canonical foliated 2-form $\omega_{\mathcal{F}} \in\Omega^2(\mathcal{F};L)$ uniquely def\/ined by
 \begin{gather} \label{eq:29}
 \omega_{\mathcal{F}}\big(\rho\big(j^1 u\big), \rho\big(j^1 v\big)\big) : = \{u,v\},
 \end{gather}
 for any $u,v \in \Gamma(L)$. An explicit form for $\omega$ is
 \begin{gather*}
 \omega_{\mathcal{F}}(\rho(u,\eta),\rho(v,\zeta)) : = \{u,v\} + \eta(\rho(v,\zeta)) - \zeta(\rho(u,\eta)),
 \end{gather*}
 for $(u,\eta), (v,\zeta) \in \Gamma(J^1 L)$;
\item the form $\omega_{\mathcal{F}}$ is closed with respect to the dif\/ferential $\de_{\mathcal{F}}$ on the complex $\Omega^* (\mathcal{F};L)$ associated to the induced f\/lat connection
 \begin{gather*}
 \bar{\nabla} \colon \ \Gamma(T \mathcal{F}) \times \Gamma(L) \to \Gamma(L)
 \end{gather*}
 on $L$ uniquely def\/ined by
 \begin{gather} \label{eq:26}
 \bar{\nabla}_{\rho(j^1u)} v := \nabla_{j^1 u}v = \{u,v\},
 \end{gather}
 for any $u,v \in \Gamma(L)$. Explicitly, the dif\/ferential $\de_{\mathcal{F}}$ is given on a foliated $k$-form $\omega$ by the standard Koszul-type formula
 \begin{gather*}
 \de_{\mathcal{F}}(\omega)(X_0,\ldots,X_k)= \sum_i(-1)^i\bar\nabla_{X_i}\big(\omega\big(X_0,\ldots,\hat X_i,\ldots, X_k\big)\big)\\
\hphantom{\de_{\mathcal{F}}(\omega)(X_0,\ldots,X_k)=}{} +\sum_{i<j}(-1)^{i+j}\omega\big([X_i,X_j], X_0,\ldots,\bar X_i,\ldots, \bar X_j,\ldots, X_k \big).
\end{gather*}
The foliated cohomology class $ [\omega_{\mathcal{F}}]\in H^2(\mathcal{F};L)$ plays a central role in the classif\/ication of CIR of a f\/ixed Jacobi manifold (cf.\ Theorem~\ref{thm:main}).
 \end{itemize}
Proofs of the above properties can be found in Appendix \ref{sec:proofs-prop-regul}.
\end{Note}

\begin{Example}\label{exm:Poisson_case}
The properties of regular Jacobi manifolds listed in Note~\ref{obs:prop_jac_even} resemble closely those satisf\/ied by Poisson manifolds; in fact, the former generalise the latter, as illustrated below. Let $\PM$ be a regular Poisson manifold, with symplectic foliation $\mathcal{F}$ and conormal bundle denoted by~$\nu^*$. The corresponding Jacobi structure is def\/ined on $L := \RR_{P}$ and the anchor $ \rho: =J^1P=T^*P\oplus \RR \to TP$ is given by $\rho(\eta, c) = \Lambda^{\sharp}(\eta)$ (cf.\ Example~\ref{obs:kir}). Thus
 \begin{gather*}
 \ker \rho = \nu^* \oplus \RR \subset \cotan P \oplus \RR,
 \end{gather*}
while the Spencer operator is given by $ \mathrm{D}(\eta,f) := \de f - \eta$, where $(\eta,f) \in \Omega^1(P) \oplus C^{\infty}(P)$ (cf.\ Example~\ref{exm:spencer_trivial}). Moreover, the $T\mathcal{F}$-connection given in Note~\ref{obs:prop_jac_even} reduces to the Bott connection
 \begin{gather*}
 \bar{\nabla} \colon \ \Gamma(T\mathcal{F}) \times C^{\infty}(P) \to
 C^{\infty}(P), \qquad (X,f) \mapsto \de f(X);
 \end{gather*}
thus the induced dif\/ferential $\de_{\mathcal{F}}$ is the restriction of the exterior derivative to the foliation~$\mathcal{F}$. The cohomology class
 \begin{gather*}
 [\omega_{\mathcal{F}}] \in \HH^2(\mathcal{F};L) = \HH^2(\mathcal{F})
\end{gather*}
def\/ined by Lemma~\ref{lemma:coho_class} is that of the foliated symplectic form def\/ined by $\Lambda$ (cf.\ Note~\ref{obs:poisson_coho}).
\end{Example}

\section{Classif\/ication of contact isotropic realisations}\label{sec:smooth-invar-cirs}
As a stepping stone towards classifying CIRs over a f\/ixed Jacobi manifold, some smooth invariants are constructed in analogy with~\cite{dd}. These are the period lattice and the Chern class. Then we proceed to study in more detail the structure that the period lattice has, giving rise to the notion of {\it transversal $\ZZ$-projective structures}. We f\/inish by providing a cohomological criteria for the existence of CIRs, which combines all the ingredients described before in this section

\subsection{The period lattice and the Chern class}

 We f\/irst show that any contact realisation of $\JM$ (cf.\ Def\/inition \ref{defn:CR}) comes equipped with an action of the Lie algebroid $J^1L \to P$. When the realisation is isotropic we called {\it period lattice} the isotropy of the action $\Sigma\subset J^1L$. Theorem~\ref{thm:local_model} shows that in this case, a CIR is in fact locally isomorphic to $\Ker\rho/\Sigma$, hence def\/ining a cohomology class, which we called its {\it Chern class}.

\begin{Definition}\label{defn:action_alg} An {\it action} of a Lie algebroid $\mathsf{A} \to P$ on the smooth manifold $M$ along a map $\phi \colon M \to P$ is a vector bundle map $\psi\colon \phi^*\mathsf{A} \to TM$ such that
 \begin{enumerate}[label=(A\arabic*), ref=(A\arabic*)]\itemsep=0pt
 \item \label{item:3} it induces a Lie algebra homomorphism $\Gamma(J^1L)\to \X(M)$;
 \item \label{item:4} for all $m \in M$, $D_m\phi\circ \psi_m=\rho_{\phi(m)}$, where $\rho \colon \mathsf{A} \to TP$ is the anchor map.
 \end{enumerate}
 If $\psi$ is injective, the action is said to be {\em faithful}.
\end{Definition}

\begin{Lemma}\label{lemma:action}
 Let $\phi\colon \CM\to\JM$ be a contact realisation. The map $\psi\colon \phi^*J^1L\to T M$ given at the level of sections by
 \begin{gather*}
 \phi^*j^1 u\mapsto R_{\phi^*u},\qquad \forall\, u\in \Gamma(L)
 \end{gather*}
defines a faithful action of $J^1L$ on $\phi \colon M \to P$.
\end{Lemma}

\begin{Note}\label{rk:action}
Note that $\psi$ can be alternatively described as the restriction of $\rho_M$ to $\phi^*J^1L \subset J^1L_M$.
\end{Note}

\begin{proof}[Proof of Lemma \ref{lemma:action}] Since $\phi$ is Jacobi, Note \ref{obs:Jac_maps} gives that
 \begin{gather}\label{action}
 D\phi\circ\rho_M\circ\phi^*=\phi^*(\rho).
 \end{gather}
 In light of Note \ref{rk:action}, it follows that property~\ref{item:4} holds. It remains to prove that $\psi$ induces a~Lie
 algebra morphism, i.e.,
 \begin{gather}\label{respect.bracket}
 \psi([\alpha,\alpha'])-[\psi(\alpha),\psi(\alpha')]=0,
 \end{gather}
for all $\alpha,\alpha' \in \Gamma(J^1L)$. This equation holds when $\alpha,\alpha'\in\Gamma(J^1L)$ are holonomic sections, i.e., of the form $j^1u,j^1v$ for $u,v \in \Gamma(L)$. This is because $\phi$ is a Jacobi map, which implies that
 \begin{gather*}
 R_{\phi^*\{u,v\}}=R_{\{\phi^*u,\phi^*v\}}=[R_{\phi^*u},R_{\phi^*v}].
 \end{gather*}
for $u,v\in \Gamma(L)$, where in the last equality property~\eqref{eq:2} of the Jacobi bracket of $\CM$ is used. In general, notice that the left hand side of equation~\eqref{respect.bracket} satisf\/ies
 \begin{gather*}
 \psi([f\alpha,\alpha'])-[\psi(f\alpha),\psi(\alpha')] = (\phi^*f)(\psi([\alpha,\alpha'])-[\psi(\alpha),\psi(\alpha')]);
 \end{gather*}
for $f \in C^{\infty}(P)$ and $\alpha,\alpha' \in \Gamma(J^1L)$. This again follows from the fact that $\phi$ is a Jacobi map, so that equation~\eqref{action} holds. Thus, since equation~\eqref{respect.bracket} holds for holonomic sections and these generate $\Gamma(J^1L)$ as a~$C^{\infty}(P)$-module, it follows that equation~\eqref{respect.bracket} holds for all sections of~$J^1L$, which proves that $\psi$ satisf\/ies property~\ref{item:3}. The proof of Lemma~\ref{lemma:regular} shows that $\psi$ is injective (this follows from property~\ref{item:IR2}).
\end{proof}

Henceforth, let $\phi \colon \CM \to \JM$ be a CIR unless otherwise stated. Denote the anchor of the Lie algebroid associated to $\JM$ by $\rho \colon J^1 L \to TP$.

\begin{Corollary}\label{obs:trans_inf_action}
The restriction of the action $\psi$ of Lemma~{\rm \ref{lemma:action}} to the bundle of abelian Lie algebras $\ker\rho \to P$ defines a~faithful action of $\ker \rho \to P$ on $M$ along $\phi \colon M \to P$ which induces a vector bundle isomorphism $\psi \colon \phi^*(\ker \rho) \to \ker D\phi$.
\end{Corollary}
\begin{proof} The only thing to check is that $\psi \colon \phi^*(\ker \rho) \to \ker D\phi$ is an isomorphism. Since $\psi$ is injective by Lemma~\ref{lemma:action}, it suf\/f\/ices to check that the two vector bundles have the same rank, but this follows from Lemma~\ref{lemma:regular}.
\end{proof}

The action $\psi \colon \phi^* \ker \rho \to T M$ should be thought of as being inf\/initesimal; since the f\/ibres of $\phi$ are compact, $\psi$ can be integrated to an action of $\pi \colon \ker \rho \to P$ considered as a bundle of abelian Lie groups. The integrated action is given by
\begin{gather}\label{eq:6}
 \Psi \colon \ \ker \rho \tensor[_\pi]{\times}{_{\phi}} M\to M, \qquad (\alpha,m) \mapsto \varphi^1_{\alpha}(m),
\end{gather}
where $ \ker \rho \tensor[_\pi]{\times}{_{\phi}} M :=\{ (\alpha,m) \in \ker \rho \times M \,|\,\pi(\alpha) = \phi(m) \}$ is a smooth manifold, and $\varphi^1_{\alpha}$ is the time-1 f\/low of $\psi(\phi^*\alpha)$.

\begin{Note}\label{obs:per_lattice} For each $p \in P$, the action of equation \eqref{eq:6} restricts to an action of the abelian Lie group $\ker \rho_{p} \cong \RR^k$ on $\phi^{-1}(p)$. Connectivity of $\phi^{-1}(p)$ and Corollary~\ref{obs:trans_inf_action} imply that the action of $\ker \rho_{p}$ is transitive. Moreover, since $\ker \rho_{p}$ is an abelian Lie group, the isotropy of the action at any two points on $\phi^{-1}(p)$ are equal. Therefore, for
\begin{gather*}
 \Sigma_p := \big\{ \alpha\in \ker \rho_{p} \,|\, \varphi^1_{\alpha} = \mathrm{id} \big\},
\end{gather*}
the isotropy subgroup at $p$, then
\begin{gather*}
 \phi^{-1}(p) \cong \ker \rho_{p}/\Sigma_p;
\end{gather*}
since $\phi^{-1}(p)$ is compact by assumption, it follows that $\phi^{-1}(p)$ is dif\/feomorphic to $\mathbb{T}^k$ and that~$\Sigma_p$ is a cocompact lattice in $\ker \rho_{p}$ and, therefore, isomorphic to $\ZZ^k$.
\end{Note}

Just as in the theory of symplectic isotropic realisations of Poisson manifolds (cf.\ Section~\ref{sec:sympl-isotr-real}), the isotropy of the action \eqref{eq:6} plays an important role in the classif\/ication of contact isotropic realisations of Jacobi manifolds and endows the foliation of the underlying Jacobi manifold with an interesting transverse geometric structure (cf.\ Section~\ref{sec:period-lattices-as}).

\begin{Definition}\label{defn:period_lattice} The subset
 \begin{gather*}
 \Sigma := \coprod\limits_{p \in P} \Sigma_p \subset \ker \rho_{p}
 \end{gather*}
is called the {\it period lattice} of the CIR $\phi \colon \CM \to \JM$.
\end{Definition}

The following theorem provides a smooth local model for a CIR; its proof is omitted as it is entirely analogous to that of \cite[Theorem~2.1]{dd}.

\begin{Theorem}\label{thm:local_model}
 Let $\phi\colon \CM\to\JM$ be a CIR with associated period lattice $\Sigma$. Then
 \begin{enumerate}[label={\rm (\Roman*)}, ref={\rm (\Roman*)}]\itemsep=0pt
 \item \label{item:1} $\Sigma$ is a closed submanifold of $\ker\rho$ with the property that the composite $\Sigma \hookrightarrow \ker \rho \to P$ is a $\ZZ^k$-bundle, i.e., it has fibre $\ZZ^k$ and structure group $\mathrm{GL}(k;\ZZ)$;
 \item \label{item:2} the quotient $\ker\rho/\Sigma$ is a smooth manifold and the projection
 \begin{gather*}
 \pi \colon \ \ker\rho/\Sigma \to P
 \end{gather*}
is a fibre bundle with fibre $\mathbb{T}^k$;
 \item \label{item:7} upon a choice of a local section $\sigma\colon U\subset M\to \phi^{-1}(U)$, the map
 \begin{gather}\label{eqn:local_isomorphism}
 \Psi_\sigma\colon \ \ker \rho/\Sigma|_U\to \phi^{-1}(U), \qquad [\alpha]\mapsto \varphi^1_{\alpha}(\sigma(\pi(\alpha)))
 \end{gather}
 is an isomorphism of fibres bundles, making the following diagram
 \begin{gather*}
 \xymatrix{
 \ker \rho/\Sigma|_U \ar[rr]^{\Psi_\sigma} \ar[rd]_{\pi} & & \phi^{-1}(U) \ar[ld]^{\phi}\\
 & U & }
 \end{gather*}
 commute;
 \item \label{item:18} the map $\phi \colon M \to P$ is a principal $\ker \rho/\Sigma$-bundle and is classified up to isomorphism by a~cohomology class $t \in \HH^1(P;\mathcal{C}^{\infty}(\ker \rho/\Sigma))$, where $\mathcal{C}^{\infty}(\ker \rho/\Sigma)$ denotes the sheaf of smooth sections of $\ker
 \rho/\Sigma \to P$.
 \end{enumerate}
\end{Theorem}

The short exact sequence of sheaves
\begin{gather*}
 1 \to \underline{\Sigma} \hookrightarrow \mathcal{C}^{\infty} (\ker \rho) \to \mathcal{C}^{\infty}(\ker \rho/\Sigma) \to 1,
\end{gather*}
where $ \underline{\Sigma}$ denotes the sheaf of sections of $\Sigma \to P$, induces a long exact sequence in cohomology whose connecting morphisms
\begin{gather*}
 \delta\colon \ \HH^{l}(P; \mathcal{C}^{\infty}(\ker \rho/\Sigma)) \to \HH^{l+1}(P;\underline{\Sigma})
\end{gather*}
are isomorphisms for all $l \geq 1$, since $\mathcal{C}^{\infty}(\ker \rho)$ is f\/ine.

\begin{Definition}\label{defn:cc}
 The cohomology class
 \begin{gather*}
 c = \delta(t) \in \HH^2(P;\underline{\Sigma})
 \end{gather*}
is called the {\it Chern class} of the CIR $\phi \colon \CM \to \JM$.
\end{Definition}

\subsection{Period lattices as transversal $\ZZ$-projective structures}\label{sec:period-lattices-as}
In order to tackle the problem of determining whether a given regular Jacobi mani\-fold $\JM\!$ all of whose leaves are even dimensional (see Lemma~\eqref{lemma:regular}) admits a CIR, it is fundamental to describe necessary conditions for a smooth submanifold $\Sigma \subset \ker \rho$ to be the period lattice of some CIR over~$\JM$. This is the aim of this subsection, which, in fact, also gives a description of these necessary conditions in terms of a geometric structure {\em transverse} to the foliation.

First, local sections of a period lattice of a CIR are characterised in the following lemma, whose statement is analogous to the corresponding condition to be a period lattice of a symplectic isotropic realisation of a Poisson manifold (cf.\ \cite[Corollaire~1]{dd}).

\begin{Lemma}\label{cor:sections_period}
If $\Sigma$ is the period lattice of $\phi \colon \CM \to \JM$ then any local section $\alpha \in \Gamma_{\mathrm{loc}}(\Sigma)$ is of the form $j^1 u$ for some $u \in \Gamma_{\mathrm{loc}}(L)$.
\end{Lemma}

For expositional reasons, the proof of Lemma \ref{cor:sections_period} is postponed until Appendix~\ref{apendixC}. Therefore, the period lattice of a CIR $\phi \colon \CM \to \JM$ is a full-rank lattice of $\ker \rho$ (by Part~\ref{item:1} of Theorem~\ref{thm:local_model}) whose local sections are holonomic (by Lemma~\ref{cor:sections_period}). Recalling that $\ker \rho = \Sigma \otimes_{\ZZ} \RR=: \Sigma^{\RR}$ f\/its in the short exact sequence of equation~\eqref{eq:38}, it follows that any period lattice of a CIR is an example of the following object.

\begin{Definition}\label{defn:tilp}
Let $(N,\mathcal{F})$ be a foliated manifold whose foliation has codimension $l \geq 0$ and f\/ix a line bundle $L \to N$. A {\it transversal $\ZZ$-projective lattice} on $(N,\mathcal{F})$ is a choice of embedded smooth submanifold $\Sigma\subset J^1L$ satisfying:
 \begin{enumerate}[label=(T\arabic*), ref=(T\arabic*)]\itemsep=0pt
 \item\label{item:8} the composite
 \begin{gather*}
 \Sigma\hookrightarrow J^1L\to N
 \end{gather*}
 is a $\ZZ^{l+1}$-bundle (cf.\ Part \ref{item:1} of Theorem~\ref{thm:local_model});
 \item \label{item:11} $\Sigma^{\RR} := \Sigma \otimes_{\ZZ} \RR$ f\/its in a short exact sequence of vector bundles
 \begin{gather*}
 0 \to \nu^* \otimes L \to \Sigma^{\RR} \to L \to 0,
 \end{gather*}
where $\nu^* \subset \cotan N$ is the conormal bundle to $\mathcal{F}$, and $\Sigma^{\RR} \to L$ is the restriction of the canonical projection $J^1 L \to L$;
 \item\label{item:9} any $\alpha \in \Gamma_{\mathrm{loc}}(\Sigma)$ is holonomic.
 \end{enumerate}
\end{Definition}

\begin{Example}\label{exm:trivial_bundle} If $L = \RR_N$, then $J^1\RR_N$ is canonically isomorphic to $T^*N \oplus \RR$. In this case conditions~\ref{item:8} and~\ref{item:11} become that $\Sigma$ is a full-rank lattice in $\nu^* \oplus \RR$, where $\nu^* :=\left(T\mathcal{F}\right)^{\circ}\subset T^*N$ is the conormal bundle to $\mathcal{F}$.
\end{Example}

The reason for the terminology in Def\/inition \ref{defn:tilp} is that a~transversal $\ZZ$-projective lattice induces a transversal $\ZZ$-projective structure in the following sense (cf.\ Proposition \ref{thm:tilp_tilpa}).

\begin{Definition}\label{defn:tipa} A {\it transversal $\ZZ$-projective structure} on a foliated manifold $(N,\mathcal{F})$ is an atlas $\mathcal{A}
 :=\{(U_i, \chi_i)\}$ of submersions $\chi_i \colon U_i \to \RR\mathrm{P}^l$ locally def\/ining $\mathcal{F}$, such that, for all~$i$,~$j$ with $U_{ij} := U_i \cap U_j \neq \varnothing$, there exists a smooth map $A_{ij}\colon U_{ij} \to \mathrm{GL}(l+1;\ZZ)$ satisfying
 \begin{itemize}\itemsep=0pt
 \item $\chi_j = [A_{ij}]\circ \chi_i$ on $U_{ij}$, where $[A_{ij}] \colon \RR \mathrm{P}^l \to \RR \mathrm{P}^l$ is the map induced by $A_{ij}$;
 \item $\{A_{ij}\}$ satisf\/ies the cocycle condition.
 \end{itemize}
The components of $\chi_j$ are called {\em local transversal $\ZZ$-projective coordinates} on $(N,\mathcal{F})$.
\end{Definition}

\begin{Note}\label{obs:trans_proj_structure}
Transversal $\ZZ$-projective structures are examples of {\it transversal projective structures} on foliated manifolds, which can be def\/ined as follows. If $(N,\mathcal{F})$ is a foliated manifold and $l \geq 0$ is the codimension of $\mathcal{F}$, then a~transversal projective structure on $(N,\mathcal{F})$ is an atlas $\mathcal{A} :=\{(U_i, \chi_i)\}$ of submersions $\chi_i \colon U_i \to \RR\mathrm{P}^l$ locally def\/ining $\mathcal{F}$, such that, for all~$i$,~$j$ with $U_{ij} := U_i \cap U_j \neq \varnothing$, there exists a smooth map $\mathsf{a}_{ij}\colon U_{ij} \to \mathrm{PGL}(l+1;\RR)$ satisfying $\chi_j = \mathsf{a}_{ij}\circ \chi_i$ on $U_{ij}$.
 \end{Note}

\begin{Proposition}\label{thm:tilp_tilpa}
 There is a {\rm 1-1} correspondence between transversal $\ZZ$-projective lattices and transversal $\ZZ$-projective structures on a foliated manifold $(N,\mathcal{F})$. In this {\rm 1-1} correspondence, locally
 \begin{gather*}
 \Sigma := \ZZ \big\langle j^1x_1,\ldots, j^1x_{l+1} \big\rangle,
 \end{gather*}
where $[x_1:\ldots:x_{l+1}]$ are local transversal $\ZZ$-projective
 coordinates.
\end{Proposition}

The proof of Proposition \ref{thm:tilp_tilpa} is given in Appendix~\ref{apendixC}.

It is instructive to see the above correspondence in a specif\/ic example, which can be considered to be `universal'.

\begin{Example}[a $\ZZ$-projective lattice and structure on $\mathbb{R}\mathrm{P}^l$]\label{ex:projective_plane}
Fix $l \geq 0$ and consider the line bundle $\pi\colon O(1) \to \RR \mathrm{P}^l$. The standard coordinates $x_1,\ldots,x_{l+1}$ on $\RR^{l+1}$ def\/ine global sections of $O(1)$, whose f\/irst jets are independent everywhere. In other words, $J^1(O(1)) \cong \RR \langle j^1x_1,\ldots,j^1x_{l+1} \rangle$. Def\/ine
\begin{gather*}
 \Sigma^{l} : = \ZZ \big\langle j^1x_1,\ldots,j^1x_{l+1} \big\rangle \subset J^1(O(1));
 \end{gather*}
it can be verif\/ied directly that $\Sigma^l$ def\/ines a~$\ZZ$-projective lattice on $\RR\mathrm{P}^l$ with respect to the foliation by points\footnote{Whenever the foliation is by points, the adjective `transversal' is omitted.}. On the other hand, there is a canonical $\ZZ$-projective structure on $\RR\mathrm{P}^l$, constructed by taking the identity $\RR\mathrm{P}^l \to \RR \mathrm{P}^l$ as the coordinate map and the maps $A_{ij} \in \mathrm{GL}(l+1;\ZZ)$ also equal to the identity. The above lattice and structure are mapped one to the other under the correspondence of Proposition~\ref{thm:tilp_tilpa}.
\end{Example}

\subsection{The realisation problem for CIRs}\label{sec:class-cir-over}
The aim of this section is to solve the {\em realisability} problem for CIRs, outlined below. Given a~regular Jacobi manifold $\JM$ all of whose leaves are even dimensional, let $\mathcal{F}$ denote the induced foliation on~$P$, and f\/ix a transversal $\ZZ$-projective lattice $\Sigma \subset \ker \rho$, where $\rho$ is the anchor map of $J^1L$ (cf.\ Def\/inition~\ref{defn:tilp}). Natural questions to address are
\begin{itemize}\itemsep=0pt
\item to determine which cohomology classes in $\HH^2(P;\Sigma)$ correspond to Chern classes of CIR of $\JM$ whose period lattice is $\Sigma$, and
\item supposing that $c \in \HH^2(P;\Sigma)$ is the Chern class of a CIR with period lattice $\Sigma$, to determine all the CIRs with period lattice $\Sigma$ and Chern class $c$ up to the following notion of isomorphism.
\end{itemize}

\begin{Definition}\label{defn:iso_cir}
Two CIRs $\phi \colon \CM \to \JM$ and $\phi'\colon (M',H') \to \JM$ are said to be {\em isomorphic} if there exists a dif\/feomorphism $I\colon M \to M'$ satisfying $\phi'\circ I = \phi$ and $I^* \theta' = \theta$, where $\theta$ and $\theta'$ are the generalised contact forms associated to~$H$ and~$H'$ respectively.
\end{Definition}

\begin{Remark}\label{rk:strict} Def\/inition \ref{defn:iso_cir} provides a notion of {\em strict} isomorphism of CIRs, i.e., the bundle component of the Jacobi morphism $I \colon \CM \to (M',H')$ is the identity. It is possible to def\/ine a more general notion of isomorphism of CIRs which allows for arbitrary bundle components, and to prove a result analogous to Part~\ref{item:20} of Theorem~\ref{thm:main} for this more general notion of isomorphism using the ideas presented below. For simplicity, only the case of strict isomorphism is considered here.
\end{Remark}

The main result of this section, Theorem \ref{thm:main}, provides cohomological criteria that solve the above problems. Recall that, associated to a Jacobi manifold $\JM$ as above, there is a canonical foliated, closed 2-form $\omega_{\mathcal{F}}$ with values in~$L$ (cf.\ Note~\ref{obs:prop_jac_even}). On the other hand, it can be shown that the Spencer operator $\mathrm{D} \colon \Gamma (J^1L) \to \Omega^1\left(P;L\right)$ associated to $\JM$ induces maps in cohomology
\begin{gather} \label{eq:39}
 \mathcal{D} \colon \ \HH^l\big(P;\mathcal{C}^{\infty}(\ker \rho/\Sigma)\big)\to \HH^l\big(P,\mathcal{Z}^1(\mathcal{F};L)\big)
\end{gather}
for all $l \geq 0$ (cf.\ the discussion following Corollary \ref{cor:zero}). For $l \geq 1$, there are isomorphisms $\HH^l(P;\mathcal{C}^{\infty}(\ker \rho/\Sigma)) \cong \HH^{l+1}(P;\underline{\Sigma})$ (cf.\ Section~\ref{sec:smooth-invar-cirs}), and $\HH^l(P,\mathcal{Z}^1(\mathcal{F};L)) \cong
\HH^{l+1}(\mathcal{F};L)$ via a standard double \v{C}ech--Lie algebroid dif\/ferential complex argument. Hence, for all $l \geq 1$, the homomorphism of equation~\eqref{eq:39} induces a homomorphism
\begin{gather*}
 \HH^{l+1}(P;\underline{\Sigma})\to \HH^{l+1}(\mathcal{F};L),
\end{gather*}
which, by abuse of notation, is also denoted by $\mathcal{D}$.

\begin{Theorem}\label{thm:main}\quad
 \begin{enumerate}[label={\rm \arabic*.}, ref={\rm (\arabic*)}]\itemsep=0pt
 \item \label{item:19} A Jacobi manifold $(P,L,\{,\})$ admits a CIR with period lattice $\Sigma$ and Chern class $c \in \HH^2(P;\Sigma)$ if and only if
 \begin{gather*} 
 \mathcal{D}c = [\omega_{\mathcal{F}}].
 \end{gather*}
 \item \label{item:20} The CIRs over $\JM$ with period lattice $\Sigma$ and Chern class $c$ are classified up to isomorphism by~$\HH^1(\mathcal{F};L)$. The correspondence is given as follows: if $\phi \colon \CM \to \JM$ is such a CIR, $\theta$ is the generalised contact form defining $H$, and $[\eta_{\mathcal{F}}] \in \HH^1(\mathcal{F};L)$, then any other CIR over~$\JM$ with period lattice $\Sigma$ and Chern class $c$ is isomorphic to
 \begin{gather*}
 \phi \colon \ (M, H' = \ker (\theta + \phi^* \eta) ) \to \JM,
 \end{gather*}
 where $\eta \in \Omega^1(P;L)$ is a 1-form which represents $[\eta_{\mathcal{F}}]$.
\end{enumerate}
\end{Theorem}

An immediate corollary of Theorem \ref{thm:main} is the classif\/ication of CIR for the zero Jacobi bracket on the trivial line bundle, obtained in~\cite{Banyaga,lerman_contact_toric,wolbert} (the last two deal with singularities coming from the existence of a global contact toric action). Given a~manifold~$P$, let $(P,0,\RR_{P})$ denote the Jacobi manifold corresponding to the zero Jacobi structure on the trivial line bundle.

\begin{Corollary}\label{cor:zero}
 Given any $\ZZ$-projective lattice $\Sigma \subset J^1P$ on the Jacobi manifold $(P,0,\RR_P)$, any cohomology class $c \in \HH^2(P;\underline{\Sigma})$ corresponds to a unique CIR of $(P,0,\RR_{P})$ with period lat\-ti\-ce~$\Sigma$ up to isomorphism.
\end{Corollary}

\begin{proof} Fix $(P,0,\RR_{P})$ and a $\ZZ$-projective lattice $\Sigma \subset J^1 P = \ker \rho$, $\rho \colon J^1 P \to TP$ being the anchor map (which is identically zero in this case!). The induced foliation $\mathcal{F}$ is by points and, therefore, the sheaf $\mathcal{Z}^1(\mathcal{F}):= \mathcal{Z}^1(\mathcal{F};\RR) = \Omega^1(P)$ is f\/ine, thus showing, in particular, that $\HH^2(\mathcal{F}) = 0$. Therefore the homomorphism $\mathcal{D} \colon \HH^2(P;\underline{\Sigma}) \to \HH^2(\mathcal{F})$ is trivial and every cohomology class can be realised as the Chern class of some CIR of~$(P,0,\RR_{P})$ with period lattice~$\Sigma$. To prove uniqueness up to isomorphism, observe that, since the foliation is by points, every closed foliated 1-form is exact.
\end{proof}

The rest of this subsection provides the necessary preparatory material as well as the proof of Theorem~\ref{thm:main}, which is at the end of the subsection. Some of the intermediate results are interesting on their own; for instance, Theorem~\ref{thm:aa} gives a local normal form for the contact structure of a CIR, which generalises \cite[Theorem~1]{Banyaga} and \cite[Theorems~4 and~5]{jovanovic}. Henceforth, f\/ix a regular Jacobi manifold $\JM$ all of whose leaves are even dimensional, let $\mathcal{F}$ denote the induced foliation on $P$, and f\/ix a transversal $\ZZ$-projective lattice
$\Sigma \subset \ker \rho$, where $\rho \colon J^1 L \to TP$ is the anchor map. First, the homomorphism $\mathcal{D}$ is constructed starting from the Spencer operator $\mathrm{D} \colon \Gamma(J^1L) \to \Omega^1(M;L)$ associated to the Lie algebroid $J^1L \to P$ (cf.\ Note~\ref{obs:prop_cso} and the ensuing discussion). Recall that there is a foliated f\/lat connection $\bar\nabla$ on~$L$ uniquely def\/ined by equation~\eqref{eq:26} (cf.\ Note~\ref{obs:prop_jac_even} and Lemma~\ref{prop:well-defined}).

\begin{Note}\label{obs:sod}
For any $\alpha \in \Gamma(\ker \rho)$, $\de_{\mathcal{F}}\mathrm{D}(\alpha) = 0$; in fact $\mathrm{D}(\alpha) = \de_{\mathcal{F}}(\pr(\alpha))$, where $\pr \colon \ker \rho \to L$ is the projection (cf.\ Note~\ref{obs:prop_jac_even} and equation~\eqref{eq:38}). The proof of this result can be found in Appendix~\ref{sec:proofs-prop-regul}.
\end{Note}

By Note \ref{obs:sod}, the Spencer operator induces a sheaf homomorphism $\bar{\mathrm{D}} \colon \mathcal{C}^{\infty}(\ker \rho) \to \mathcal{Z}^1(\mathcal{F};L)$, where $\mathcal{Z}^1(\mathcal{F};L)$ denotes the sheaf of closed, foliated forms with values in~$L$. If $\mathcal{C}^{\infty}(\ker \rho/\Sigma)$ denotes the sheaf of smooth sections of $\ker \rho/\Sigma \to P$, Lemma~\ref{cor:sections_period} implies that
there is a morphism of sheaves
\begin{gather} \label{eq:20}
 \hat{\mathrm{D}} \colon \ \mathcal{C}^{\infty}(\ker \rho/\Sigma) \to \mathcal{Z}^{1}(\mathcal{F};L),
\end{gather}
which induces the homomorphisms at the level of cohomology of equation~\eqref{eq:39}. The strategy of the proof of Theorem~\ref{thm:main} hinges upon a~strengthened version of Theorem~\ref{thm:local_model}, i.e., on a~{\em contact} local normal form for suf\/f\/iciently small open $\phi$-saturated neighbourhoods. One of the main ingredients of this result is the following.

\begin{Note}\label{obs:contact form} Consider the canonical contact form $\theta_{\mathrm{can}}\in \Omega^1(J^1 L; \pr^* L)$ described in Example~\ref{ex:jet_bundle}. Translations by elements of $\Sigma$, which are holonomic by Lemma~\ref{cor:sections_period}, preserve $\theta_{\mathrm{can}}$ (cf.\ Note~\ref{obs:prop_cso}), and therefore its restriction to $\ker\rho\subset J^1L$ descends to a 1-form
 \begin{gather*}
 \theta_0\in \Omega^1(\ker\rho/\Sigma;\pi^*L),
 \end{gather*}
which does not necessarily def\/ine a contact distribution (unless $\ker\rho=J^1 L$).
\end{Note}

Fix a CIR $\phi\colon \CM \to (P,L,\{,\})$ and let $\theta \in \Omega^1(M,L_M)$ denote the generalised contact form associated to~$\CM$.

\begin{Theorem}\label{thm:aa} Given a local section $\sigma\colon U\subset P\to \CM$, then
 \begin{gather*}
 \Psi_\sigma^*\theta=\theta_0+\pi^*\sigma^*\theta,
 \end{gather*}
where $\Psi_\sigma\colon \ker \rho/\Sigma|_U\to \phi^{-1}(U)$ is defined by equation~\eqref{eqn:local_isomorphism}. Moreover, $\de_{\mathcal{F}} (\sigma^*\theta ) = \omega_{\mathcal{F}}$.
\end{Theorem}

The proof of Theorem \ref{thm:aa} is postponed to Appendix~\ref{apendixC}. The following result (stated without proof as it can be checked using local coordinates) provides a partial converse to Theorem~\ref{thm:aa}.

\begin{Proposition}\label{prop:local_cir}
Let $U \subset P$ be an open set and fix $\beta\in\Omega^1(U,L)$ satisfying $\de_{\mathcal{F}}\beta=\omega_{\mathcal{F}}$. Then
\begin{gather*}
 \pi\colon \ (\ker\rho/\Sigma|_U,\theta_0+\pi^*\beta)\to (U,L_U,\{\cdot,\cdot\})
\end{gather*}
is a contact isotropic realisation.
\end{Proposition}

Finally, it is possible to proceed with the proof of the main result of this section.

\begin{proof}[Proof of Theorem~\ref{thm:main}] The strategy is similar to that of \cite[Theorems~4.2 and~4.3]{dd}. Throu\-ghout the proof, let $\mathcal{U} := \{U_i\}$ denote a good open cover of~$P$.

{\bf Part \ref{item:19}.} Suppose f\/irst that $\phi\colon \CM \to \JM$ is a CIR with Chern class $c \in \HH^2(P;\Sigma)$, and let $\theta$ be the generalised contact form associated to $\CM$. For each $i$, let $\sigma_i \colon U_i \to M$ be a local section of $\phi$; as in Theorem~\ref{thm:local_model}, the section $\sigma_i$ induces a local trivialisation $\Psi_i \colon \ker \rho/\Sigma|_{U_i} \to \phi^{-1}(U_i)$. By Def\/inition~\ref{defn:cc}, a \v{C}ech cocycle representing $c$ is given by the smooth maps $t_{ij} \colon U_{ij} \to \ker \rho /\Sigma$ def\/ined by
 \begin{gather*}
 \Psi_j^{-1} \circ \Psi_i([\alpha]) = [\alpha] + t_{ij}(\pi(\alpha)),
 \end{gather*}
where $\pi\colon \ker \rho/\Sigma \to P$ is the projection (cf.\ the discussion following Theorem~\ref{thm:local_model}). Thus a~\v{C}ech cocycle representing $\mathcal{D}c$ is given by $\{t_{ij}^*\theta_0 = \hat{\mathrm{D}}(t_{ij})\}$~-- cf.\ Note~\ref{obs:contact form}. By def\/inition of~$t_{ij}$, $\Psi_i \circ t_{ij} = \sigma_j$ on $U_{ij}$; thus
 \begin{gather} \label{eq:13}
 \sigma_j^* \theta =t_{ij}^* \circ \Psi_i^*\theta = t_{ij}^*(\theta_0 + \pi^* \sigma_i^*\theta) = t_{ij}^*\theta_0 + \sigma_i^* \theta,
 \end{gather}
where the second equality uses Theorem~\ref{thm:aa}, and the last follows by noticing that $\pi \circ t_{ij} = {\rm id}$. Equation~\eqref{eq:13} gives that $t_{ij}^*\theta_0 = \sigma_j^*\theta - \sigma_i^*\theta$. By Theorem~\ref{thm:aa}, $\de_{\mathcal{F}}(\sigma_j^* \theta) = \omega = \de_{\mathcal{F}}(\sigma_i^*\theta)$; this implies that the cohomology class of $\{\sigma_j^*\theta - \sigma_i^* \theta\}$ equals $[\omega_{\mathcal{F}}]$, as it is the dif\/ference of two primitives of~$\omega_{\mathcal{F}}$.

Conversely, suppose that $\mathcal{D}c = [\omega_{\mathcal{F}}]$. Choose the good open cover $\mathcal{U}$ so that, for each $i$, there exists $\beta_i \in \Omega^1(\mathcal{F}|_{U_i};L)$ with $\de_{\mathcal{F}}\beta_i = \omega_{\mathcal{F}}|_{U_i}$. Since $\mathcal{F}$ is regular, the foliated 1-forms~$\beta_i$ can be extended to elements of $\Omega^1(U_i;L)$, which are henceforth also denoted by~$\beta_i$ by abuse of notation. Since $\mathcal{D}c = [\omega_{\mathcal{F}}]$, there exists a~\v{C}ech cocycle $\{t_{ij}\}$ representing~$c$ which satisf\/ies $t_{ij}^* \theta_0 = \beta_j - \beta_i$ on~$U_{ij}$. For each $i$, consider the CIR $\pi \colon (\ker \rho/\Sigma|_{U_i},\theta_0 + \pi^*\beta_i) \to (U_i,L|_{U_i}, \{\cdot,\cdot\}|_{U_i})$ (cf.\ Proposition~\ref{prop:local_cir}). The principal $\ker \rho/\Sigma$-bundle over $P$ with Chern class $c$ is constructed (up to isomorphism) by glueing the above local models using the translations $t_{ij}$. The condition $t_{ij}^* \theta_0 = \beta_j - \beta_i$ ensures that the locally def\/ined generalised contact forms $\theta_0 + \pi^*\beta_i$ patch together to give a globally def\/ined generalised contact form $\theta \in \Omega^1(M;\phi^* L)$, where $\phi \colon M \to P$ is the principal $\ker \rho/\Sigma$-bundle over $P$ with Chern class $c$ constructed above. Setting $H := \ker \theta$, have that $\phi\colon \CM \to \JM$ is a CIR as it suf\/f\/ices to check that the properties~\ref{item:IR1},~\ref{item:IR2} and~\ref{item:IR3} hold locally and they do by construction, since the bundle is obtained by glueing CIR.

{\bf Part \ref{item:20}.} Suppose that $\mathcal{D}c = [\omega_{\mathcal{F}}]$; by Part \ref{item:19}, there exists a CIR $\phi \colon \CM \to \JM$ with period lattice $\Sigma$ and Chern class $c$, which is henceforth f\/ixed as a~`reference'. Suppose that $\phi'\colon (M',H') \to \JM$ is another CIR with period lattice $\Sigma$ and Chern class $c$. By composing with an isomorphism of principal $\ker \rho/\Sigma$-bundles, it may be assumed that $\phi' = \phi$ and that $M' = M$. Denote the generalised contact forms associated to $H$ and $H'$ by $\theta$ and~$\theta'$ respectively. For each $i$, consider the form $\eta_i:=\sigma^*_i\left(\theta - \theta'\right)$. By Theorem~\ref{thm:aa}, $\de_{\mathcal{F}} \eta_i = 0$. Moreover, on $U_{ij}$,
 \begin{gather*}
 \eta_i - \eta_j = \sigma^*_i (\theta - \theta' ) - \sigma^*_j(\theta - \theta')= (\sigma^*_i \theta - \sigma^*_j \theta) - (\sigma^*_i \theta' -
 \sigma^*_j \theta') = t_{ij}^* \theta_0 - t_{ij}^* \theta_0 = 0;
\end{gather*}
thus the collection $\{\eta_i\}$ def\/ines a globally def\/ined foliated closed 1-form $\eta_{\mathcal{F}} = \eta_{\mathcal{F}}(\theta,\theta')\in \HH^0(P;\mathcal{Z}^1(\mathcal{F};L)) = \mathrm{Z}^1(\mathcal{F};L)$. In other words, there is a~function (depending on the reference CIR $\phi \colon \CM \to \JM$!)
 \begin{gather*}
 b_{\theta}\colon \ \{ \text{CIRs with period lattice } \Sigma \text{ and Chern class } c\} \to \mathrm{Z}^1(\mathcal{F};L), \nonumber \\
 \left( \phi' \colon \ (M',H') \to \JM\right) \mapsto \eta_{\mathcal{F}}(\theta,\theta'),
\end{gather*}
where $\mathrm{Z}^1(\mathcal{F};L)$ denotes the vector space of globally def\/ined closed, foliated 1-forms with values in~$L$. Next it is shown that $b_{\theta}$ is surjective. Let $\eta_{\mathcal{F}} \in \mathrm{Z}^1(\mathcal{F};L)$ and let $\eta \in \Omega^1(P;L)$ be any 1-form extending it; then Proposition~\ref{prop:local_cir} ensures that $\phi\colon (M, H') \to \JM$ is a CIR with period lattice $\Sigma$ and Chern class $c$, where $H' = \ker \left(\theta + \phi^*\eta\right)$. Setting $\theta' = \theta + \phi^* \eta$, it is simple to check that $\eta_{\mathcal{F}}(\theta,\theta') = \eta_{\mathcal{F}}$, thus showing that $b_{\theta}$ is onto. However, $b_{\theta}$ does not descend to a function on the set of isomorphism classes of CIRs with period lattice $\Sigma$ and Chern class $c$. For, if $\phi \colon \CM \to \JM$ and $\phi \colon (M,H') \to \JM$ are isomorphic in the sense of Def\/inition~\ref{defn:iso_cir}, then $\eta_{\mathcal{F}}(\theta,\theta') \in \hat{\mathrm{D}}\left(\Gamma\left( \ker \rho/\Sigma\right)\right)$; this can be proved as follows. Using Theorem~\ref{thm:aa}, for each $i$ there exists an isomorphism of CIRs
 \begin{gather*}
 T_i \colon \ (\ker \rho/\Sigma|_{U_i},\, \theta_0 + \pi^* \sigma^*_i \theta ) \to (\ker \rho/\Sigma|_{U_i}, \, \theta_0 + \pi^* \sigma^*_i
 \theta' ).
 \end{gather*}
Fix $i$. Since an isomorphism of CIRs is, in particular, an isomorphism of principal $\ker \rho/\Sigma$-bundles, it follows that there exists a~section $\tau_i \in \Gamma (\ker \rho/\Sigma|_{U_i} )$ such that $T_i(\alpha + \Sigma) = \alpha + \tau_i +\Sigma$. Using the fact that $T_i^* (\theta_0 +\pi^*\sigma^*_i \theta' ) =\theta_0 +\pi^*\sigma^*_i \theta'$, i.e., $T_i$ is an isomorphism of CIRs, it follows that $\hat{\mathrm{D}}(\tau_i) = \eta_i$, as $\tau_i^* \theta_0 = \hat{\mathrm{D}}(\tau_i)$. This holds for all $i$; however, since the collection of isomorphisms $\{T_i\}$ patches together to give a {\em global} isomorphism of CIRs, it follows that for all~$i$,~$j$ with $U_{ij} \neq \varnothing$, $T_j \circ t_{ij} = t_{ij} \circ T_i$, which is equivalent to $\tau_i = \tau_j$ on~$U_{ij}$. This shows that the collection $\{\tau_i\}$ def\/ines a~globally def\/ined section $\tau \in \HH^0(P;\mathcal{C}^{\infty}(\ker \rho/\Sigma)) = \Gamma(\ker \rho/\Sigma)$; since, for all $i$, $\hat{\mathrm{D}}(\tau_i) = \eta_i$, it follows that the 1-form $\eta_{\mathcal{F}}(\theta,\theta') \in \hat{\mathrm{D}} (\Gamma (\ker \rho/\Sigma))$. In fact, something more general is true: if $\phi \colon (M,H') \to \JM$ and $\phi \colon (M,H'') \to \JM$ are isomorphic CIRs, then $\eta_{\mathcal{F}}(\theta,\theta') - \eta_{\mathcal{F}}(\theta,\theta'') \in \hat{\mathrm{D}}(\Gamma(\ker \rho/\Sigma))$, where $\theta$, $\theta'$, $\theta''$ are the generalised contact forms associated to $H$, $H'$ and $H''$ respectively. This follows from the fact that
 \begin{gather*}
 \eta_{\mathcal{F}}(\theta',\theta'') = \eta_{\mathcal{F}}(\theta,\theta'') - \eta_{\mathcal{F}}(\theta,\theta'),
 \end{gather*}
which is a direct consequence of the def\/inition of the forms $\eta_{\mathcal{F}}(\theta',\theta'')$, $\eta_{\mathcal{F}}(\theta,\theta'')$, and $ \eta_{\mathcal{F}}(\theta, \theta')$~-- to def\/ine $\eta_{\mathcal{F}}(\theta',\theta'')$, $\phi \colon (M,H') \to \JM$ is taken to be the `reference' CIR. The above argument proves that $\eta_{\mathcal{F}}(\theta',\theta'') \in \hat{\mathrm{D}}(\Gamma(\ker \rho/\Sigma))$, which yields the claimed result. Moreover, if $\eta_{\mathcal{F}} = \hat{\mathrm{D}}\tau$ for some $\tau \in \Gamma(\ker \rho/\Sigma)$, reverse the above reasoning to obtain that the CIRs $\phi \colon \CM \to \JM$ and $\phi \colon (M,H') \to \JM$ are isomorphic, where $H' = \ker (\theta + \phi^*\eta)$, where $\eta$ is a 1-form extending $\eta_{\mathcal{F}}$. In conclusion, the above discussion implies that there exists a well-def\/ined bijection
 \begin{gather} \label{eq:18}
\frac{\left\{ \text{CIRs with period lattice } \Sigma \text{ and Chern
 class } c\right\}}{\sim} \to \frac{\mathrm{Z}^1(\mathcal{F};L)}{\hat{\mathrm{D}}(\Gamma(\ker \rho/\Sigma))}, \\
\left[ \phi' \colon (M',H') \to \JM\right] \mapsto \eta(\theta,\theta') + \hat{\mathrm{D}}(\Gamma(\ker \rho/\Sigma)),\nonumber
 \end{gather}
where $\sim$ denotes the equivalence relation def\/ined by the notion of isomorphism of Def\/inition~\ref{defn:iso_cir}. To complete the proof, it suf\/f\/ices to prove that $\hat{\mathrm{D}}(\Gamma(\ker \rho/\Sigma))$ is precisely the vector space of exact, globally def\/ined foliated 1-forms with values in~$L$. By Note~\ref{obs:sod}, if $\alpha + \Sigma \in \Gamma (\ker \rho/\Sigma )$, then $\hat{\mathrm{D}}(\alpha + \Sigma) = \de_{\mathcal{F}} (\operatorname{pr}(\alpha) )$. Therefore, $\hat{\mathrm{D}}(\Gamma(\ker \rho/\Sigma))$ consists of exact, foliated 1-forms with values in~$L$. Suppose that $\de_{\mathcal{F}}u$ is an exact foliated 1-form with values in~$L$, where $u \in \Gamma(L)$. Since $\operatorname{pr}\colon \ker \rho \to L$ is onto (cf.\ Note~\ref{obs:prop_jac_even}), it follows that there exists a section $\alpha \in \Gamma(\ker \rho)$ with $u = \operatorname{pr}(\alpha)$. Then $\hat{\mathrm{D}}(\alpha + \Sigma) = \de_{\mathcal{F}}u$ by construction, which shows that all exact, foliated 1-forms with values in $L$ are contained in $\hat{\mathrm{D}} (\Gamma (\ker \rho/\Sigma))$. Hence,
 \begin{gather*}
 \frac{\mathrm{Z}^1(\mathcal{F};L)}{\hat{\mathrm{D}}(\Gamma(\ker \rho/\Sigma))} = \HH^1 (\mathcal{F};L ),
 \end{gather*}
therefore completing the proof.
\end{proof}

\section[The case of Poisson manifolds: comparing symplectic and contact isotropic\\ realisations]{The case of Poisson manifolds: comparing symplectic\\ and contact isotropic realisations}\label{sec:sympl-vers-cont}

\subsection{Symplectic isotropic realisations of Poisson manifolds: a reminder}\label{sec:sympl-isotr-real}
For completeness, the theory of symplectic isotropic realisations of Poisson manifolds is recalled here (cf.~\cite{dd} for further details and proofs of all results stated below).

\begin{Definition}\label{defn:sir}
A {\it symplectic realisation} of a Poisson manifold $\PM$ is a~symplectic manifold $\SM$ together with a surjective submersion $\Phi \colon \SM \to \PM$ which is a Poisson morphism. If, in addition, the f\/ibres of $\Phi$ are isotropic submanifolds of $\SM$, the realisation is said to be {\it isotropic}.
\end{Definition}

Symplectic realisations play an important role in the study of Poisson manifolds (cf.~\cite{cdw}); while symplectic isotropic realisations ({\it SIR} for short) appear naturally when considering Poisson manifolds of compact types (cf.~\cite{PMCT1,PMCT2}). Henceforth, all SIRs considered have compact f\/ibres unless otherwise stated.

It is well-known that a Poisson manifold $\PM$ admits a SIR only if it is regular (cf.~\cite{dd} and compare with Lemma~\ref{lemma:regular}).

\begin{Note}\label{obs:poisson_coho}
There are two important cohomology classes that can be attached to a~regular Poisson manifold $\PM$ with symplectic foliation~$\mathcal{F}$. On the one hand, the Poisson bivector $\Lambda$ induces a foliated cohomology class $[\omega_{\mathcal{F}}] \in \HH^2(\mathcal{F})$ (cf.\ Example~\ref{exm:Poisson_case}). On the other, if $\omega \in \Omega^2(P)$ is a~2-form $\omega \in \Omega^2(P)$ extending $\omega_{\mathcal{F}}$, $\de \omega$ is a closed 3-form whose restriction to the leaves of~$\mathcal{F}$ is zero; therefore, it def\/ines a cohomology class $[\de \omega ] \in \HH^3_{\mathrm{rel}}(P;\mathcal{F})$, where $\HH^*_{\mathrm{rel}}(P;\mathcal{F})$ is the cohomology of the subcomplex of forms that vanish along $\mathcal{F}$ (this is sometimes referred to as the cohomology of $P$ {\em relative} to~$\mathcal{F}$). This cohomology class is independent of the choice of $\omega$; the class $\upsilon:=[\de \omega]$ is known as the {\it characteristic class} of~$\PM$.
\end{Note}

Fix a SIR $\Phi \colon \SM \to \PM$, let $\mathcal{F}$ and $\nu^*$ denote the regular symplectic foliation on $\PM$ and its conormal bundle respectively. The given SIR induces a full-rank lattice $\Xi \subset \nu^*$ with the following properties:
\begin{itemize} \itemsep=0pt
\item $\Phi \colon \SM \to \PM$ is a principal $\nu^*/\Xi$-bundle, which is classif\/ied by an element (called {\it Chern class}) in $\HH^1 (P;\mathcal{C}^{\infty} (\nu^*/\Xi ) ) \cong \HH^2(P;\underline{\Xi})$, where $\mathcal{C}^{\infty}(\nu^*/\Xi)$ and $\underline{\Xi}$ denote the sheaves of sections of $\nu^*/ \Xi \to P$ and of $\Xi \to P$ respectively (cf.\ \cite[Corollary~2]{dd});
\item any locally def\/ined section of $\Xi \to P$ is a closed 1-form (cf.\ \cite[Corollary~1]{dd}).
\end{itemize}

The full-rank lattice $\Xi \subset \nu^*$ encodes an important geometric structure transverse to the foliation $\mathcal{F}$.

\begin{Definition}\label{defn:tias}
Let $(N,\mathcal{F})$ be a foliated manifold, where the codimension of the foliation is $l \geq 0$. A {\it transversal $\ZZ$-affine structure} on $(N,\mathcal{F})$ is an atlas $\mathcal{A} :=\{(U_i, \chi_i)\}$ of submersions $\chi_i \colon U_i \to \RR^l$ locally def\/ining $\mathcal{F}$, such that, for all $i,j$ with $U_{ij} := U_i \cap U_j \neq \varnothing$, there exists a~smooth map $h_{ij}\colon \chi_i(U_{ij}) \to \chi_j(U_{ij})$ which is (the restriction of) an $\ZZ$-af\/f\/ine transformation in $\operatorname{Af\/f}_{\ZZ}(\RR^l) := \mathrm{GL}(l;\ZZ) \ltimes \RR^l$ such that $\chi_j|_{U_{ij}} = h_{ij} \circ \chi_i|_{U_{ij}}$.
\end{Definition}

\begin{Note}\label{obs:trans_aff_structure}
Allowing the maps $h_{ij}$ to be (restrictions of) elements in $\operatorname{Af\/f}(\RR^l) := \mathrm{GL}(l;\RR) \ltimes \RR^l$, obtain the notion of {\it transversal affine structure} on a foliated manifold, which stands to its integral counterpart as transversal projective structures stand to their integral analogues (cf.\ Def\/initions~\ref{defn:tipa} and~\ref{obs:trans_proj_structure}).
\end{Note}

The following proposition, stated without proof, establishes the analogue of Proposition \ref{thm:tilp_tilpa} for transversal $\ZZ$-af\/f\/ine structures (cf.\ \cite[Proposition~3.2.4]{PMCT2}).

\begin{Proposition}\label{prop:correspondence}
If $(N,\mathcal{F})$ is a manifold together with a foliation of codimension $l$, there is a~{\rm 1-1} correspondence between
 \begin{itemize}\itemsep=0pt
 \item transversal $\ZZ$-affine structures on $(N,\mathcal{F})$;
 \item full-rank lattices $\Xi \subset \nu^*$ whose local sections are closed.
 \end{itemize}
In this correspondence, $\Xi$ is locally given by
 \begin{gather*}
 \Xi = \ZZ \langle \de x_1,\ldots , \de x_l \rangle,
 \end{gather*}
where $x_1,\ldots,x_l$ are local transversal $\ZZ$-affine coordinates on $(N,\mathcal{F})$.
\end{Proposition}

Henceforth, any full-rank lattice corresponding to a transversal $\ZZ$-af\/f\/ine structure as in Proposition~\ref{prop:correspondence} is referred to as a~{\it transversal $\ZZ$-affine lattice}. The realisation problem for SIRs has been solved in~\cite{dd}. Suppose that $\PM$ is a regular Poisson manifold endowed with a transversal $\ZZ$-af\/f\/ine lattice $\Xi \subset \nu^*$. Since the sections of $\Xi$ are closed, the standard exterior derivative induces a~homomorphism of sheaves
\begin{alignat*}{3}
 & \hat{\de} \colon \ && \mathcal{C}^{\infty}\left(\nu^* /\Xi \right) \to \mathcal{Z}^2(\nu^*), & \\
 &&& \alpha + \Xi \mapsto \de \alpha,&
\end{alignat*}
where $\mathcal{Z}^2(\nu^*)$ denotes the sheaf of closed 2-forms which vanish when restricted to $\mathcal{F}$. For any $l \geq 0$, this induces homomorphisms in cohomology
\begin{gather} \label{eq:16}
 \mathfrak{d}\colon \ \HH^l\big(P;\mathcal{C}^{\infty} (\nu^* /\Xi )\big) \to \HH^l\big(P;\mathcal{Z}^2(\nu^*)\big).
\end{gather}
For $l \geq 1$, there exist isomorphisms
\begin{gather*}
 \HH^l\big(P;\mathcal{C}^{\infty} (\nu^* /\Xi )\big) \cong \HH^{l+1} (P;\underline{\Xi} )
\end{gather*}
and
\begin{gather*}
 \HH^l\big(P;\mathcal{Z}^2(\nu^*)\big) \cong \HH_{\mathrm{rel}}^{l+2}(P;\mathcal{F}).
\end{gather*}
The former arise from considering the short exact sequence of sheaves
\begin{gather*} 
 0 \to \underline{\Xi} \to \mathcal{C}^{\infty}(\nu^*) \to \mathcal{C}^{\infty}(\nu^*/\Xi) \to 0,
\end{gather*}
and by observing that $\mathcal{C}^{\infty}(\nu^*)$ is f\/ine; the latter are explained in \cite[Section~4, Part~(a)]{dd}. The homomorphisms of equation \eqref{eq:16} induce, for all $l \geq 1$, homomorphisms $\HH^{l+1}(P;\Xi) \to \HH^{l+2}_{\mathrm{rel}}(P;\mathcal{F})$, which are also denoted by $\mathfrak{d}$ by abuse of notation. With the above homomorphisms at hand, the following results can be stated (cf.\ \cite[Theorems~4.2 and~4.3]{dd} for a~proof, and compare with Theorem~\ref{thm:main}).

\begin{Theorem}\label{thm:real_sirs}\quad
 \begin{enumerate}[label={\rm \arabic*.}, ref={\rm (\arabic*)}]\itemsep=0pt
\item \label{item:16} A regular Poisson manifold $\PM$ admits a~symplectic isotropic realisation with period lattice~$\Xi$ and Chern class $c \in \HH^2(P;\underline{\Xi})$ if and only if
 \begin{gather} \label{eq:40}
 \mathfrak{d}c = \upsilon,
 \end{gather}
 where $ \upsilon \in \HH_{\mathrm{rel}}^3(P;\mathcal{F})$ is the characteristic class of $\PM$.
 \item \label{item:17} Suppose that $c \in \HH^2(P;\Xi)$ satisfies equation \eqref{eq:40}, then the set of isomorphism classes of symplectic isotropic realisations of $\PM$ with period lattice $\Xi$ and Chern class $c$ is in bijection with
 \begin{gather*}
 \frac{\HH^0\big(P; \mathcal{Z}^2(\nu^*)\big)}{\mathfrak{d} \big(\HH^0 (P;\mathcal{C}^{\infty}(\nu^*/\Xi)
 )\big)} = \frac{\mathrm{Z}^2(\nu^*)}{\hat{\de}(\Gamma(\nu^*/\Xi))},
 \end{gather*}
 where $\mathrm{Z}^2(\nu^*)$ denotes the vector space of globally defined closed $2$-forms which vanish along the symplectic foliation $\mathcal{F}$.
 \end{enumerate}
\end{Theorem}

\begin{Note}\label{obs:alt_map}
The map $\mathfrak{d} \colon \HH^2(P;\underline{\Xi}) \to \HH^3_{\mathrm{rel}}(P;\mathcal{F})$ can be constructed alternatively as follows. There is a~commutative diagram of short exact sequence of sheaves
 \begin{gather} \label{eq:44}
 \begin{split} &
 \xymatrix{ 1 \ar[r] & \underline{\Xi} \ar@{^{(}->}[d] \ar@{^{(}->}[r] &
 \mathcal{C}^{\infty}(\nu^*) \ar[r] \ar@{=}[d] &
 \mathcal{C}^{\infty}(\nu^*/\Xi)
 \ar[d]^-{\hat{\de}} \ar[r] & 1\\
 1 \ar[r] & \mathcal{Z}^{1}(\nu^*)
 \ar@{^{(}->}[r] & \mathcal{C}^{\infty} (\nu^*)
 \ar[r]_-{\de} & \mathcal{Z}^2(\nu^*) \ar[r] & 1,}
 \end{split}
 \end{gather}
 which induces a commutative diagram
 \begin{gather} \label{eq:45}
 \begin{split}& \xymatrix{\HH^1(P;\mathcal{C}^{\infty}(\nu^*
 /\Sigma))\ar[d]_-{\mathfrak{d}} \ar[r]^-{\cong} & \HH^2(P;\underline{\Xi})
 \ar[d] \\
 \HH^1(P;\mathcal{Z}^{2}(\nu^*)) \ar[r]^-{\cong} & \HH^2(P;\mathcal{Z}^1(\nu^*)),}
 \end{split}
 \end{gather}
where the vertical maps are induced by the outer vertical maps in equation \eqref{eq:44} and the horizontal maps are the connecting morphisms in the long exact sequences in cohomology induced by the short exact sequences of equation~\eqref{eq:44}, as $\mathcal{C}^{\infty}(\nu^*)$ is a f\/ine sheaf. Under the identif\/ications of equation \eqref{eq:45}, the map $\mathfrak{d}$ can be identif\/ied with the homomorphism $\HH^2(P;\underline{\Xi}) \to \HH^2(P;\mathcal{Z}^1(\nu^*)) \cong \HH^3_{\mathrm{rel}}(P;\mathcal{F})$ induced by the inclusion $\underline{\Xi} \hookrightarrow \mathcal{Z}^1(\nu^*)$.
\end{Note}

\subsection{Contact isotropic realisations of Poisson manifolds}\label{sec:cont-isotr-real-1}
It is worth rephrasing some of the results of Section~\ref{sec:class-cir-over} in the case in which the underlying Jacobi manifold is, in fact, Poisson, as many of the above objects simplify (cf.\ Section~\ref{sec:sympl-vs-cont} for more applications). Using the notation of Example~\ref{obs:kir}, given a Poisson manifold $\PM$, denote the induced Jacobi structure on the trivial bundle $\RR_P$ also by $\PM$. Suppose that it is regular and denote its symplectic foliation and corresponding conormal bundle by $\mathcal{F}$, $\nu^*$ respectively. Recall that if $\rho \colon J^1P \cong T^*P \oplus \RR \to TP$ is the anchor, then $\ker \rho = \nu^* \oplus \RR$ (cf.\ Example~\ref{exm:Poisson_case}). Fix a transversal $\ZZ$-projective lattice $\Sigma \subset \ker \rho$; by def\/inition, any local section of $\Sigma$ is of the form $j^1 f$, for some $f \in C^{\infty}(P)$. However, since $\Sigma \subset \nu^* \oplus \RR$,
\begin{gather} \label{eq:21}
 j^1 f \in \Gamma(\Sigma) \quad \Rightarrow \quad \de f \in \Gamma(\nu^*) \quad \Leftrightarrow \quad f \text{ is a Casimir}.
\end{gather}

\begin{Remark}\label{rk:general_ker_Casimir}
Let $\JM$ be a regular Jacobi manifold all whose leaves are even di\-men\-sio\-nal. Say that a (locally def\/ined) section $u \in \Gamma(L)$ is a {\it Casimir}, if, for all $v \in \Gamma(L)$, $\{u,v\}=0$. Then, using the fact that the structure is regular and all its leaves are even di\-men\-sio\-nal, it can be shown that $u \in \Gamma(L)$ is a Casimir if and only if $j^1u \in \Gamma(\ker \rho)$.
\end{Remark}

The equality $\ker \rho = \nu^* \oplus \RR$ can be exploited to give another way to construct the homomorphism of equation \eqref{eq:18} which underpins the cohomological criterion of Theorem~\ref{thm:main}. Denote the sheaf of {\it basic} smooth functions on~$(P,\mathcal{F})$ by $\mathcal{C}^{\infty}_{\mathrm{basic}}(P;\mathcal{F})$, i.e., it consists of smooth functions which are locally constant on the leaves of~$\mathcal{F}$. Equivalently, this can be def\/ined as the sheaf of Casimirs of~$\PM$, seeing as, in this case, a~functions is basic for~$\mathcal{F}$ if and only if it is a Casimir. There is a short exact sequence of sheaves
\begin{gather*}
 1 \to\mathcal{C}^{\infty}_{\mathrm{basic}}(P;\mathcal{F}) \hookrightarrow\mathcal{C}^{\infty} (P) \xrightarrow{\de_{\mathcal{F}}}\mathcal{Z}^1(\mathcal{F}) \to 1,
\end{gather*}
where $\mathcal{C}^{\infty} (P)$ and $\mathcal{Z}^1(\mathcal{F})$ are the sheaves of smooth functions on $P$ and of closed foliated 1-forms respectively.

\begin{Lemma}\label{lemma:poisson}
 The following is a commutative diagram of short exact sequences of sheaves
 \begin{gather} \label{eq:19}
 \begin{split} &
 \xymatrix{ 1 \ar[r] & \underline{\Sigma} \ar[d]^-{\pr} \ar@{^{(}->}[r] &
 \mathcal{C}^{\infty}(\ker \rho) \ar[r] \ar[d]^-{\pr} &
 \mathcal{C}^{\infty}(\ker \rho/\Sigma)
 \ar[d]^-{\hat{\mathrm{D}}} \ar[r] & 1\\
 1 \ar[r] & \mathcal{C}^{\infty}_{\mathrm{basic}}(P;\mathcal{F})
 \ar@{^{(}->}[r] & \mathcal{C}^{\infty} (P)
 \ar[r]_-{\de_{\mathcal{F}}} & \mathcal{Z}^1(\mathcal{F}) \ar[r] & 1,}
 \end{split}
 \end{gather}
where $pr$ denotes the homomorphism of sheaves induced by the projection $\pr \colon \cotan P \oplus \RR \to \RR_{P}$, $\hat{\mathrm{D}}$ is the homomorphism of sheaves defined by equation~\eqref{eq:20}, and~$\underline{\Sigma}$ denotes the sheaf of smooth sections of~$\Sigma \to P$.
\end{Lemma}
\begin{proof}
First, observe that equation \eqref{eq:21} implies that the image $\pr(\Sigma)$ of the sheaf homomorphism $\pr \colon \mathcal{C}^{\infty}(\ker \rho) \to \mathcal{C}^{\infty}(P)$ lies in the sheaf $\mathcal{C}^{\infty}_{\mathrm{basic}}(P;\mathcal{F})$. The only non-trivial fact that needs checking is commutativity on the right hand side of the diagram~\eqref{eq:19}. Let $[\alpha] \in \Gamma(\ker \rho/\Sigma)$ and set $\alpha = (\eta,f) \in \Gamma(\ker \rho)$ be a lift of $[\alpha]$. Then
 \begin{gather*}
\hat{\mathrm{D}}([\alpha]) = \mathrm{D}(\eta,f)|_{\mathcal{F}} = (\de f - \eta)|_{\mathcal{F}} = \de_{\mathcal{F}} f = \de_{\mathcal{F}} \circ \pr(\alpha),
 \end{gather*}
where the third equality follows from the fact that $\ker \rho = \nu^* \oplus \RR$.
\end{proof}

Since both $\mathcal{C}^{\infty}(P)$ and $\mathcal{C}^{\infty}(\ker \rho)$ are f\/ine sheaves, there is a commutative diagram
\begin{gather} \label{eq:22}
\begin{split} &
 \xymatrix{\HH^1(P;\mathcal{C}^{\infty}(\ker\rho /\Sigma))\ar[d]_-{\mathcal{D}} \ar[r]^-{\cong} & \HH^2(P;\underline{\Sigma}) \ar[d] \\
 \HH^1(P;\mathcal{Z}^{1}(\mathcal{F})) \ar[r]^-{\cong} & \HH^2(\mathcal{F}),}
\end{split}
\end{gather}
where the vertical maps are induced by the outer vertical maps of diagram~\eqref{eq:19} (cf.\ equation~\eqref{eq:18} also for the def\/inition of~$\mathcal{D}$), and the horizontal isomorphisms are the connecting morphisms induced by the short exact sequences in equation~\eqref{eq:19}. Therefore, the following corollary holds.

\begin{Corollary}\label{cor:po}
 The map in cohomology
 \begin{gather*}
 \HH^2(P;\underline{\Sigma}) \to \HH^2\big(P;\mathcal{C}^{\infty}_{\mathrm{basic}}(P;\mathcal{F})\big) \cong \HH^2(\mathcal{F})
 \end{gather*}
induced by the sheaf homomorphism $\pr \colon \underline{\Sigma} \to \mathcal{C}^{\infty}_{\mathrm{basic}}(P;\mathcal{F})$ equals $\mathcal{D}$ via the identifications of equation \eqref{eq:22}.
\end{Corollary}

Corollary \ref{cor:po} indicates how to calculate the map~$\mathcal{D}$, which is central to understanding whether a regular Poisson manifold $\PM$ with given transversal $\ZZ$-projective lattice $\Sigma \subset \nu^* \oplus \RR$ admits a CIR (cf.\ Part~\ref{item:19} of Theorem~\ref{thm:main}).

\begin{Example} \label{exm:trivial_tips}
Suppose that $\Sigma = \ZZ \langle j^1 f_1 ,\ldots, j^1 f_k \rangle$ for some Casimirs $f_i$, so that it induces the trivial $\ZZ^k$-system of coef\/f\/icients. Then (modulo torsion), an element in $\HH^2(P;\Sigma)$ can be written as
 \begin{gather*}
 \sum\limits_{i=1}^k [\omega_i] \otimes j^1 f_i,
 \end{gather*}
where, for each $i$, $\omega_i \in \Omega^2(P)$ is a 2-form with integral cohomology class. Since the map $\mathcal{D} \colon \HH^2(P;\underline{\Sigma}) \to \HH^2(\mathcal{F})$ is induced by the projection $\pr \colon \Sigma \to \mathcal{C}^{\infty}_{\mathrm{basic}}(P;\mathcal{F})$ and $\Sigma$ is a trivial bundle,
\begin{gather*}
 \mathcal{D} \left(\sum\limits_{i=1}^k [\omega_i] \otimes j^1 f_i\right) = \sum\limits_{i=1}^k [f_i \omega_i],
 \end{gather*}
where, for each $i$, $f_i\omega_i$ is a closed foliated 2-form since $f_i$ is basic.
\end{Example}

\subsection{Strong transversal $\ZZ$-af\/f\/ine structures are $\ZZ$-projective}\label{sec:strong-transv-zz}
Before tackling the main problem of this section, it is necessary to show that, under some integrality condition (cf.\ Def\/inition~\ref{defn:strong_Z_aff}), transversal $\ZZ$-af\/f\/ine structures induce transversal $\ZZ$-projective structures. Let $(N,\mathcal{F})$ be a foliated manifold and recall the def\/initions of transversal af\/f\/ine and projective structures on~$(N,\mathcal{F})$ (cf.\ Notes~\ref{obs:trans_proj_structure} and~\ref{obs:trans_aff_structure}). The whole idea of this subsection hinges on the following well-known result, whose proof is included for completeness.

\begin{Lemma}\label{lemma:affine_as_proj}
Any transversal affine structure on $(N,\mathcal{F})$ induces a~transversal projective structure on $(N,\mathcal{F})$.
\end{Lemma}
\begin{proof}
Let $l \geq 0$ denote the codimension of $\mathcal{F}$. If $l =0$, there is nothing to prove, so suppose that $l \geq 1$. Consider the smooth map
 \begin{alignat*}{3}
 & \mathcal{I} \colon \ && \RR^l \to \RR\mathrm{P}^l, & \\
 &&& \left(x_1,\ldots,x_l\right) \mapsto [x_1:\ldots:x_l :1];&
 \end{alignat*}
this map is a dif\/feomorphism onto its image. Moreover, if
 \begin{gather} \label{eq:41}
 \begin{split}
 \mathsf{I} \colon \ \operatorname{Af\/f}(\RR^l) &\hookrightarrow \mathrm{GL}(l+1;\RR), \\
 (A,\mathbf{b}) &\mapsto
 \begin{pmatrix}
 A & \mathbf{b} \\
 0 & 1
 \end{pmatrix}
 \end{split}
 \end{gather}
denotes the standard inclusion as groups, then for all $(A,\mathbf{b}) \in \operatorname{Af\/f}(\RR^l)$ and for
 all $\mathbf{x} \in \RR^l$,
 \begin{gather} \label{eq:43}
 \mathcal{I} (A \mathbf{x} + \mathbf{b} ) = [\mathsf{I}(A,\mathbf{b})](\mathcal{I}(\mathbf{x})),
 \end{gather}
where $[\cdot] \in \mathrm{PGL}(l+1;\RR)$ denotes the equivalence class of $\cdot \in \mathrm{GL}(l+1;\RR)$. Suppose that $\mathcal{A}=\{(U_i, \chi_i)\}$ is a transversal af\/f\/ine structure on $(N,\mathcal{F})$, and suppose that, for $i,j$ with $U_{ij} \neq \varnothing$, $h_{ij} \in \operatorname{Af\/f}(\RR^l)$ are the transversal af\/f\/ine changes of coordinates as in Note~\ref{obs:trans_aff_structure}. Then equation~\eqref{eq:43} implies that $\bar{\mathcal{A}}:=\{(U_i,\bar{\chi}_i)\}$, where $\bar{\chi}_i := \mathcal{I} \circ \chi_i$ def\/ines a transversal projective structure on~$(N,\mathcal{F})$, the transversal projective changes of coordinates being $\bar{h}_{ij}=[\mathsf{I}(h_{ij})]$.
\end{proof}

Lemma \ref{lemma:affine_as_proj} implies that any transversal $\ZZ$-af\/f\/ine structure induces a transversal projective structure; however, the latter need not be integral! Fundamentally, the reason is that, for any $l \geq 1$, $\mathsf{I}(\operatorname{Af\/f}_{\ZZ}(\RR^l))\not\subset \mathrm{GL}(l+1;\ZZ)$, the issue being that the translational components of $\ZZ$-af\/f\/ine transformations are not necessarily integral. This leads to the following def\/inition\footnote{The terminology in Def\/inition~\ref{defn:strong_Z_aff} is not standard. For instance, \cite{gross_siebert} refers to the notion of strong (transversal) $\ZZ$-af\/f\/ine structure simply as (transversal) $\ZZ$-structure.}, which can be thought of as an `integrality condition' for a transversal $\ZZ$-af\/f\/ine structure.

\begin{Definition}\label{defn:strong_Z_aff}
A transversal $\ZZ$-af\/f\/ine structure $\mathcal{A}=\{(U_i,\chi_i)\}$ on $(N,\mathcal{F})$ is said to be {\it strong} if, for all~$i$,~$j$ with $U_{ij} \neq \varnothing$, $h_{ij}$ is (a restriction of) an element in $\operatorname{Af\/f}(\ZZ^l):= \mathrm{GL}(l;\ZZ) \ltimes \ZZ^l$, where $l$ is the codimension of $\mathcal{F}$.
\end{Definition}

\begin{Corollary}\label{cor:strong_Z_proj}
Any strong transversal $\ZZ$-affine structure on $(N,\mathcal{F})$ induces a transversal $\ZZ$-projective structure on $(N,\mathcal{F})$.
\end{Corollary}
\begin{proof}
Let $\mathcal{A}=\{(U_i,\chi_i)\}$ be a strong transversal $\ZZ$-structure on $(N,\mathcal{F})$ and, as in the proof of Lemma~\ref{lemma:affine_as_proj}, let $h_{ij} \in \operatorname{Af\/f}(\ZZ^l)$ be the transversal $\ZZ$-af\/f\/ine changes of coordinates. Setting $A_{ij}:= \mathsf{I}(h_{ij}) \in \mathrm{GL}(l+1;\ZZ)$, observe that $\{A_{ij}\}$ satisf\/ies the cocycle condition, for $h_{ij}$ does. This fact can be used, together with Lemma~\ref{lemma:affine_as_proj}, to prove that the atlas $\bar{\mathcal{A}}$ constructed in the proof of Lemma~\ref{lemma:affine_as_proj} def\/ines a~transversal $\ZZ$-projective structure (cf.\ Def\/inition~\ref{defn:tipa}).
\end{proof}

Henceforth, f\/ix a strong transversal $\ZZ$-af\/f\/ine structure $\mathcal{A}=\{(U_i,\chi_i)\}$ on $(N,\mathcal{F})$ which, by Corollary~\ref{cor:strong_Z_proj} induces a transversal $\ZZ$-projective structure $\bar{\mathcal{A}} = \{(U_i,\bar{\chi}_i)\}$, which is also f\/ixed. By Proposition~\ref{prop:correspondence}, $\mathcal{A}$ corresponds to a full-rank lattice $\Xi \subset \nu^*$ whose sections are closed, where $\nu^*$ denotes the conormal bundle to~$\mathcal{F}$. On the other hand, Proposition~\ref{thm:tilp_tilpa} ensures the existence of a line bundle $L \to N$ and a transversal
$\ZZ$-projective lattice $\Sigma \subset J^1L$ satisfying conditions \ref{item:8}--\ref{item:9} in Def\/inition~\ref{defn:tilp} which corresponds to $\bar{\mathcal{A}}$. Seeing as $\bar{\mathcal{A}}$ is induced by $\mathcal{A}$, it is natural to ask what relation there is between $\Xi$ and $\Sigma$. To this end, it is useful to recall that, for each~$i$,
\begin{itemize}\itemsep=0pt
\item $\Xi|_{U_i}= \chi_i^* \Xi^l$, where $\Xi^l:=\ZZ \langle \de x_1,\ldots , \de x_l \rangle \subset T^* \RR^l$ is the standard $\ZZ$-af\/f\/ine structure on $\RR^l$ (cf.\ Proposition~\ref{prop:correspondence});
\item $L|_{U_i} = \bar{\chi}_i^*(O(1))$ and $\Sigma|_{U_i} = \bar{\chi}_i^* \Sigma^l$, where $O(1) \to \RR\mathrm{P}^l$ is the dual of the tautological line bundle and $\Sigma^l \subset J^1(O(1))$ is the $\ZZ$-projective lattice of Example~\ref{ex:projective_plane} (cf.\ Proposition~\ref{thm:tilp_tilpa}). In fact, since $\bar{\chi}_i = \mathcal{I} \circ \chi_i$ (cf.\ the proof of Lemma~\ref{lemma:affine_as_proj}), $L|_{U_i} = \chi^*_i(\mathcal{I}^*(O(1)))$ and $\Sigma|_{U_i} = \chi^*_i(\mathcal{I}^*(\Sigma^l))$.
\end{itemize}

\begin{Claim}\label{claim:L_trivial}
 The line bundle $L \to N$ is trivialisable.
\end{Claim}
\begin{proof}
First, observe that the line bundle $\mathcal{I}^*(O(1)) \to \RR^l$ is trivialisable, for the restriction of $O(1)$ to $\mathcal{I}(\RR^l) = \{[x_1,\ldots,x_l,x_{l+1}] \,|\, x_{l+1} \neq 0\}$ admits a nowhere vanishing section $\zeta$, namely the functional
 \begin{gather*}
 \begin{split}
 & \RR\langle x_1,\ldots,x_{l+1}\rangle \to \RR, \\
& \mathbf{y} = (y_1,\ldots,y_{l+1}) \mapsto y_{l+1}.\end{split}
 \end{gather*}
The restriction of the natural action of $\mathrm{GL}(l+1;\RR)$ on~$\Gamma(O(1))$ to $\mathsf{I}(\operatorname{Af\/f}(\ZZ^l))$ f\/ixes $\zeta$: this is essentially because of the def\/inition of the homomorphism $\mathsf{I}$ (cf.\ equation~\eqref{eq:41}). For each~$i$, the section $\chi^*_i \circ \mathcal{I}^*\zeta$ of $L|_{U_i}$ is nowhere vanishing. If $U_{ij} \neq \varnothing$, then
 \begin{gather*}
 \big(\chi^*_j \circ \mathcal{I}^*\zeta\big)|_{U_{ij}} = \big(\chi^*_i \circ h^*_{ij} \circ
 \mathcal{I}^*\zeta\big)|_{U_{ij}} =\big(\chi^*_i \circ \mathcal{I}^* \circ \mathsf{I}(h_{ij})^* \circ
 \zeta\big)|_{U_{ij}} =\big(\chi^*_i \circ \mathcal{I}^*\zeta\big)|_{U_{ij}},
 \end{gather*}
where the second equality follows from equation \eqref{eq:43} and the last one by the fact that $\zeta$ is f\/ixed by $\mathsf{I} (\operatorname{Af\/f}(\ZZ^l) )$. Therefore, $L \to N$ admits a globally def\/ined nowhere vanishing section, thus proving that it is trivialisable.
\end{proof}

Henceforth, f\/ix the trivialisations $\mathcal{I}^*(O(1)) \cong \RR_{\RR^l}$ and $L \cong \RR_{N}$ induced as in the proof of Claim~\ref{claim:L_trivial} unless otherwise stated. Since $L \to N$ is trivial, $\Sigma \subset J^1\RR_N$ is a full-rank lattice of $\nu^* \oplus \RR$ (cf.\ Example~\ref{exm:trivial_bundle}) all of whose sections are holonomic (cf.\ condition~\ref{item:9}).

\begin{Proposition}\label{prop:sigma_proj_xi}
If $\mathsf{P} \colon \nu^* \oplus \RR \to \nu^*$ denotes projection onto the first factor, then $\mathsf{P}(\Sigma) = \Xi$.
\end{Proposition}
\begin{proof}
First, it is shown that if $\mathsf{P}_{\RR^l} \colon T^* \RR^l \oplus \RR =J^1\RR_{\RR^l} \to T^*\RR^l$ is the projection onto the f\/irst factor, then
 $\mathsf{P}_{\RR^l}(\mathcal{I}^*(\Sigma^l)) = \Xi^l$. This follows basically by unravelling the def\/initions of~$\Sigma^l$ and~$\Xi^l$. By def\/inition of $\Sigma^l$ and of $\mathcal{I}$ (cf.\ Example~\ref{ex:projective_plane} and the proof of Lemma~\ref{lemma:affine_as_proj}),
 \begin{gather*}
 \mathcal{I}^*\Sigma^l= \ZZ\big \langle j^1x_1,\ldots,j^1x_l, j^1 1\big\rangle,
 \end{gather*}
where, for each $j$, $x_j$ denotes the functional which assigns to a~point in $\RR^l$ its $j$-th coordinate. Writing, for each $j$, $j^1x_l = (\de x_l,x_l) \in \Gamma(T^* \RR^l \oplus \RR)$, and $j^11=(0,1)$, it is clear that \smash{$\mathsf{P}_{\RR^l}(\mathcal{I}^*(\Sigma^l))= \Xi^l$}. The general case follows immediately by observing that, for all $i$, the restrictions $\mathsf{P}|_{U_i} \colon (\nu^* \oplus \RR)|_{U_i} \to \nu^*|_{U_i}$, $\Sigma|_{U_I}$ and $\Xi|_{U_i}$ are simply the pull-backs along $\chi_i$ of $\mathsf{P}_{\RR^l}$, $\mathcal{I}^* \Sigma^l$ and $\Xi^l$ respectively, where, as above, $\mathcal{A} = \{(U_i,\chi_i)\}$ is the f\/ixed transversal $\ZZ$-af\/f\/ine structure on $(N,\mathcal{F})$.
\end{proof}

\begin{Note}\label{obs:trivial_sub_lattice}
A closer look at the proof of Proposition \ref{prop:sigma_proj_xi} yields a slightly stronger result. Observe that $\mathcal{I}^* \Sigma^l$ contains a one-dimensional sub-lattice $\ZZ\langle j^11\rangle$ which is precisely the kernel of $\mathsf{P}_{\RR^l}\colon \mathcal{I}^* \Sigma^l \to \Xi^l$. The locally def\/ined one-dimensional sub-lattices $\chi^*_i\left(\ZZ\langle j^11\rangle\right) \subset \Sigma|_{U_i}$ patch together to yield a~globally def\/ined, trivial one-dimensional sub-lattice of~$\Sigma$ which is precisely the kernel of $\mathsf{P}\colon \Sigma \to \Xi$ (where, by abuse of notation, the restrictions of $\mathsf{P}_{\RR^l}$ and of $\mathsf{P}$ to their respective lattices are also denoted by $\mathsf{P}_{\RR^l}$ and $\mathsf{P}$ respectively). One invariant way to see this is to observe that one globally def\/ined section of $\Sigma$ is the f\/irst jet of the nowhere vanishing section of $L \to N$. Under the trivialisation $L \cong \RR_N$ induced by Claim~\ref{claim:L_trivial}, this section corresponds to the function $1$; in other words, $\ker \mathsf{P} = \ZZ \langle j^11 \rangle$. The projection~$\mathsf{P}$ and its restriction to $\Sigma$ induces the following commutative diagram of short exact sequences of bundle maps
 \begin{gather} \label{eq:42}\begin{split} &
 \xymatrix{ & 0 \ar[d] & 0 \ar[d] & 0 \ar[d] &\\
 0 \ar[r] & \ZZ \langle j^11 \rangle \cong \ZZ_N \ar[r] \ar[d] & \Sigma \ar[r]^-{\mathsf{P}}
 \ar[d] & \Xi \ar[d] \ar[r] & 0 \\
 0 \ar[r] & \RR\langle j^11 \rangle \cong \RR_N \ar[r] \ar[d] & \nu^* \oplus \RR
 \ar[r]^-{\mathsf{P}} \ar[d] & \nu^* \ar[d] \ar[r] & 0 \\
 0 \ar[r] & \left(\RR\langle j^11 \rangle / \ZZ \langle j^11
 \rangle \right) \cong S^1_{N} \ar[r] \ar[d] & \left(\nu^* \oplus
 \RR\right)/\Sigma \ar[r]^-{\bar{\mathsf{P}}} \ar[d] &
 \nu^*/\Xi \ar[r] \ar[d] & 0, \\
 & 0 & 0 & 0 &
 }\end{split}
 \end{gather}
where $S^1_N \to N$ is the trivial $S^1$-bundle over $N$ and $\bar{\mathsf{P}} \colon (\nu^* \oplus \RR)/\Sigma \to \nu^*/\Xi$ is the projection induced by~$\mathsf{P}$.
\end{Note}

By Proposition \ref{prop:sigma_proj_xi} and its proof, there is a homomorphism of sheaves
\begin{alignat*}{3} 
& \underline{\mathsf{P}} \colon \ && \underline{\Sigma} \to \underline{\Xi}, & \\
&&& j^1f \mapsto \de f, &
 \end{alignat*}
where $\underline{\Sigma}$ and $\underline{\Xi}$ denote the sheaves of sections of $\Sigma \to N$ and $\Xi \to N$ respectively.

\begin{Lemma}\label{lemma:ses_relation}
There is a commutative diagram of short exact sequences of sheaves
 \begin{gather} \label{eq:47}\begin{split}&
 \xymatrix{ 0 \ar[r] & \underline{\ZZ} \ar@{^{(}->}[r] \ar@{^{(}->}[d] &
 \underline{\Sigma} \ar[d]_-{pr} \ar[r]^-{\underline{\mathsf{P}}}
 & \underline{\Xi} \ar@{^{(}->}[d] \ar[r] & 0 \\
 0 \ar[r] & \underline{\RR} \ar@{^{(}->}[r] &
 \mathcal{C}^{\infty}_{\mathrm{basic}}(N;\mathcal{F})
 \ar[r]^-{\de} & \mathcal{Z}^1(\nu^*) \ar[r] & 0,}\end{split}
 \end{gather}
where $\underline{\ZZ}, \underline{\RR}$ are the sheaf of sections of $\ZZ_N$ and $\RR_N$ respectively, $\mathcal{C}^{\infty}_{\mathrm{basic}}(N;\mathcal{F})$ is the sheaf of basic functions, $\mathcal{Z}^1(\nu^*)$ is the sheaf of closed $1$-forms which vanish along $\mathcal{F}$, $\pr \colon \underline{\Sigma} \to \mathcal{C}^{\infty}_{\mathrm{basic}}(N;\mathcal{F})$ is the homomorphism induced by the projection $\nu^* \oplus \RR \to \RR_N$, and $\de \colon \mathcal{C}^{\infty}_{\mathrm{basic}}(N;\mathcal{F}) \to \mathcal{Z}^1(\nu^*) $ is induced by taking derivatives.
\end{Lemma}

\begin{proof} The top row is a short exact sequence as it is induced by the top short exact sequence of bundles of equation~\eqref{eq:42}, while it is well-known that the bottom row is a short exact sequence. The projection $\pr \colon \nu^* \oplus \RR \to \RR_N$ induces a homomorphism of sheaves of sections $\pr \colon \mathcal{C}^{\infty} (\nu^* \oplus \RR) \to \mathcal{C}^{\infty}(N)$; therefore, {\it a priori}, the codomain of its restriction to $\underline{\Sigma}$ is $\mathcal{C}^{\infty}(N)$. However, $\Sigma$ is a transversal $\ZZ$-projective structure on $(N,\mathcal{F})$ and, therefore, all its locally def\/ined sections are holonomic. Since $\Sigma \subset \nu^* \oplus \RR$, it follows that if $j^1f \in \Gamma_{\mathrm{loc}}(\Sigma)$, then $\pr(j^1f) = f$ is basic (cf.\ the implications in equation~\eqref{eq:21}). Therefore the middle vertical homomorphism of equation~\eqref{eq:47} is well-def\/ined. Commutativity of the diagram of equation~\eqref{eq:47} follows by def\/inition of~$\underline{\mathsf{P}}$.
\end{proof}

Consider the long exact sequences in cohomology associated to the top and bottom rows of equation~\eqref{eq:47}. The one induced by the bottom short exact sequence is well-known to be
\begin{gather*}
 \xymatrix@1{\cdots \ar[r] & \HH^2(N;\RR) \ar[r] & \HH^2(\mathcal{F}) \ar[r] &
 \HH^3_{\mathrm{rel}}(N;\mathcal{F}) \ar[r] & \HH^3(N;\RR) \ar[r] & \cdots,}
\end{gather*}
where the map $\HH^2(\mathcal{F}) \to \HH^3_{\mathrm{rel}}(N;\mathcal{F})$ sends the cohomology class of a~foliated 2-form to the cohomology class of the exterior dif\/ferential of any of its extensions as in Note~\ref{obs:poisson_coho} (cf.~\cite{CrainicFernandes:jdd,kacimi} for details of the above short exact sequence). Lemma~\ref{lemma:ses_relation} implies that there is the following commutative diagram
\begin{gather} \label{eq:48}\begin{split} &
 \xymatrix{\cdots \ar[r] & \HH^2(N;\ZZ) \ar[r] \ar[d] &
 \HH^2(N;\underline{\Sigma}) \ar[r]^-{\mathcal{P}}
 \ar[d]_-{\mathcal{D}} & \HH^2(N;\underline{\Xi})
 \ar[d]^-{\mathfrak{d}} \ar[r]^-{\delta} & \HH^3(N;\ZZ) \ar[d]
 \ar[r] & \cdots \\
 \cdots \ar[r] & \HH^2(N;\RR) \ar[r] & \HH^2(\mathcal{F}) \ar[r] &
 \HH^3_{\mathrm{rel}}(N;\mathcal{F}) \ar[r] & \HH^3(N;\RR) \ar[r] & \cdots,}\end{split}
\end{gather}
where
\begin{itemize}\itemsep=0pt
\item $\HH^*(N;\underline{\ZZ}), \HH^*(N;\underline{\RR})$ are identif\/ied with the singular cohomology groups $\HH^*(N;\ZZ)$ and $ \HH^*(N;\RR)$ respectively, and the maps $\HH^*(N;\ZZ) \to \HH^*(N;\RR) \cong \HH^*(N;\ZZ) \otimes_{\ZZ} \RR$ send $\HH^*(N;\ZZ)$ to $\HH^*(N;\ZZ) \otimes_{\ZZ} 1$;
\item $\mathcal{P} \colon \HH^2(N;\underline{\Sigma}) \to \HH^2(N;\underline{\Xi})$ is induced by $\underline{\mathsf{P}} \colon \underline{\Sigma} \to \underline{\Xi}$;
\item the maps $\mathcal{D} \colon \HH^2(N;\underline{\Sigma}) \to \HH^2(\mathcal{F})$ and $\mathfrak{d} \colon \HH^2(N;\underline{\Xi}) \to \HH^3_{\mathrm{rel}}(N;\mathcal{F})$ are the homomorphisms constructed in Sections~\ref{sec:class-cir-over} and~\ref{sec:sympl-isotr-real} respectively (cf.\ Corollary~\ref{cor:po} and Note~\ref{obs:alt_map}).
\end{itemize}

\begin{Note}\label{obs:geom_inter_P}
The homomorphism $\mathcal{P} \colon \HH^2(N;\underline{\Sigma}) \to \HH^2(N;\underline{\Xi})$ has a clear geometric interpretation. Suppose that $\phi \colon M \to N$ is a principal $(\nu^* \oplus \RR)/\Sigma$-bundle with Chern class $c \in \HH^2(N;\underline{\Sigma})$; then there is a free and proper action $(\nu^* \oplus \RR)/\Sigma \curvearrowright M$ along $\phi$. Restricting this action to the subbundle $(\RR \langle j^11\rangle/ \ZZ \langle j^11 \rangle ) \cong S^1_N \subset (\nu^* \oplus \RR)/\Sigma$ induces a free and proper $S^1$-action on $M$. Set $S:=M/S^1$. Since the $S^1$-action is tangent to the f\/ibres of $\phi$, there is a uniquely def\/ined smooth map $\Phi \colon S \to N$ along which the quotient $((\nu^* \oplus \RR)/\Sigma)/(\RR \langle j^11\rangle/ \ZZ \langle j^11 \rangle ) = \nu^*/\Xi$ acts freely and properly. Thus $\Phi \colon S \to N$ is a principal $\nu^*/\Xi$-bundle whose Chern class is precisely~$\mathcal{P}(c)$.
\end{Note}

\subsection{Symplectic vs contact isotropic realisations of Poisson manifolds}\label{sec:sympl-vs-cont}
Suppose that $\PM$ is a regular Poisson manifold and let $\Xi \subset \nu^*$ be a~{\it strong} transversal $\ZZ$-af\/f\/ine lattice on its symplectic foliation~$\mathcal{F}$. Using the notation of Section~\ref{sec:strong-transv-zz}, let $\Sigma \subset \nu^*\oplus \RR$ denote the transversal $\ZZ$-projective lattice induced by Corollary~\ref{cor:strong_Z_proj}. Then the following realisability criterion holds.

\begin{Theorem}\label{thm:cir_vs_sir}
The regular Poisson manifold $\PM$ admits a CIR with period lattice $\Sigma$ if and only if it admits a SIR $\Phi \colon \SM \to \PM$ with period lattice~$\Xi$ with the property that $\omega$ is integral.
\end{Theorem}

\begin{proof}
\underline{$(\Rightarrow)$}: Suppose that $\PM$ admits a CIR $\phi \colon \CM \to \PM$ with period lattice $\Sigma$. As in Note~\ref{obs:geom_inter_P}, consider the free and proper $S^1$-action on $\CM$ arising from restricting the $(\nu^*\oplus \RR)/\Sigma$-action on~$\CM$ to~$\RR \langle j^11\rangle /\ZZ \langle j^11\rangle \cong S^1_P$, and denote the resulting quotient map by $\Pi \colon M \to S := M/S^1$. There is an induced surjective submersion $\Phi \colon S \to P$ (which is a~principal $\nu^*/\Xi$-bundle by Note~\ref{obs:geom_inter_P}). The aim is to show that $S$ inherits an integral symplectic form from $\CM$. Since $\PM$ def\/ines a Jacobi structure on the trivial line bundle $\RR_P$, it follows that there is a contact form $\theta \in \Omega^1(M)$ with $H = \ker \theta$. Moreover, $\theta$ can be chosen so that the associated Reeb vector f\/ield $R_1$ is precisely the generator of the above principal $S^1$-action; this is because of the def\/inition of the action $(\nu^*\oplus \RR)/\Sigma \curvearrowright \CM$ and since the $S^1$-action arises by restricting to $\RR \langle j^11\rangle /\ZZ \langle j^11\rangle$. Thus $\theta$ is a connection 1-form for the principal bundle $\Pi\colon M \to S$ and, since it is a~contact form, its curvature $\omega \in \Omega^2(S)$ is, in fact, an integral symplectic form. To summarise the above discussion, the following diagram commutes
 \begin{gather} \label{eq:51}\begin{split} &
 \xymatrix{\CM \ar[rr]^-{\Pi} \ar[dr]_-{\phi}& & \SM
 \ar[dl]^-{\Phi} \\
 & \PM. &}\end{split}
 \end{gather}
It remains to show that $\Phi \colon \SM \to \PM$ is a SIR with period lattice $\Xi$. First, observe that~$\phi$ and~$\Pi$ being Jacobi maps, and $\Pi$ being a submersion, imply that $\Phi$ is Poisson map. To see that the f\/ibres of $\Phi$ are isotropic, observe that $D\Pi( \ker D\phi \cap H) = \ker D \Phi$. The condition of $\phi \colon \CM \to \PM$ being a contact {\it isotropic} realisation translates into $\ker D\phi \cap H \subset H$ being isotropic for $\de \theta$ (cf.\ Note~\ref{item:po}). The fact that $\de \theta = \Pi^* \omega$ thus implies that the f\/ibres of $\Phi$ are isotropic. If $\rho_M \colon J^1M \to TM$ and $\rho_S \colon T^*S \to TS$ denote the anchor maps associated to the contact and symplectic structures on $M$ and $S$ respectively, then the fact that $\Pi$ is a Jacobi map and commutativity of the diagram in equation~\eqref{eq:51} imply that, for all $(\alpha,f) \in \Gamma(\nu^* \oplus \RR)$,
 \begin{gather*}
 D \Pi(\rho_M(\phi^*(\alpha,f))) = \rho_S(\Phi^*\alpha) \circ \Pi.
 \end{gather*}
In other words, the vector f\/ields $\rho_M(\phi^*(\alpha,f))$ and $\rho_S(\Phi^*\alpha)$ are $\Pi$-related; therefore, for all $t \in \RR$, $\Pi \circ \varphi^t_{(\alpha,f)} = \varphi^t_{\alpha} \circ \Pi$, where $\varphi^t_{(\alpha,f)}$ and $\varphi^t_{\alpha}$ denote the f\/lows at time $t$ of $\rho_M(\phi^*(\alpha,f))$ and $\rho_S(\Phi^*\alpha)$ respectively. Thus if $(\alpha,f) \in \Gamma(\Sigma)$, then $\alpha$ is a section of the isotropy bundle of the $\nu^*$-action on $\Phi \colon \SM \to \PM$; however, $\alpha = \mathsf{P}(\alpha,f)$, where $\mathsf{P} \colon \nu^* \oplus \RR \to \nu^*$ is the projection of equation~\eqref{eq:47}. Therefore, $\Xi = \mathsf{P}(\Sigma)$ (cf.\ Proposition~\ref{prop:sigma_proj_xi}) is contained in the isotropy $\hat{\Xi}$ of the action $\nu^* \hookrightarrow \SM$ along $\Phi$. To show that $\Xi = \hat{\Xi}$, observe that if $\alpha$ is a local section of $\hat{\Xi}$ def\/ined on a suf\/f\/iciently small domain $U \subset P$, then $\alpha = \de f$ for some locally def\/ined smooth function $f$ on $P$. Then have that $\Pi \circ \varphi^1_{(\de f,f)} = \varphi^1_{\de f} \circ \Pi = \Pi$, so that $\varphi^1_{(\de f,f)}(m) \in \Pi^{-1}(\Pi(m))$ for all $m \in \phi^{-1}(U)$. Since $\ker D\Pi = \RR \langle R_1 \rangle$, there exists a smooth function $\tau \colon \phi^{-1}(U) \to \RR/\ZZ$ such that, for all $m \in
 \phi^{-1}(U)$, $\varphi^1_{(\de f,f)}(m) = \varphi^{\tau(m)}_{(0,1)}(m)$.

 \begin{Claim}\label{claim:basic}
 The function $\tau$ is $\phi$-basic.
 \end{Claim}
\begin{proof}[Proof of Claim \ref{claim:basic}] The above statement is proved below by showing that $\tau$ is constant along the f\/ibres of $\phi$. Fix $m \in \phi^{-1}(U)$ and let $m' \in \phi^{-1}(\phi(m))$. Since the action of $\nu^* \oplus \RR$ on $\phi \colon \CM \to \PM$ is transitive along the f\/ibres (cf.\ Note~\ref{obs:per_lattice}), there exists an element $\xi \in (\nu^* \oplus \RR )_{\phi(m)}$ such that $m' = \varphi^1_{\xi}(m)$. Then
 \begin{gather*} 
 \varphi^{\tau (\varphi^1_{\xi}(m) )}_{(0,1)} \big(\varphi^1_{\xi}(m) \big) = \varphi^1_{(\de f, f)}\circ \varphi^1_{\xi}(m) =
 \varphi^1_{\xi} \circ \varphi^1_{(\de f, f)}(m) = \varphi^1_{\xi} \circ \varphi^{\tau(m)}_{(0,1)}(m) \\
\hphantom{\varphi^{\tau (\varphi^1_{\xi}(m) )}_{(0,1)} \big(\varphi^1_{\xi}(m) \big)}{} = \varphi^1_{\xi} \circ \varphi^{1}_{(0,\tau(m))}(m) =
 \varphi^{1}_{(0,\tau(m))} \big(\varphi^1_{\xi}(m)\big) \\
\hphantom{\varphi^{\tau (\varphi^1_{\xi}(m) )}_{(0,1)} \big(\varphi^1_{\xi}(m) \big)}{} = \varphi^{\tau(m)}_{(0,1)} \big(\varphi^1_{\xi}(m) \big) ,
 \end{gather*}
where the second and fourth equalities follow from the fact that the vector f\/ields $\rho_M(\phi^*\xi)$, $\rho_M(\phi^*(\de f,f))$ and $\rho_M(\phi^*(0,1))$ commute, as $\xi, (\de f,f)_{\phi(m)}, (0,1)_{\phi(m)}$ are elements of $(\nu^* \oplus \RR)_{\phi(m)}$, which is an {\it abelian} Lie algebra.
 \end{proof}

Let $\bar{\tau}$ be the unique smooth function such that $\tau = \phi^*\bar{\tau}$. Since $R_1$ is tangent to the f\/ibres of $\phi$, it follows that $\varphi^{\tau}_{(0,1)} = \varphi^1_{(0,\bar{\tau})}$. Thus, on $\phi^{-1}(U)$, $\varphi^1_{(\de f,f)} \circ \varphi^1_{(0,-\bar{\tau})} = id.$, thus showing that $(\de f, f - \bar{\tau})$ is a local section of $\Sigma$; this immediately implies that $\de f$ is a local section of $\mathsf{P}(\Sigma) =\Xi$, thus proving that $\Xi = \bar{\Xi}$.

\underline{$(\Leftarrow)$}: Conversely, suppose that $\Phi \colon \SM \to \CM$ is a SIR with period lattice $\Xi$, where $\omega$ is an integral symplectic form. Let $\Pi \colon M \to S$ be a principal $S^1$-bundle with Chern class equal to $[\omega] \in \HH^2(S;\ZZ)$; a connection 1-form $\theta \in \Omega^1(M)$ with curvature $\omega$ def\/ines a contact structure $H=\ker \theta$ on $M$ with respect to which $\Pi \colon \CM \to \SM$ is a CIR (cf.\ Example~\ref{sec:motiv-preq-circle}). The aim is to show that the composite $\phi:= \Phi \circ \Pi \colon \CM \to \PM$ is a CIR with period lattice equal to~$\Sigma$. Since $\phi = \Phi \circ \Pi$ and $\Pi$ and $\Phi$ are Jacobi maps, $\phi$ is a Jacobi map. Given that $\ker D \Pi = \RR\langle R_1\rangle$ and that $\ker D \Pi \subset \ker D \phi$, it follows that $\ker D \phi$ is transversal to $H$. Observe that, by def\/inition of~$\theta$, for any $m \in M$, $D_m \Pi|_{H_m}\colon H_m \to T_{\Pi(m)}S$ is an isomorphism of symplectic vector spaces which identif\/ies $\ker D\phi \cap H$ with $\ker D\Phi$. Since the latter is isotropic, so is the former; by Note~\ref{item:po}, it follows that $\phi \colon \CM \to \PM$ is a CIR. It remains to show that the period lattice of~$\phi$ is precisely $\Sigma$. Denote the isotropy of the $\nu^* \oplus \RR$-action on $\phi \colon \CM \to \PM$ by $\hat{\Sigma}$; observe that, by construction, $\ZZ\langle j^11\rangle \subset \hat{\Sigma}$, as the Reeb vector f\/ield $R_1$ is the inf\/initesimal generator of the free and proper $S^1$-action on~$M$. Using notation from $(\Rightarrow)$-part of the proof, have that, for all $(\alpha,f) \in \Gamma(\nu^* \oplus \RR)$, and for all $t \in \RR$,
 \begin{gather} \label{eq:53}
 \Pi \circ \varphi^t_{(\alpha,f)} = \varphi^t_{\alpha} \circ \Pi.
 \end{gather}
Thus, if $j^1f \in \Gamma_{\mathrm{loc}}(\hat{\Sigma})$, then $\mathsf{P}(j^1f) = \de f \in \Gamma(\Xi)$. Hence $\hat{\Sigma} \subset \Sigma$. To show that, in fact $\hat{\Sigma} = \Sigma$, the idea is to argue as above. Consider a local section $\xi \in \Sigma$; since $\mathsf{P}(\Sigma) = \Xi$ and by equation~\eqref{eq:53}, there is a~locally def\/ined smooth function $\tau$ such that $\varphi^1_{\xi} = \varphi^{\tau}_{(0,1)}$. Arguing as in Claim~\ref{claim:basic}, there exists a locally def\/ined smooth function $\bar{\tau}$ on $P$ with $\tau = \phi^*\bar{\tau}$. Thus, $(0,\bar{\tau}) \in \Gamma(\nu^* \oplus \RR)$ and $\xi - (0,\bar{\tau}) \in \Gamma(\hat{\Sigma})$. Using the fact $\hat{\Sigma} \subset \Sigma$, obtain that $\bar{\tau} \in \ZZ$, which therefore implies that $\xi \in \Gamma(\hat{\Sigma})$, since $\ZZ\langle j^11\rangle \subset \hat{\Sigma}$. This concludes the proof.
\end{proof}

\begin{Note}\label{obs:interpretation}
Suppose that $\PM$ admits a CIR with period lattice $\Sigma$, say with Chern class $c$; commutativity of the diagram of equation \eqref{eq:48} implies that $\mathcal{P}(c)$ is the Chern class of a SIR of $\PM$ with period lattice~$\Xi$. However, Theorem~\ref{thm:cir_vs_sir} shows that the symplectic form on the total space of this SIR can be chosen to be integral, something that, to the best of our knowledge, cannot be inferred simply from the commutativity of the diagram~\eqref{eq:48}. Conversely, suppose that $c' \in \HH^2(P;\underline{\Xi})$ is the Chern class of a SIR of $\PM$ with period lattice $\Xi$. Commutativity of the diagram of equation~\eqref{eq:48} implies that a necessary condition for the existence of a CIR of $\PM$ with period lattice $\Sigma$ is that $\delta(c') = 0$, as the characteristic class $\upsilon = \mathfrak{d}(c')$ is def\/ined by an exact 3-form, thus being mapped to zero under $\HH^3_{\mathrm{rel}}(P;\mathcal{F}) \to \HH^3(P;\RR)$. This is, however, a very mild restriction: if $\HH^3(P;\ZZ)$ is torsion-free, $\delta(c') =0$ is automatic, as the vertical map $\HH^3(P;\ZZ) \to \HH^3(P;\RR)$ in equation \eqref{eq:48} is the natural inclusion. If, on the one hand, $\delta(c') =0$ implies that there exists $c \in \HH^2(P;\underline{\Sigma})$ with $\mathcal{P}(c) = c'$ by exactness of the top row of \eqref{eq:48}, on the other, this is not suf\/f\/icient to conclude that $c$ is the Chern class of a CIR of $\PM$ with period lattice $\Sigma$. Intuitively, what is needed is that the foliated 2-form $\omega_{\mathcal{F}}$ be, in some sense, {\it integral}. For instance, if $\upsilon = 0$, i.e., there exists a globally def\/ined closed 2-form $\Omega \in \Omega^2(P)$ extending $\omega_{\mathcal{F}}$, then integrality of $\Omega$ suf\/f\/ices to guarantee, together with $\delta(c') = 0$, that there exists a CIR of $\PM$ with period lattice $\Sigma$.
\end{Note}

\begin{Example}\label{exm:cpc_ss_Lie_algebras} (For details about this example, see \cite[Section~4.5.1]{PMCT2}.) Let $G$ be a compact, simply connected Lie group and set $\mathfrak{g} = \operatorname{Lie}(G)$. The subset $\mathfrak{g}^*_{\mathrm{reg}} \subset \mathfrak{g}^*$ consisting of coadjoint orbits whose stabiliser is a maximal torus is a regular Poisson manifold, whose bracket is the restriction of the linear bracket on $\mathfrak{g}^*$. In fact, if $\mathbb{T} \subset G$ is a maximal torus with $\mathfrak{t} = \operatorname{Lie}(\mathbb{T})$, and $\mathfrak{c} \subset \mathfrak{t}^*$ is the interior of a Weyl chamber, then there is a dif\/feomorphism $\mathfrak{g}^*_{\mathrm{reg}} \cong G/\mathbb{T} \times \mathfrak{c}$, which identif\/ies $\mathfrak{c}$ as the leaf space of the symplectic foliation induced by the above Poisson structure. The kernel of $\exp \colon \mathfrak{t} \to \mathbb{T}$ induces a~$\ZZ$-af\/f\/ine structure on $\mathfrak{t}^*$, which, upon suitable identif\/ications, corresponds to the standard $\ZZ$-af\/f\/ine structure on $\RR^l \cong \mathfrak{t}^*$, where $l = \operatorname{rk} G$. This $\ZZ$-af\/f\/ine structure is intimately connected to the Poisson geometry of $\mathfrak{g}_{\mathrm{reg}}^*$ (cf.\ \cite[Section~4.5.1]{PMCT2}). Fix $\xi_0 \in \mathfrak{c}$ and denote by $S_0 \subset \mathfrak{g}^*_{\mathrm{reg}}$ the coadjoint orbit through $\xi_0$. Viewing $\mathfrak{c} $ as an open subset of $ \mathfrak{t}^*$, obtain a strong transversal $\ZZ$-af\/f\/ine lattice $\Xi$ on~$\mathfrak{g}^*_{\mathrm{reg}}$ and denote by $\Sigma$ the transversal $\ZZ$-projective lattice on $\mathfrak{g}^*_{\mathrm{reg}}$ induced as in Corollary~\ref{cor:strong_Z_proj}. If $\xi^1,\ldots,\xi^l\colon \mathfrak{c} \to \RR$ denote $\ZZ$-af\/f\/ine coordinates on $\mathfrak{c}$, then $\Xi = \ZZ \langle \de \xi^1,\ldots \de \xi^l\rangle $ and $\Sigma = \ZZ \langle j^1\xi^1,\ldots,j^1\xi^l,j^11\rangle$, where, by abuse of notation, $\xi^1,\ldots \xi^l$ are seen as functions on $\mathfrak{g}^*_{\mathrm{reg}}$. On the other hand, since $G$ is compact, if~$\mathcal{F}$ denotes the symplectic foliation on $\mathfrak{g}^*_{\mathrm{reg}}$, then $\HH^2(\mathcal{F}) \cong \HH^2(S_0;\RR) \otimes C^{\infty}(\mathfrak{c})$ (cf.\ \cite{kacimi}). If $c_1,\ldots,c_l \in \HH^2(S_0;\ZZ)$ are torsion-free generators, then \cite[Remark~4.3.5]{PMCT2} implies that the foliated class of the foliated symplectic form
 $\omega_{\mathcal{F}}$ is given by
 \begin{gather*}
 [\omega_{\mathcal{F}}] = [\omega_0] + \sum\limits_{i=1}^l \xi^i \otimes c_i,
 \end{gather*}
where $[\omega_0] \in \HH^2(S_0;\RR)$ is the cohomology class of the symplectic form on $S_0$ (cf.\ \cite[Section~4.5.1]{PMCT2} for a reason why $\operatorname{rk} \HH^2(S_0;\ZZ) = \operatorname{rk} G$). Seeing as $S_0$ can be chosen so that $[\omega_0]$ is integral, Example~\ref{exm:trivial_tips} implies that $\mathfrak{g}^*_{\mathrm{reg}}$ admits a CIR with period lattice $\Sigma$ and, equivalently by Theorem~\ref{thm:cir_vs_sir}, a SIR with period lattice $\Xi$ whose total space has an integral symplectic form.

Observe that choosing $S_0$ so that it has an integral symplectic form is equivalent to having $\left(\xi^1(\xi_0),\ldots, \xi^l(\xi_0)\right) \in \ZZ^l$. Let $U \subset \mathfrak{c}$ be an open, connected subset not containing any integral point (with respect to the above strong $\ZZ$-af\/f\/ine structure!) and consider the open Poisson submanifold $(P_U,\Lambda_U) \subset \mathfrak{g}^*_{\mathrm{reg}}$ obtained by considering the union of coadjoint orbits corresponding to points in~$U$. Denote the induced transversal strong $\ZZ$-af\/f\/ine and $\ZZ$-projective lattices by $\Xi_U$ and $\Sigma_U$ respectively. Observe that $\Sigma_U$ equals the transversal $\ZZ$-projective lattice obtained from $\Xi_U$ as in Corollary~\ref{cor:strong_Z_proj}. Example~\ref{exm:trivial_tips} implies that $\left(P_U,\Lambda_U\right)$ admits no CIR with period lattice~$\Sigma_U$, but it admits a SIR with period lattice $\Xi_U$ as it can be readily verif\/ied using \cite[Section~4]{dd}.
\end{Example}

\appendix
\section[Properties of regular Jacobi manifolds all of whose leaves are even dimensional]{Properties of regular Jacobi manifolds\\ all of whose leaves are even dimensional}\label{sec:proofs-prop-regul}

\begin{Lemma}\label{prop:well-defined} Let $(P,L,\{\cdot, \cdot\})$ be a regular Jacobi manifold with even dimensional leaves. The connection $\bar\nabla$ of equation \eqref{eq:26} is a well-defined flat $T \mathcal{F}$-connection.
\end{Lemma}

\begin{proof}Assuming that $\bar\nabla$ is well-def\/ined, f\/latness of $\bar\nabla$ follows directly from f\/latness of~$\nabla$. It suf\/f\/ices to show that $\nabla$ restricted to $\ker\rho$ is zero. For this purpose, some useful properties of the anchor and $\nabla$ are derived to compute $\nabla_\alpha$ for $\alpha\in\Gamma(\ker\rho)$.

First it is shown that that the restriction of $\rho$ to $T^*P\otimes L$ is antisymmetric, i.e.,
\begin{gather*}\zeta(\rho(\eta))=-\eta(\rho(\zeta))\end{gather*}
for any two elements $\zeta,\eta\in\Omega^{1}(P;L)$. For $f\in C^\infty(P)$ and $u\in\Gamma(L)$, write $df\otimes u=fj^1u-j^1(fu)$. The def\/ining property~\ref{item:A} of the anchor shows that $L_{\rho(df\otimes u)}(g)v=-L_{\rho(dg\otimes v)}(f)u$ for any other $g\in C^\infty(P)$ and any $v\in\Gamma(L)$. This shows the antisymmetry and also shows that if $\nu^*\subset \cotan P$ is the conormal bundle of the regular foliation $\mathcal{F}$, then $\nu^*\otimes L\subset\ker\rho$, as if $f$ is constant along the leaves of $\mathcal{F}$, then $L_{\rho(df\otimes u)}(g)v=-L_{\rho(dg\otimes v)}u=0.$

Next, it is shown that, in fact,
 \begin{gather*}\nu^*\otimes L=\ker\rho\cap \big(T^*P\otimes L\big).\end{gather*}
This can by checked by dimension counting: as the corank of the Jacobi manifold is $k=\dim\ker\rho$, a straightforward computation shows that the rank of $\nu^*$ is $k-1$. Hence, $k-1\leq \dim\ker\rho(T^*P\otimes L)\leq k$, or equivalently $ 2n - 1 \leq \dim \rho(\cotan P \otimes L) \leq 2n. $ On the other hand, the antisymmetry of $\rho\colon T^*P\otimes L\to TP$ implies that the subspace $\rho(T^*P\otimes L)\subset\mathcal{F}$ is even dimensional. Hence, $\dim \rho(\cotan P \otimes L)=2n$, and $\dim\left(\ker\rho\cap (T^*P\otimes L)\right)=k-1.$

A straightforward computation shows that, for any $f \in C^{\infty}(P)$ and $u,v \in \Gamma(L)$, $\nabla_{df\otimes u}(v)=-df(\rho(j^1v))u$, hence
\begin{gather*}\nabla_\zeta(v)=-\zeta\big(\rho\big(j^1v\big)\big), \qquad \zeta\in\Omega^1(P;L).\end{gather*}
With the above properties in mind, write $\alpha$ as $\alpha=(u,\zeta)\in\Gamma(L)\oplus\Omega^1(P;L)$. If $\alpha$ is a section of $\ker\rho\subset J^1L$, then $\nabla_\alpha\colon \Gamma(L)\to\Gamma(L)$ is $C^{\infty}(M)$-linear, hence it def\/ines a section of $\hom(L;L)$. As $\hom(L;L)$ is trivial, there exists a function $G\in C^\infty(P)$ with the property that
\begin{gather*}\nabla_\alpha(v)(p)=G(p)v(p),\end{gather*}
for any $v\in\Gamma(L)$ and any $p\in P$. With this,
\begin{gather*}G(p)u(p)=\nabla_\alpha(u)(p)=\{u,u\}-\zeta\big(\rho\big(j^1u\big)\big)=\zeta(\rho(\zeta))=0,\end{gather*}
where in the third equation we use that $\rho(j^1u-\zeta)=0$, and in the last one the antisymmetry of~$\rho$. If $u(p)\neq 0$, then $G(p)=0$, and $\nabla_\alpha(p)=0$; otherwise $u(p)=0$, $\alpha(p)\in (T^*P\otimes L)_p\cap\ker\rho_p=\nu^*_p\otimes L_p$ and
\begin{gather*}\nabla_\alpha(v)(p)=\alpha(p)\big(\rho\big(j^1_pv\big)\big)=0.\tag*{\qed}\end{gather*}\renewcommand{\qed}{}
\end{proof}

\begin{Lemma}\label{lemma:abelian}
If $(P,L,\{,\})$ is a regular Jacobi manifold all of whose leaves are even dimensional then $\ker\rho\subset J^1L$ is bundle of abelian Lie algebras.
\end{Lemma}

\begin{proof} Using the Spencer decomposition of Note \ref{obs:prop_cso}, the Lie bracket of two sections $\alpha,\beta\in\Gamma(J^1L)$ equals $(\pr([\alpha,\beta]),\mathrm{D}[\alpha,\beta])$. If $\alpha$, $\beta$ are sections of $\ker\rho\subset J^1L$, the compatibility conditions~\eqref{horizontal} and~\eqref{vertical} imply that
 \begin{gather*}
 [\alpha,\beta]=\big({-}\nabla_\beta(\pr(\alpha)), \nabla_\alpha (\mathrm{D}_\cdot(\beta))-\nabla_\beta(D_{\cdot}(\alpha))\big).
 \end{gather*}
The proof of Lemma~\ref{prop:well-defined} shows that $\nabla$ restricted to $\ker\rho$ is zero, thus proving the result.
\end{proof}

The f\/lat connection $\bar\nabla$ induces the usual Koszul-like dif\/ferential $\de_{\mathcal{F}}$ on the complex of smooth foliated forms of $\mathcal{F}$ with values in $L$. The 2-form $\omega_{\mathcal{F}}\in\Omega^2(\mathcal{F};L)$ of equation \eqref{eq:29} def\/ines a~canonical cohomology class in $\HH^2(\mathcal{F};L)$, as the next lemma shows.

\begin{Lemma}\label{lemma:coho_class} Under the hypotheses of Lemma~{\rm \ref{prop:well-defined}}, $\omega_{\mathcal{F}}$ is well defined, and $\de_{\mathcal{F}} \omega_{\mathcal{F}} = 0$.
\end{Lemma}

\begin{proof} Since sections of the form $j^1 u$, for $u \in \Gamma(L)$ form a $C^{\infty}(P)$-basis of $\Gamma(J^1 L)$ and $T \mathcal{F} = \rho(J^1 L)$, equation~\eqref{eq:29} def\/ines a unique map $\Gamma(T \mathcal{F}) \times \Gamma(T \mathcal{F}) \to \Gamma(L)$ which is manifestly antisymmetric. Let $f \in C^{\infty}(P)$ and f\/ix $u,v \in \Gamma(L)$; then
 \begin{gather*}
 \omega_{\mathcal{F}}\big(f\rho\big(j^1 u\big), \rho\big(j^1 v\big)\big) = \omega_{\mathcal{F}}\big(\rho(fu,\de f \otimes u), \rho\big(j^1v\big)\big) = \{fu,v\} + \de f \otimes u\big(\rho\big(j^1v\big)\big) \\
\hphantom{\omega_{\mathcal{F}}\big(f\rho\big(j^1 u\big), \rho\big(j^1 v\big)\big)}{}
=f\{u,v\} - \mathcal{L}_{\rho(j^1v)}f u + \mathcal{L}_{\rho(j^1v)}f u = f \omega_{\mathcal{F}}\big(\rho\big(j^1 u\big), \rho\big(j^1 v\big)\big),
\end{gather*}
where the f\/irst equality follows from the $C^{\infty}(P)$-structure on $\Gamma(J^1 L)$ arising from the Spencer decomposition (cf.\ Note~\ref{obs:prop_cso}), the second by def\/inition of $\omega_{\mathcal{F}}$, and the third by the characterising property~\ref{item:A} of the anchor~$\rho$ (cf.\ Proposition~\ref{prop:alg}). The above calculation shows that equation~\eqref{eq:29} indeed def\/ines a foliated 2-form.

 To check that $\de_{\mathcal{F}}\omega_{\mathcal{F}} =0$ it suf\/f\/ices to check that
 \begin{gather*}
 \de_{\mathcal{F}}\omega_{\mathcal{F}}\big(\rho\big(j^1u\big),\rho\big(j^1v\big),\rho\big(j^1w\big)\big) =0
 \end{gather*}
for any $u,v,w \in \Gamma(L)$. This is because $\de_{\mathcal{F}}\omega_{\mathcal{F}}$ is $C^{\infty}(P)$-linear in each entry and sections of the form $\rho(j^1u)$ form a $C^{\infty}(P)$-basis of $\Gamma(T \mathcal{F})$. Then
 \begin{gather*}
\de_{\mathcal{F}}\omega_{\mathcal{F}}\big(\rho\big(j^1u\big),\rho\big(j^1v\big),\rho\big(j^1w\big)\big)\\
\qquad {} = \bar{\nabla}_{\rho(j^1u)}\big(\omega_{\mathcal{F}}\big(\rho\big(j^1v\big),\rho\big(j^1w\big)\big)\big)+ \mathrm{c.p.}
 - \big(\omega_{\mathcal{F}}\big(\big[\rho\big(j^1u\big),\rho\big(j^1v\big)\big],\rho\big(j^1w\big)\big) + \mathrm{c.p.}\big)\\
 \qquad{} = \{u,\{v,w\}\} + \mathrm{c.p.} - \big(\omega_{\mathcal{F}}\big(\rho\big(\big[j^1u,j^1v\big]\big), \rho\big(j^1w\big)\big) + \mathrm{c.p.}\big) \\
 \qquad {} = -\big(\omega_{\mathcal{F}}\big(\rho\big(j^1\{u,v\}\big),\rho\big(j^1w\big)\big) + \mathrm{c.p.}\big) =-(\{\{u,v\},w\} + \mathrm{c.p.}) = 0,
 \end{gather*}
where c.p.\ stands for cyclic permutation, the second equality follows from the fact that $\rho$ is a~map of Lie algebroids, the third by the def\/ining property~\ref{item:B} of the Lie bracket on $J^1 L$ (cf.\ Proposition~\ref{prop:alg}) and by the Jacobi identity for $\{\cdot,\cdot\}$, which also implies the last equality. This shows that $\de_{\mathcal{F}}\omega_{\mathcal{F}} = 0$, as required.
\end{proof}

Using the above ideas, we can prove Note \ref{obs:sod}.

\begin{proof}[Proof of Note \ref{obs:sod}]
Since $\Gamma(T\mathcal{F})$ is generated by elements of the form $\rho(j^1u)$ for $u \in \Gamma(J^1L)$, it suf\/f\/ices to check that $\mathrm{D}_{\rho(j^1u)}(\alpha) = \de_{\mathcal{F}}(\pr(\alpha))(\rho (j^1u))$ for any $u \in \Gamma(J^1L)$. Fix such a~section. Then
\begin{gather*}
 \mathrm{D}_{\rho(j^1u)}(\alpha) = \nabla_{\alpha} u + \pr\big(\big[j^1u,\alpha\big]\big) = \pr\big(\big[j^1u,\alpha\big]\big) = - \pr\big(\big[\alpha,j^1u\big]\big) \\
\hphantom{\mathrm{D}_{\rho(j^1u)}(\alpha)}{} = \nabla_{j^1u}(\pr(\alpha)) - \mathrm{D}_{\rho(\alpha) }\big(j^1u\big)
 = \{u,\pr(\alpha)\} = \de_{\mathcal{F}} (\pr(\alpha) )\big(\rho\big(j^1u\big)\big),
 \end{gather*}
where the f\/irst and fourth equalities use the compatibility condition~\eqref{horizontal} for the Spencer operator, the second uses that $\alpha \in \Gamma(\ker \rho)$, the third exploits anti-symmetry of the Lie bracket, and the last two follow by def\/inition of $\nabla$ (cf.\ equation~\eqref{eq:nabla}).
\end{proof}

\section{Proofs of the main results of Section \ref{sec:motiv-famil-exampl}}\label{sec:proofs-sect-refs}

\begin{proof}[Proof of Proposition \ref{lemma:mu_misses_zero}] Suppose that the statement does not hold. Then there exists $\alpha \in \SYM \subset T^* M$ with $\mu(\alpha) = 0$. Recall that $\mu$ is the restriction of the moment map of the cotangent lift which preserves the Liouville 1-form on~$T^*M$ (cf.\ Note~\ref{obs:contact_action_lift}); unravelling the def\/initions, it follows that, for any $\xi \in \mathfrak{g}$,
\begin{gather*}
\langle \mu(\alpha), \xi \rangle = \alpha\big(D_{\alpha}\operatorname{pr} (\xi_{\SYM}(\alpha) )\big),
 \end{gather*}
where $\xi_{\SYM}(\alpha)$ is the image of $\xi$ under the inf\/initesimal action of $G$ on $\SYM$ evaluated at~$\alpha$. Since $\operatorname{pr} \colon \SYM \to M$ is $G$-equivariant,
 \begin{gather*}
 D_{\alpha}\operatorname{pr} (\xi_{\SYM}(\alpha) ) = \xi_M(p),
 \end{gather*}
where $p = \operatorname{pr}(\alpha)$. Thus, for any $\xi \in \mathfrak{g}$, $\alpha(\xi_M(p)) = 0$, which implies that $\xi_M(p) \in H_p$, in turn yielding that $T_p ( G \cdot p ) \subset H_p$.

Let $\theta$ be a contact 1-form def\/ined on a $G$-invariant neighbourhood $U$ of $p$, i.e., $H = \ker \theta$ locally and $\de \theta|_H$ is symplectic. By \cite[Lemma~2.6]{lerman_contact_toric}, it may be assumed that $\theta$ is $G$-invariant. For any $\xi, \xi' \in \mathfrak{g}$,
 \begin{gather*}
 \de \theta_p\big(\xi_M(p), \xi_M'(p) \big) = \theta_p \big([\xi_M , \xi_M'](p) \big) = \theta_p\big([\xi, \xi']_M (p)\big) = 0,
 \end{gather*}
 where the f\/irst and third equalities follow from the fact that $T_p(G \cdot p) \subset H_p$. Therefore $T_p(G \cdot p )$ is an isotropic subspace of $(H_p,\de \theta_p)$.

The 1-form $\theta$ determines a trivialisation of $L^*|_U \cong U \times \RR$. Henceforth, identify any $\alpha' \in L^*|_U$ with the pair $(p',t') \in U \times \RR$, where $\operatorname{pr} (\alpha' ) = p '$ and $\alpha' = t' \theta_p$. Since $\theta$ is $G$-invariant, the induced $G$-action on $L^*|_U$ is given in this trivialisation by
 \begin{gather*}
 g \cdot (p', t') = (g \cdot p', t').
 \end{gather*}
Setting $\alpha = (p,t)$ (where $t \neq 0$), $G$-invariance of $\theta$ implies that there is a splitting as $G_{\alpha} = G_p$-symplectic vector spaces
 \begin{gather*} 
 \big(T_{\alpha} ( \SYM ),\Omega_{\alpha}\big) = (H_p, t\de \theta_p ) \oplus \left(\RR \left\langle R_p,
 \frac{\partial}{\partial t} \right\rangle , \de t \wedge \theta_p\right)=: (H_p,t \de \theta_p) \oplus (V, \omega_V ),
 \end{gather*}
where $R$ is the locally def\/ined Reeb vector f\/ield associated to the contact 1-form $\theta$, i.e., $\theta(R) = 1$ and $\de \theta(R,-) = 0$, and the $G_{\alpha}$-action on $V$ is trivial. Moreover,
 \begin{gather} \label{eq:17}
 T_{\alpha} (G \cdot \alpha) = T_{(p,t)} (G \cdot (p,t) ) = T_p (G \cdot p ) \oplus 0;
 \end{gather}
since $T_{p} (G \cdot p) \subset (H_p, \de \theta_p)$ is isotropic, equation \eqref{eq:17} implies that $T_{\alpha} ( G \cdot \alpha ) \subset (T_{\alpha} ( \SYM ),\Omega_{\alpha})$ is isotropic. In fact,
 \begin{gather*}
(T_{\alpha} (G \cdot \alpha))^{\Omega_{\alpha}} = (T_{p} (G \cdot p ))^{\de \theta_p} \oplus V.
 \end{gather*}
Setting
 \begin{gather*}
 N_{\alpha}:= \frac{(T_{\alpha} (G \cdot \alpha))^{\Omega_{\alpha}}}{T_{\alpha} (G \cdot \alpha)}\qquad \text{and} \qquad N_p:= \frac{(T_{p} (G \cdot p))^{\de\theta_p}}{T_{p} (G \cdot p)},
 \end{gather*}
obtain that, as $G_{\alpha}=G_p$-symplectic vector spaces,
 \begin{gather} \label{eq:24}
(N_{\alpha}, \omega_{N_{\alpha}}) \cong (N_{p}, \omega_{N_{p}})\oplus (V,\omega),
 \end{gather}
where the $G_p$-action on the right hand side is the product of the induced one on the f\/irst factor and the trivial one on the second.

Since $G$ is compact, the Marle--Guillemin--Sternberg local normal form for Hamiltonian actions holds (cf.~\cite{guillemin_sternberg,marle_ham}); as $T_{\alpha} ( G \cdot \alpha )$ is isotropic, this reduces to the following. Set $\mu(\alpha) = \eta \in \mathfrak{g}^*$; $G$-equivariance of $\mu$ gives the following inclusion of stabilisers $G_{\alpha} \subset G_{\eta}$, where $G \curvearrowright \mathfrak{g}^*$ acts by the coadjoint action. This induces an analogous inclusion for the corresponding Lie algebras $\mathfrak{g}_{\alpha} \subset \mathfrak{g}_{\eta}$. Fix a $G_{\alpha}$-invariant inner product on~$\mathfrak{g}$, which induces a $G_{\alpha}$-invariant orthogonal direct sum decomposition $\mathfrak{g}_{\eta} = \mathfrak{g}_{\alpha} \oplus \mathfrak{m}$. Then a $G$-invariant neighbourhood of $G\cdot \alpha$ in $\SYM$ is isomorphic, as a Hamiltonian $G$-space, to a $G$-invariant neighbourhood of $G \cdot [e,0,0]$ in
 \begin{gather*}
(Y,\omega_Y) : = \big(G \times_{G_{\alpha}} (\mathfrak{m}^* \times N_{\alpha}), \omega_Y \big),
 \end{gather*}
where the right hand side is, as a smooth manifold, the quotient of $G \times ( \mathfrak{m}^* \times N_{\alpha})$ by the antidiagonal action of~$G_{\alpha}$, and the action of $G$ on $Y$ descends from the left action of $G$ on the f\/irst factor of $G \times (\mathfrak{m}^* \times N_{\alpha})$. Following~\cite{ortega_ratiu,sl}, the closed 2-form $\omega_Y$ (which is symplectic only near $G \cdot [e,0,0]$) can be constructed using Marsden-Weinstein reduction on the presymplectic manifold $(G \times \mathfrak{g}^*_{\eta} \times N_{\alpha}, \sigma \oplus \omega_{\alpha})$ endowed with the induced
 Hamiltonian $G_{\alpha}$-action, where $\sigma$ is a closed 2-form on $G \times \mathfrak{g}^*_{\eta}$ constructed as in~\cite{ortega_ratiu}. By the splitting of equation~\eqref{eq:24}, it follows that, as Hamiltonian a $G_{\alpha}$-space, $(G \times \mathfrak{g}^*_{\eta} \times N_{\alpha}, \sigma \oplus \omega_{\alpha})$ equals
 \begin{gather*}
 \big(G \times \mathfrak{g}^*_{\eta} \times N_{p} \times V, \sigma \oplus \omega_{\alpha} \oplus \omega_V \big),
 \end{gather*}
 where the $G_{\alpha}$-action on the last factor is trivial. Unravelling the construction of in~\cite{ortega_ratiu,sl}, it follows that, as a~presymplectic
 Hamiltonian $G$-space,
 \begin{gather*}
 (Y,\omega_Y) \cong \big(W \times V, \omega_W \oplus \omega_V\big):=\big(\big(G \times_{G_{\alpha}} (\mathfrak{m}^* \times N_{p} )\big) \times V, \omega_W \oplus \omega_V \big),
 \end{gather*}
where the $G$-action on $V$ is trivial; this also implies that the $G$-action on the presymplectic manifold $(W,\omega_W)$ is Hamiltonian. Restricting to a $G$-invariant neighbourhood of the image of $G \cdot \alpha$ on which $\omega_W \oplus \omega_V$ is symplectic, the Hamiltonian action of $G$ is multiplicity-free; it follows, in particular, that the action of $G$ on $W$ is locally free at some point $w \in W$; by the principal orbit theorem (cf.\ \cite[Theorem~2.8.5]{duistermaat_kolk}), there is an open dense set of~$W$ whose stabiliser is discrete. Seeing as $(W,\omega_W)$ is symplectic on a~$G$-invariant neighbourhood $U$ of $G\cdot[e,0,0]$, it follows that the Hamiltonian $G$-action on $(U, \omega_W|_U)$ is locally free at some point in~$U$, which implies that $\dim U = \dim W \geq \dim G + \operatorname{rk} G$ by \cite[Theorem~5.1.6]{gui_sja}. However, $\dim W = \dim Y - \dim V = \dim \SYM - 2$; this contradicts the fact that the action of~$G$ on~$\CM$ is multiplicity-free, thus completing the proof.
\end{proof}

\begin{proof}[Proof of Theorem \ref{lemma:properties_Jac_mm}] First it is shown that $F_{\phi} \colon \phi^*(O(1)) \to L$ as def\/ined in the statement is, in fact, a vector bundle isomorphism. Consider $H_{\phi} \colon L \to \phi^*(O(1))$ def\/ined by
 \begin{gather*}
 u \in L_p \mapsto (p,\eta_u),
 \end{gather*}
where $\eta_u$ is def\/ined as follows. Let $\alpha \in L^*_p \setminus \{0\}$; then $\eta_u \colon \RR \langle \mu_u(\alpha)\rangle \to \RR$ is given by $t\mu_u(\alpha) \mapsto t \langle \alpha , u \rangle$. It is straightforward to check that the def\/inition of~$\eta_u$ does not depend on the choice of $\alpha \in L^*_p \setminus \{0\}$. Smoothness of $H_{\phi}$ can be checked by writing the map explicitly once a local trivialisation of $\operatorname{pr} \colon L^*\setminus \{0\} \to M$ is f\/ixed. Fibrewise linearity of~$H_{\phi}$ follows immediately from the def\/inition, while its f\/ibrewise injectivity follows from the fact that $\alpha \in L^*_p \setminus \{0\}$. Seeing as~$H_{\phi}$ is a map between line bundles, it follows that it is an isomorphism. The map $F_{\phi}$ is easily checked to be its inverse, therefore proving that it is a~vector bundle isomorphism. Checking that $\phi$ is Jacobi with bundle component~$F_{\phi}$ can be done as follows. Unravelling the various def\/initions, obtain that, for all $s \in \Gamma(O(1))$, following equality of maps on $O(1)$,
 \begin{gather} \label{eq:34}
 F_{\operatorname{pr}} \circ \operatorname{pr}^* \circ F_{\phi} \circ \phi^*(s) = \mu^* \circ F_{\pi} \circ \pi^*(s),
 \end{gather}
where $\pi \colon \mathfrak{g}^* \setminus \{0\} \to \pg$ is the projection, and $F_{\operatorname{pr}} \colon \operatorname{pr}^*L \to \RR_{L^* \setminus \{0\}}$ and $F_{\pi} \colon \pi^*(O(1)) \to \RR_{\mathfrak{g}^* \setminus \{0\}}$ are the isomorphisms def\/ined in Examples~\ref{exm:symplectisation} and~\ref{exm:proj_dual_lie_algebra} respectively. Since $\mu$ is Jacobi and $\pi$ is Jacobi with bundle component~$F_{\pi}$, it follows that, for any $s_1,s_2 \in \Gamma(O(1))$,
 \begin{gather*}
 \mu^* \circ F_{\pi} \circ \pi^* \{s_1,s_2\} = \{\mu^* \circ F_{\pi} \circ \pi^*(s_1), \mu^* \circ F_{\pi} \circ \pi^*(s_2)\};
 \end{gather*}
equation \eqref{eq:34} implies that
 \begin{gather}
 F_{\operatorname{pr}} \circ \operatorname{pr}^* \circ F_{\phi} \circ \phi^*\{s_1,s_2\} = \{ F_{\operatorname{pr}} \circ \operatorname{pr}^* \circ F_{\phi} \circ \phi^*(s_1), F_{\operatorname{pr}} \circ \operatorname{pr}^* \circ F_{\phi} \circ \phi^*(s_2)\} \nonumber\\
\hphantom{F_{\operatorname{pr}} \circ \operatorname{pr}^* \circ F_{\phi} \circ \phi^*\{s_1,s_2\}}{}
=F_{\operatorname{pr}} \circ \operatorname{pr}^* \{F_{\phi} \circ \phi^*(s_1), F_{\phi} \circ \phi^*(s_2)\},\label{eq:35}
\end{gather}
where the last equality follows from the fact that $\operatorname{pr}$ is Jacobi with bundle component $F_{\operatorname{pr}}$. Equation~\eqref{eq:35} implies that $\phi$ is Jacobi with bundle component $F_{\phi}$ by noticing that $F_{\operatorname{pr}}$ is an isomorphism and that $\operatorname{pr}^*$ is injective.

Henceforth, f\/ix a point $p \in M_{\mathrm{prin}}$. Since $\operatorname{pr} \colon L^* \setminus \{0\} \to M$ is $G$-equivariant, the stabiliser of any $\alpha \in L_p^* \setminus \{0\}$ is also discrete. Fix such an $\alpha$. By \cite[Corollary~5.1.2]{gui_sja}, $D_{\alpha} \mu$ is surjective; since $\pi \circ \mu = \phi \circ \operatorname{pr}$ and $\pi$ is a submersion, property~\ref{item:5} follows. Suppose that property~\ref{item:6} does not hold, then, since $H \subset TM$ has codimension 1, it follows that $\ker D_p \phi \subset H_p$. First, observe that $D_{\alpha} \operatorname{pr}$ restricted to $\ker D_{\alpha} \mu$ is an isomorphism onto $\ker D_p \phi$: one inclusion follows from $\pi \circ \mu = \phi \circ \operatorname{pr}$, and the other by dimension counting and the fact that property~\ref{item:5} holds. In fact, since $\operatorname{pr}\colon L^* \setminus \{0\} \to M$ is an $\RR^*$-principal bundle, write
 \begin{gather*}
 T_{\alpha} L^* \setminus \{0\} = T_p M \oplus \RR \left\langle \frac{\partial}{\partial t} \right\rangle = H_p \oplus \RR\langle R_p
 \rangle \oplus \RR \left\langle \frac{\partial}{\partial t} \right\rangle,
 \end{gather*}
where $\alpha(R_p) = 1$ and $\de \alpha(R_p,-) = 0$, and $D_{\alpha}\operatorname{pr}$ is simply projection onto the f\/irst two components. Since by assumption $\ker D_p \phi \subset H_p$, it follows that $\ker D_{\alpha} \mu \cap 0 \oplus \RR\langle R_p \rangle \oplus 0 = \{0\}$. Since $\ker D_{\alpha} \mu = ( T_{\alpha}(G \cdot \alpha))^{\Omega}$ by \cite[Corollary~5.1.4]{gui_sja}, it follows that $\RR \langle \frac{\partial}{\partial t} \rangle \subset T_{\alpha}(G \cdot \alpha)$, which is a~contradiction, since $\dim G \cdot \alpha = \dim G \cdot p$ as both stabiliser subgroups are discrete. Therefore~\ref{item:6} holds.

It remains to prove condition \ref{item:14}. First, it is shown that $\ker D_{\alpha}\mu$ is isotropic. Since $\mu $ is $G$-equivariant, the range of $D_{\alpha} \mu$ restricted to $T_{\alpha}(G \cdot \alpha)$ is $T_{\mu(\alpha)} (G \cdot \mu(\alpha))$. The dimension of the latter is at most $\dim G - \operatorname{rk} G$, while the dimension of $D_{\alpha} \mu (T_{\alpha} (G \cdot \alpha))$ is at most $\dim G - \dim (\ker D_{\alpha}\mu) = \dim G - (\dim (L^* \setminus \{0\}) - \dim G) = \dim G - \operatorname{rk} G$, where the f\/irst equality follows from the fact that $D_{\alpha}\mu$ is surjective, and the second from def\/inition of multiplicity-free action. Therefore $\ker D_{\alpha} \mu =(T_{\alpha} (G \cdot \alpha))^{\Omega} \subset T_{\alpha}(G \cdot \alpha)$, which proves that $\ker D_{\alpha} \mu$ is isotropic. Using the fact that $D_{\alpha} \mu$ is surjective, this condition is equivalent to the following inclusion (suppressing the dependence on~$\alpha$ to simplify notation)
\begin{gather} \label{eq:36}
 \ker D \mu \subset \rho_{L^{*} \setminus \{0\}} \big(\mu^*J^1(\mathfrak{g}^* \setminus \{0\})\big),
\end{gather}
where $\rho_{L^{*} \setminus \{0\}}$ is the anchor of the Lie algebroid associated to the Jacobi manifold $(\SYM,\Omega)$. In fact, since
$(\SYM,\Omega)$ is Poisson, the right hand side equals $ \rho_{L^{*} \setminus \{0\}} (\mu^*( T^*(\mathfrak{g}^* \setminus \{0\}) \oplus 0))$, where $J^1(\mathfrak{g}^* \setminus \{0\}) = T^*(\mathfrak{g}^* \setminus \{0\}) \oplus \RR$ (cf.\ Note~\ref{obs:trivial_case}). Unravelling the def\/initions, it can be checked that the map $F_{\pi} \circ \pi^* \colon O(1) \to J^1(\mathfrak{g}^* \setminus \{0\})$ is transverse to $0 \oplus \RR$. Therefore, equation~\eqref{eq:36} can be written as
\begin{gather*}
 \ker D \mu \subset \rho_{L^{*} \setminus \{0\}} (\mu^*\circ F_{\pi} \circ \pi^*(O(1)))= \rho_{L^{*} \setminus \{0\}}(F_{\operatorname{pr}} \circ
 \operatorname{pr}^* \circ F_{\phi} \circ \phi^*(O(1)) ),
\end{gather*}
where the last equality follows from equation~\eqref{eq:34}. Applying $D_{\alpha} \operatorname{pr}$ to both sides and observing that $D_{\alpha} \operatorname{pr} (\ker D_{\alpha} \mu) = \ker D_p \phi$, have that
 \begin{gather*}
 \ker D \phi \subset D \operatorname{pr} \circ \rho_{L^{*} \setminus \{0\}} (F_{\operatorname{pr}} \circ \operatorname{pr}^* \circ F_{\phi} \circ \phi^*(O(1)) ).
 \end{gather*}
Since $\operatorname{pr}$ is Jacobi with bundle component $F_{\operatorname{pr}}$, the right hand side equals
 $\rho(F_{\phi}(\phi^*(O(1))))$~-- cf.\ Note~\ref{obs:Jac_maps}~-- which equals $(\ker D_p \phi)^{\perp}$ by def\/inition. Thus property~\ref{item:14} holds.
\end{proof}

\section{Missing proofs from Sections \ref{sec:period-lattices-as} and \ref{sec:class-cir-over}}\label{apendixC}

\subsection*{Proofs of results of Section \ref{sec:period-lattices-as}}
The proof of Lemma \ref{cor:sections_period} and the associated preliminary results use in a crucial fashion the geometry of the Spencer operator associated to $\JM$. Fix a CIR $\phi\colon \CM\to (P,L,\{,\})$ and let $\theta$ denote the generalised contact form associated to $\CM$. Recall that for $\alpha\in\Gamma(\ker\rho)$, $\psi( \phi^*\alpha) \in \Gamma(\ker D\phi)$, where $\psi \colon \phi^* \ker \rho \to TM$ is def\/ined as in Lemma~\ref{lemma:action}. Therefore, the f\/low $\varphi^t_{\alpha}\colon M\to M$ of $\psi(\phi^*\alpha)$ preserves the f\/ibres of $\phi$, i.e., $\phi\circ \varphi^t_{\alpha}=\phi$ for all $t \in \RR$. Hence, $(\varphi^t_{\alpha})^* \phi^* L \cong \phi^* L$ canonically. This allows to def\/ine a `Lie
derivative' type operator on $\Omega^*(M;\phi^*L)$ by
\begin{gather}\label{eq:lie_derivative}
(\mathcal{L}_{\alpha}\omega)_m=\frac{d}{dt}\Big\rvert_{t=0}\big(\big(\varphi^t_{\alpha}\big)^*\omega\big)_m\end{gather}
for any $m\in M$ and $\omega \in \Omega^l(M,\phi^*L)$. It obeys rules which are analogous to those of the standard Lie derivative, for instance,
\begin{gather} \label{eq:8}
 \frac{d}{dt}\Big\rvert_{t=s}\big(\big(\varphi^t_{\alpha}\big)^*\omega\big)_m = ((\varphi^s_{\alpha})^*(\mathcal{L}_{\alpha} \omega))_m,
\end{gather}
as well as the following property proved in \cite[Lemma 3.8]{Maria}.

\begin{Lemma}\label{lemma:imt}
 For any $\omega\in\Omega^l(M;\phi^*L)$ and $X\in\X(M),$
 \begin{gather*}
 [i_X,\mathcal{L}_{\alpha}]\omega= i_{[X,\psi(\phi^*\alpha)]}\omega.
 \end{gather*}
\end{Lemma}

With the above preliminary results, the following can be proved.

\begin{Lemma}\label{lemma:fundamental}
 For any $\alpha \in \Gamma(\ker \rho)$,
 \begin{gather*}
 \big(\varphi^1_{(u,\eta)}\big)^*\theta - \theta = \phi^*\eta,
 \end{gather*}
where $\alpha = (u, \eta) \in \Gamma(L) \oplus \Omega^1(P;L)$ in the Spencer decomposition.
\end{Lemma}

\begin{proof}
Let $(u,\eta) \in \Gamma(\ker \rho)$ and $X\in\X(M)$. It suf\/f\/ices to show that
 \begin{gather}\label{eq:equiv}
 \left( \frac{\de}{\de t} \big(\big(\varphi^t_{(u,\eta)}\big)^* \theta\big)\right)(X) = (\phi^*\eta)(X);
 \end{gather}
for the lemma follows by integrating from $t=0$ to $t=1$. Computing the left hand side of equation~\eqref{eq:equiv} using equation~\eqref{eq:8} and Lemma~\ref{lemma:imt}, obtain that
 \begin{gather} \label{eq:9}
 \left(\frac{\de}{\de t} \big(\big(\varphi^t_{(u,\eta)}\big)^* \theta\big)\right)(X)=
 \mathcal{L}_{(u,\eta)}\big(\theta\big(D \varphi^t_{(u,\eta)}(X)\big)\big) + \theta\big[D \varphi^t_{(u,\eta)}(X),\psi(\phi^* u, \phi^*\eta)\big].
 \end{gather}
Set
 \begin{gather}\label{eq:11} D\varphi^t_{(u,\eta)}(X) = \sum\limits_i f_{t,i} A_{t,i} + X_t, \end{gather}
where $f_{t,i} \in C^{\infty}(M)$, $A_{t,i} := \psi(\phi^*u_{t,i},\phi^*\eta_{t,i}) \in \Gamma(\ker D\phi)$ with $(u_{t,i},\eta_{t,i})\in\Gamma(\ker\rho)$, and $X_t \in \Gamma(H)$ are all time-dependent, with $D\phi(X_t) = D\phi(X)$. This decomposition follows from property~\ref{item:IR2} and the fact that the action $\psi \colon \Gamma(\phi^* \ker \rho) \to \Gamma(\ker D\phi)$ is an isomorphism. Without loss of generality, assume $i=1$. On the one hand,
 \begin{gather}\label{eq:1}
 \theta\big(D\varphi^t_{(u,\eta)}(X)\big) = \theta ( f_{t}\psi(\phi^*u_{t},\phi^* \eta_{t})) =f_t\theta(\rho_M(\phi^*u_{t}, \phi^* \eta_{t}))=
 f_{t} \phi^*u_{t},
 \end{gather}
since the anchor $\rho_M$ of a contact manifold can be easily computed to be $ \rho_M(\phi^*u_{t}, \phi^* \eta_{t}) = R_{\phi^*u_{t}} + c^{\sharp}(\phi^* \eta_{t}\rvert_H)$. By def\/inition of $c^{\sharp}$ and of Reeb vector f\/ields, $c^{\sharp}( \phi^* \eta_{t}\rvert_H) \in \Gamma(H)$ and $\theta(R_{\phi^*u_{t}})=\phi^*u_{t}$ for all $t$ (cf.\ Example~\ref{exm:ctc_jacobi}). Using equation~\eqref{eq:1} and the fact that $\phi\circ\varphi^s_{(u,v)}=\phi$, have that $\phi^*u_{t,\varphi^s_{(u,v)}(m)}=u_{t,\phi(m)}$ for all~$s$; this in turn implies that
 \begin{gather} \label{eq:12}
 \mathcal{L}_{(u,\eta)}\big( \theta\big(D \varphi^t_{(u,\eta)}(X)\big)\big)_m =\mathcal{L}_{\psi(\phi^*u,\phi^*\eta)}(f_t)_mu_{t,m}.
 \end{gather}
On the other hand, using equation \eqref{eq:11}
 \begin{gather}
\theta \big(\big[D\varphi^t_{(u,\eta)}(X), \psi(\phi^*u,\phi^*\eta)\big]\big) = \theta([X_t,\psi(\phi^*u,\phi^*\eta)]) - \mathcal{L}_{\psi(\phi^*u,\phi^*\eta)}(f_{t})\theta(A_{t})\nonumber\\
 \hphantom{\theta \big(\big[D\varphi^t_{(u,\eta)}(X), \psi(\phi^*u,\phi^*\eta)\big]\big)}{}
 = \theta([X_t,\psi(\phi^*u,\phi^*\eta)]) - \mathcal{L}_{\psi(\phi^*u,\phi^*\eta)}(f_{t}) \phi^* u_{t}, \label{eq:15}
\end{gather}
where the f\/irst equality follows from the fact that $A_{t}, \psi(\phi^*u,\phi^*\eta) \in \Gamma(\ker D\phi)$ and this is a~bundle of {\it abelian} Lie algebras, and the second by def\/inition of $\theta(A_{t})$. Plugging equations~\eqref{eq:12} and~\eqref{eq:15} into equation~\eqref{eq:9}, obtain that
 \begin{gather} \label{eq:27}
 \left( \frac{\de}{\de t} \big(\varphi^t_{(u,\eta)}\big)^* \theta\right)(X)
 = \big(\varphi^t_{(u,\eta)}\big)^*(\theta([X_t,\psi(\phi^*u,\phi^*\eta)]));
 \end{gather}
the right hand side of equation \eqref{eq:27} can be computed to be
 \begin{gather*}
 \theta([X_t,c^{\sharp}( \phi^* \eta \rvert_H)]) = \phi^*\eta \rvert_H(X_t) = \phi^*(\eta(D\phi(X_t)))) = \phi^*(\eta(D\phi(X)),
 \end{gather*}
where the last equality uses that $D\phi(X_t) = D\phi(X)$. Therefore, equation \eqref{eq:27} becomes
 \begin{gather*}
 \left( \frac{\de}{\de t} \big(\varphi^t_{(u,\eta)}\big)^* \theta\right)(X) = \big(\varphi^t_{(u,\eta)}\big)^* \phi^*(\eta
 D\phi(X)) = \phi^*(\eta (D\phi X)) = (\phi^*\eta)(X), \end{gather*}
where the second equality follows from $\phi \circ \varphi^t_{(u,\eta)} = \phi$.
\end{proof}

Lemma \ref{cor:sections_period} is a obtained as a corollary of Lemma~\ref{lemma:fundamental}.

\begin{proof}[Proof of Lemma \ref{cor:sections_period}] Let $\phi\colon \CM\to(P,L,\{\cdot,\cdot\})$ be a CIR whose period net is~$\Sigma$. Use the Spencer decomposition to write a section as $(u,\eta) \in \Gamma(\Sigma)$ (possibly locally def\/ined). As $\Sigma$ is the isotropy of the action def\/ined by~\eqref{eq:6}, $\varphi^1_{(u,\eta)} = \mathrm{id}$. Applying Lemma~\ref{lemma:fundamental} to the section $(u,\eta)$, obtain that $\phi^*\eta= 0$, hence $\eta=0$. This means precisely that $(u,\eta)$ is holonomic, i.e., of the form $j^1 u$, as required.
\end{proof}

Having established properties of the above geometric properties of the period bundle of a CIR, we turn to the proof of the result relating transversal $\ZZ$-projective structures with transversal $\ZZ$-projective lattices.

\begin{proof}[Proof of Proposition \ref{thm:tilp_tilpa}] Throughout this proof, f\/ix a foliated manifold $(N,\mathcal{F})$, where $\mathcal{F}$ has codimension $l$. The proof f\/irst shows how to construct a transversal $\ZZ$-projective structure starting from a lattice and then the opposite construction is described. Checking that the constructions are inverse to one another is left as an exercise for the reader.

{\bf From lattices to structures.} Fix a line bundle $\pi\colon L \to N$ and a transversal $\ZZ$-projective lattice $\Sigma \subset J^1L$. Choose an open cover $\mathcal{U} = \{U_i\}$ of $N$ with the following properties
 \begin{itemize}\itemsep=0pt
 \item there exists a nowhere vanishing section $z_i \colon U_i \to \pi^{-1}(U_i)$;
 \item there exists a trivialisation of $\Sigma|_{U_i}$ with local frame $\alpha^i_1,\ldots, \alpha^i_{l+1}\in \Gamma(\Sigma|_{U_i})$.
 \end{itemize}
Fix the above trivialisations and let $c_{ij}\colon U_{ij} \to \mathrm{GL}(1;\RR)$ and $A_{ij}\colon U_{ij} \to \mathrm{GL}(l+1;\ZZ)$ denote the corresponding transition functions. In fact, without loss of generality, it may be assumed that, for all~$i$,~$j$, $c_{ij}$ takes values in $\mathrm{GL}(1;\ZZ) = \mathrm{O}(1;\RR)$. Property~\ref{item:9} implies that for all $i$ and each $r = 1,\ldots, l+1$, there exist smooth functions
 \begin{gather*}
 g^i_r\colon \ U_i\to \RR\qquad \text{with}\quad \alpha^i_r = j^1\big(g_r^iz_i\big).
 \end{gather*}
The transversal $\ZZ$-projective structure is constructed using these functions~$g^i_r$; before proceeding to the construction, two preparatory claims are proved.

\begin{Claim}\label{claim:trans}
The map $\bar\chi_i\colon U_i\to \RR^{l+1}$, $\bar\chi_i:=(g^i_1,\ldots,g^i_{l+1})$, takes values in $\RR^{l+1}\setminus 0$ and is transversal to the Euler vector field $E=\sum\limits_rx_r\frac{\partial}{\partial x_r}$.
 \end{Claim}

\begin{proof}[Proof of Claim \ref{claim:trans}] As $j^1(g^i_1z_i),\ldots,j^1(g^i_{l+1}z_i)$ are a frame of $\Sigma|_{U_i}$ and $\Sigma^{\RR}\to L$ is onto by property~\ref{item:11}, it follows that $\bar\chi_i \neq 0$. Secondly, to prove that $\bar\chi_i$ is transversal to $E$, i.e., $\Im D\bar\chi_i + \RR \langle E \rangle=T\RR^{l+1}$, it suf\/f\/ices to show that $\bar\chi_i^*\colon \Ann(E)\to \nu^*$ is injective, where $\nu^*$ is the conormal bundle to $\mathcal{F}$. For, if this is the case, then $D\bar{\chi}_i|_{\nu} \colon \nu \to \Ann(E)^*$ is onto and the claim follows. First, it is shown that $\bar{\chi}^*_i(\Ann(E)) \subset \nu^*$; to this end, observe that a transversal $\ZZ$-projective lattice induces a f\/lat $T\mathcal{F}$-connection on $L \to N$, which is def\/ined as follows. Let $v \in \Gamma(L)$ and f\/ix $p \in N$; by property~\ref{item:11}, there exists an open neighbourhood $V \subset N$, and a~section $j^1u \in \Gamma(\Sigma|_V)$ with with $u(q) \neq 0$ for all $q \in U$. Then, locally, there exists $f \in C^{\infty}(U)$ with $v =f u$. For $X \in \Gamma(T\mathcal{F}|_U)$, def\/ine
 \begin{gather*}
 \nabla_X v := \mathcal{L}_X(f) u.
 \end{gather*}
It can be checked that this is well-def\/ined, that is, indeed, a f\/lat $T\mathcal{F}$-connection on $L \to N$, and that if $j^1u \in \Gamma(\Sigma)$, then $\nabla u =0$. Using this connection, have that there exists a closed, foliated 1-form $\beta_i$ such that, for all~$r$, $dg^i_r|_{T\mathcal{F}}=g^i_r\beta_i|_{T\mathcal{F}}$. Hence
 \begin{gather}\label{eq:euler}
 D\bar\chi_i|_{T\mathcal{F}}=\beta_i \otimes E,
 \end{gather}
which implies that $\bar\chi_i^*(\Ann(E))\subset \nu^*$. To show that $\bar\chi_i^*\colon \Ann(E)\to \nu^*$ is injective, it suf\/f\/ices to show that it is surjective (which is, of course, equivalent by dimension counting). Let $\gamma \in \Gamma(\nu^*|_{U_i})$, by property~\ref{item:11} there exist
 $f^1,\ldots,f^{l+1}\in C^\infty(U_i)$ with
 \begin{gather*}
 \gamma=\sum_rf^rj^1g^i_r\in \nu^*\subset \Sigma^\RR|_{U_i}\subset J^1L|_{U_i}.
 \end{gather*}
Using the Spencer decomposition, it follows that $\sum_rf^rg^i_r=0$ and $\gamma=\sum_rf^rdg^i_r$. Therefore, $\bar\chi_i^*(\gamma')=\gamma$ for
 \begin{gather*}
 \gamma'=\sum\limits_rf^r(p)\de x_r\in \Gamma\big(\Ann(E)|_{\bar\chi_i(U_i)}\big),
 \end{gather*}
 and hence surjectivity follows.
 \end{proof}

 \begin{Claim}\label{claim:commutative}
 The diagram
 \begin{gather*}\xymatrix{
 & U_{ij} \ar[dl]_{\bar\chi_i} \ar[dr]^{\bar\chi_j} & \\
 \RR^{l+1}\setminus 0\ar[rr]^{c_{ij}A_{ij}} & & \RR^{l+1}\setminus 0
 }
 \end{gather*}
 commutes.
 \end{Claim}
 \begin{proof}[Proof of Claim \ref{claim:commutative}] This follows by def\/inition of the maps $\bar{\chi}_i$.
 \end{proof}

Using Claims \ref{claim:trans} and \ref{claim:commutative}, a~transversal $\ZZ$-projective structure can be def\/ined as follows. Set $\chi_i\colon U_i\to \RR\mathrm{P}^l$ to be the composite
 \begin{gather*}
 \chi_i\colon \ U_i\overset{\bar\chi_i}{\to} \RR^{l+1}\setminus 0\overset{q}{\to}\RR\mathrm{P}^l,
 \end{gather*}
where $q\colon \RR^{l+1}\setminus 0\to \RR\mathrm{P}^l$ is the quotient map. By Claim~\ref{claim:trans}, $\chi_i$ is a submersion, and Claim~\ref{claim:commutative} implies that $\chi_j = [A_{ij}] \circ \chi_i $ on $U_{ij}$. Moreover, since $Dq(E)=0$, equation~\eqref{eq:euler} implies that $D\chi_i|_{T\mathcal{F}}=0$, hence the leaves of $\mathcal{F}$ are tangent to the f\/ibres of $\chi_i$ on the one hand, and on the other, dimension counting shows that for $p\in U_i$, $T\mathcal{F}_p=\ker D_p\chi_i$.

{\bf From structures to lattices.} The idea is to pull-back the $\ZZ$-projective lattice on $\RR\mathrm{P}^l$ constructed in Example~\ref{ex:projective_plane} using the transversal $\ZZ$-projective structure. Fix a transversal $\ZZ$-projective structure~$\mathcal{A} = \{(U_i,\chi_i)\}$ with cocycle $A_{ij} \colon U_{ij} \to \mathrm{GL}(l+1;\ZZ)$. On each $i$, set $L_i:= \chi^*_i(O(1))$ and $\Sigma_i = \chi^*_i \Sigma^l$, where $\Sigma^l \subset J^1(O(1))$ is as in Example~\ref{ex:projective_plane}. The line bundle $\pi \colon O(1) \to \RR\mathrm{P}^l$ is $\mathrm{GL}(l+1;\ZZ)$-linearisable in the sense of \cite{kklv}, i.e., there exists an action of $\mathrm{GL}(l+1;\ZZ)$ on $O(1)$ which is linear on the f\/ibres of $\pi$ and makes $\pi$ $\mathrm{GL}(l+1;\ZZ)$-equivariant. The induced action is the standard $\mathrm{GL}(l+1;\ZZ)$-action on the vector space $\RR \langle x_1,\ldots, x_{l+1} \rangle$ and, thus, it preserves $\Sigma_l$. Suppose that $U_{ij} \neq \varnothing$; then
 \begin{gather*}
 L_j|_{U_{ij}} = \chi^*_j (O(1)) = \big([A_{ij}] \circ \chi_i\big)^* (O(1)) = \chi_i^* (O(1)) = L_i|_{U_{ij}},
 \end{gather*}
where the third equality uses the fact that $[A_{ij}]^*(O(1)) = O(1)$ as $O(1) \to \RR\mathrm{P}^l$ is $\mathrm{GL}(l+1;\ZZ)$-linearised. Therefore, the above construction yields a line bundle $L \to N$. Similarly, it can be shown that $\Sigma_j|_{U_{ij}} = \Sigma_i|_{U_{ij}}$, thus obtaining a $\ZZ^{l+1}$-bundle $\Sigma \to N$ whose total space is an embedded submanifold of $J^1L$. Therefore property~\ref{item:8} holds. Property~\ref{item:11} holds since the submersions $\chi_i$ locally def\/ine $\mathcal{F}$ (thus showing that $\nu^* \otimes L \hookrightarrow \Sigma^{\RR}$) and because $\left(\Sigma^l\right)^{\RR}$ surjects onto $O(1)$, thus showing that $\Sigma^{\RR} \to L$ is also onto. Local sections of $\Sigma$ are holonomic by construction, thus showing that property \ref{item:9} holds.
\end{proof}

\subsection*{Proofs of results from Section \ref{sec:class-cir-over}}
\begin{proof}[Proof of Theorem \ref{thm:aa}] First it is shown that the dif\/ference $\Psi_\sigma^*\theta-\theta_0\in\Omega^1(\ker\rho/\Sigma|_U;\pi^*L)$
is basic, i.e.,
\begin{gather}\label{eqn:contact form}
\Psi_\sigma^*\theta-\theta_0=\pi^*\beta
\end{gather}
for some 1-form $\beta\in\Omega^1(U;L)$. This is the case if $i_Z(\Psi_\sigma^*\theta-\theta_0)=0$ and $\mathcal{L}_Z(\Psi_\sigma^*\theta-\theta_0)=0$ for any $Z$ tangent to $\ker(\pi\colon \ker\rho/\Sigma|_U\to U) =: T^{\pi}(\ker \rho/\Sigma)$ (where $\mathcal{L}_Z\omega$ is
def\/ined as in equation~\eqref{eq:lie_derivative}). In order to prove this, identify (canonically) $T^\pi_{z}(\ker_{p}/\Sigma_p)$ with $\ker\rho_{\pi(z)}$ via the isomorphism $T^\pi(\ker\rho/\Sigma)\simeq \pi^*\ker\rho$. As $\Psi_\sigma$ comes from the action of equation~\eqref{eq:6}, a~straightforward computation shows that for $(u',\eta'),(u,\eta)\in \Gamma(\ker\rho)$,
\begin{gather*}
 \Psi_\sigma((u',\eta')+t(u,\eta))=\varphi^t_{(u,\eta)}(\Psi_\sigma(u',\eta'))
\end{gather*}
for any $t\in\RR$, where $\varphi^t_{(u,\eta)} \colon M\to M$ is the f\/low of $\psi(\phi^*(u,\eta))$ and the Spencer decomposition has been used. Equivalently,
\begin{gather}\label{eq:conmuta}\Psi_\sigma\circ \varphi^t_{(u,v)}=\varphi^t_{(u,v)}\circ\Psi_\sigma,\end{gather}
where on the left hand side $\varphi^t_{(u,v)}\colon \ker\rho/\Sigma\to\ker\rho/\Sigma$ stands for the f\/low $z\mapsto z+t(u,\eta)(\pi(z))$ of the vertical vector f\/ield $\pi^*(u,\eta)$. Dif\/ferentiating equation~\eqref{eq:conmuta}, obtain that
\begin{gather*}
 D_{z}\Psi_\sigma(\pi^*(u,\eta))=\psi(\phi^*u,\phi^*\eta)_{\Psi_\sigma(z)}.
\end{gather*}
Let $Z=\pi^*(u,\eta)$; by def\/inition of $\psi$ and the fact that for contact manifold $\rho_M(\alpha,\gamma)=R_\alpha+c^\sharp(\gamma)$, for $(\alpha,\gamma)\in\Gamma(L_M)\oplus\Omega^1(M;L_M),$ for $R_\alpha$ the Reeb vector f\/ield of $\alpha$ (cf.\ Example~\ref{exm:ctc_jacobi}), and~$c$ the curvature map of equation~\eqref{curvature}, the above yields
\begin{gather*}
 i_Z\Psi^*_\sigma\theta=\theta(\psi(\phi^*u,\phi^*\eta))=\bar\theta(R_{F\phi^*u}+c(F\phi^*\eta|_H))=u.
\end{gather*}
On the other hand, $\theta_0$ restricted to $\pi^*\ker\rho$ is equal to the projection $\pr\colon \ker\rho\to L$ which implies that $i_Z\theta_0=\theta_{0}(u,\eta)=u$. Thus $i_Z(\Psi_\sigma^*\theta-\theta_0)=0$ follows. To compute the Lie derivative along $Z$, observe that
\begin{gather*}
 \mathcal{L}_Z(\Psi^*_\sigma\theta )= \frac{d}{dt}\big(\varphi^t_{(u,\eta)}\big)^*(\Psi^*_\sigma\theta)|_{t=0} =\frac{d}{dt}\Psi^*_\sigma\circ
\big(\varphi^t_{(u,\eta)}\big)^*\theta|_{t=0},
\end{gather*}
where the second equality uses equation \eqref{eq:conmuta}. By equation~\eqref{eq:equiv}, have that
\begin{gather*}
 \frac{d}{dt}\big(\varphi^t_{(u,\eta)}\big)^*\theta|_{t=0}=\phi^*\eta,
\end{gather*}
which implies that $\mathcal{L}_Z\Psi^*_\sigma\theta=\pi^*\eta$. On the other hand,
\begin{gather*} \mathcal{L}_Z\theta_0=\frac{d}{dt}\big(\varphi^t_{(u,\eta)}\big)^*\theta_0|_{t=0}=\frac{d}{dt}\big({\rm id}^*+t\pi^*\circ(u,\eta)^*\big)\theta_0|_{t=0}=\frac{d}{dt}(\theta_0+t\pi^*\eta)|_{t=0}=\pi^*\eta,
\end{gather*}
where the third equality uses that for a section $(u,\eta)\colon P\to J^1L$, $(u,\eta)^*\theta_{\mathrm{can}}=\eta$. With this, $\mathcal{L}_Z(\Psi_\sigma^*\theta-\theta_0)=0$. As vector f\/ields of the form $\pi^*(u,v)$ generate $\pi^*\ker\rho\simeq T^\pi(\ker\rho/\Sigma)$ as a~$C^\infty(\ker\rho/\Sigma)$-module, this implies that equation~\eqref{eqn:contact form} holds.

To show that $\beta=\sigma^*\theta$, consider the section $z\colon P\to \ker\rho/\Sigma, p\mapsto [\Sigma_p]=0$ of $\pi\colon \ker/\Sigma\to P$. Then $z^*\operatorname{pr}^*\beta=\beta$, $z^*\Psi_\sigma^*\theta=\sigma^*\theta$ as $\Psi_\sigma\circ z=\sigma$, and $z^*\theta_0=0$ as for any $s\in\Gamma(\Sigma)$, $s^*\theta_{\mathrm{can}}=0$ (see Note~\ref{obs:contact form}). Therefore, $\beta=\sigma^*\theta$.

It remains to show that $\de_\mathcal{F}(\sigma^*\theta)=\omega_{\mathcal{F}}$. As $\Gamma(\mathcal{F})$ is generated by elements of the form $\rho(j^1v),v\in\Gamma(L)$, it suf\/f\/ices to check the statements for vector f\/ields of this form. By def\/inition,
\begin{gather}
 \de_{\mathcal{F}}(\sigma^*\theta)\big(\rho\big(j^1u\big),\rho\big(j^1v\big)\big)
 = -\sigma^*\theta\big(\big[\rho\big(j^1u\big), \rho\big(j^1v\big)\big]\big)\nonumber\\
 \hphantom{\de_{\mathcal{F}}(\sigma^*\theta)\big(\rho\big(j^1u\big),\rho\big(j^1v\big)\big)=}{} +
 {\nabla}_{\rho(j^1u)}\big(\sigma^*\theta\big(\rho\big(j^1v\big)\big)\big) - {\nabla}_{\rho(j^1v)}\big(\sigma^*\theta\big(\rho\big(j^1u\big)\big)\big). \label{eq:30}
\end{gather}
First, consider $ {\nabla}_{\rho(j^1u)}(\sigma^*\theta(\rho(j^1v))) ={\nabla}_{\rho(j^1u)}(\theta(D\sigma(\rho(j^1v))))$; since $D\phi \circ D\sigma = {\rm id}$ and $\phi$ is a~Jacobi map, the diagram of Note~\ref{obs:Jac_maps} commutes. Thus
\begin{gather*}
 D\sigma\big(\rho\big(j^1v\big)\big) = \rho_M\big(j^1\big(\phi^*v\big)\big) + X_v,
\end{gather*}
where $X_v \in \Gamma(\ker D\phi|_{\sigma(U)})$ (and similarly for $u$).

As $\phi$ satisf\/ies property~\ref{item:IR3}, $\ker D\phi =\rho_M (\phi^* \ker \rho)$, hence there exist local sections $(w,\eta), (z,\zeta)$ $\in \Gamma(\ker \rho)$ with
\begin{gather*}
 X_u = \rho_M(\phi^*w,\phi^*\eta) \qquad \text{and} \qquad X_v = \rho_M(\phi^*z,\phi^*\zeta).
\end{gather*}
Computing $\rho_M$ for a contact manifold,
\begin{gather*}
 \theta\big(D\sigma\big(\rho\big(j^1v\big)\big)\big) =\theta(R_{\phi^*v}+R_{\phi^*z}+c^\sharp(\phi^*\zeta|_H)) = \phi^*(v + z),
\end{gather*}
and, similarly, $\theta(D\sigma(\rho(j^1u))) = \phi^*(u+w)$. Therefore, by def\/inition of ${\nabla}$,
\begin{gather} \label{eq:31}
 {\nabla}_{\rho(j^1u)}\big(\sigma^*\theta\big(\rho\big(j^1v\big)\big)\big) - {\nabla}_{\rho(j^1v)}\big(\sigma^*\theta\big(\rho\big(j^1u\big)\big)\big) = \{u,v+z\} - \{v, u + w\}.
\end{gather}
On the other hand,
\begin{gather}
(\sigma^*\theta)\big(\big[\rho\big(j^1u\big), \rho\big(j^1v\big)\big]\big) =
 \sigma^*\big(\theta\big(D\sigma \big(\big[\rho\big(j^1u\big), \rho\big(j^1v\big)\big]\big) \big)\big) \nonumber\\
\hphantom{(\sigma^*\theta)\big(\big[\rho\big(j^1u\big), \rho\big(j^1v\big)\big]\big)}{}
 = \sigma^*\big(\theta\big(\big(\big[D\sigma \big(\rho\big(j^1u\big)\big), D\sigma\big(\rho\big(j^1v\big)\big)\big]\big)\big)\big) \nonumber\\
\hphantom{(\sigma^*\theta)\big(\big[\rho\big(j^1u\big), \rho\big(j^1v\big)\big]\big)}{}
 = \sigma^*\big(\theta\big(\big[R_{\phi^*(u+w)} +
 c^{\sharp}(\phi^*\eta|_H), R_{\phi^*(v+z)} + c^{\sharp}(\phi^*\zeta|_H)\big]\big)\big) \nonumber\\
\hphantom{(\sigma^*\theta)\big(\big[\rho\big(j^1u\big), \rho\big(j^1v\big)\big]\big)}{}
 = \sigma^*\big(\theta\big(R_{\{\phi^*(u+w), \phi^*(v+z)\}_M} + \big[c^{\sharp}(\phi^*\eta|_H), c^{\sharp}(\phi^*\zeta|_H)\big]\big)\big) \nonumber\\
\hphantom{(\sigma^*\theta)\big(\big[\rho\big(j^1u\big), \rho\big(j^1v\big)\big]\big)}{}
 = \sigma^*(\phi^*(\{u+w,v+z\}_M) - \phi^*\eta(\rho_M( \phi^*\zeta)))\label{eq:32}\\
\hphantom{(\sigma^*\theta)\big(\big[\rho\big(j^1u\big), \rho\big(j^1v\big)\big]\big)}{}
 = \{u+w,v+z\} - \eta(D\phi (\rho_M(
 \phi^*\zeta))) = \{u+w,v+z\} - \eta(\rho(\zeta)), \nonumber
\end{gather}
where the fourth equality follows from the def\/ining property of Reeb vector f\/ields (cf.\ Example~\ref{exm:ctc_jacobi}), the f\/ifth by def\/inition of Reeb vector f\/ields and the curvature map of equation~\eqref{curvature}, and the last from the fact that $\phi$ is a Jacobi map. Observe that since $(w,\eta),(z,\zeta) \in \Gamma(\ker \rho)$,
\begin{gather*}
 \{w,z\} - \eta (\rho(\zeta)) = \{w,z\} + \eta\big(\rho\big(j^1z\big)\big) = \omega (\rho(w,\eta),\rho(z,\zeta)) = 0,
\end{gather*}
where the general def\/inition of $\omega_{\mathcal{F}}$ is used (cf.\ Note~\ref{obs:prop_jac_even}). Thus equation \eqref{eq:32} yields that
\begin{gather*}
(\sigma^*\theta)\big(\big[\rho\big(j^1u\big), \rho\big(j^1v\big)\big]\big) = \{u,v\} + \{u,z\} + \{w,v\}.
\end{gather*}
Using \looseness=-1 this identity together with equation~\eqref{eq:31} in equation~\eqref{eq:30} yields the required result.
\end{proof}

\subsection*{Acknowledgements}
We would like to thank two anonymous referees for the suggestions that helped improve signif\/icantly the content and its presentation. Furthermore, we would like to thank Camilo Arias Abad for interesting conversations. M.A.S.\ would like to thank IMPA, CRM and MPIM Bonn for hospitality at various stages of the project. M.A.S.\ was partly supported by the DevMath programme of the Centre de Recerca Matem\`atica and by the Max Planck Institute for Mathematics in Bonn. D.S.\ was partly supported by ERC starting grant 279729, by the NWO Veni grant 639.031.345 and by CNPq.

\LastPageEnding

\end{document}